\documentclass[12pt]{jhrs}
\usepackage{amsmath}
\usepackage{graphics}
\usepackage[all]{xy}
\xyoption{v2}
\xyoption{knot}

\let\url=\undefined
\usepackage
[bookmarks,colorlinks,linkcolor=blue,citecolor=blue,urlcolor=blue]
{hyperref}

\title{Operads in iterated monoidal categories}
\author {Stefan Forcey, Jacob Siehler and E. Seth Sowers}

\thanks{Thanks to {\Xy-pic} for the diagrams. }
\email{sforcey@tnstate.edu}
\address{Department of Physics and Mathematics\\
Tennessee State University\\
Nashville, TN 37209 \\
USA\\
}


\swapnumbers
\theoremstyle{plain}
\newtheorem{theorem}{Theorem}[section]
\newtheorem{lemma}[theorem]{Lemma}
\newtheorem{corollary}[theorem]{Corollary}

\theoremstyle{definition}
\newtheorem{definition}[theorem]{Definition}
\newtheorem{example}[theorem]{Example}

\theoremstyle{remark}
\newtheorem{remark}[theorem]{Remark}

\def\cal#1{\mathcal{#1}}
\newcommand{\bcal}[1]{\mbox{\boldmath${\cal {#1}}$}}

\def\nat{\mathbb{N}}
\def\th{^\text{th}}
\def\dash{\text{---}}
\DeclareMathOperator{\concat}{Concat}
\DeclareMathOperator{\seq}{Seq}
\DeclareMathOperator{\modseq}{ModSeq}
\DeclareMathOperator{\col}{Col}
\DeclareMathOperator{\mon}{Mon}
\DeclareMathOperator{\oper}{Oper}

\begin{document}

\begin{abstract}
The structure of a $k$-fold monoidal category as introduced by
Balteanu, Fiedor\-owicz, Schw\"anzl and Vogt in \cite{Balt} can be
seen as a weaker structure than a symmetric or even braided monoidal
category.  In this paper we show that it is still sufficient to
permit a good definition of ($n$-fold) operads in a $k$-fold monoidal
category which generalizes the definition of operads in a braided
category. Furthermore, the inheritance of structure by the category
of operads is actually an inheritance of iterated monoidal structure,
decremented by at least two iterations.  We prove that the category of $n$-fold operads in a
$k$-fold monoidal category is itself a $(k-n)$-fold monoidal, strict
$2$-category, and  show
that $n$-fold operads are automatically $(n-1)$-fold operads.
We also introduce a family of simple examples of
$k$-fold monoidal categories and classify operads in the example
categories.
\end{abstract}
\received{Month Day, Year}   
\revised{Month Day, Year}    
\published{Month Day, Year}  
\submitted{J. Stasheff}  

\volumeyear{2006} 
\volumenumber{1}  
\issuenumber{1}   

\startpage{1}     
\maketitle
\tableofcontents

\section{Introduction}
In this introductory section we will give a brief, non-chronological overview of the relationship
between operads, higher category theory, and topology. This will serve to motivate the study of
iterated monoidal categories and their operads that comprises the remaining sections.  In the
second section, in order to be self contained, we repeat the definition of the iterated monoidal
categories first set down in \cite{Balt}.  In the fourth section we seek to fill a gap in the
literature which currently contains few good examples of that definition.  Thus our first new
contribution consists of a series of simple and very geometric iterated monoidal categories based
on totally ordered monoids. By simple we mean that axioms are largely fulfilled due to
relationships between max, plus, concatenation, sorting and lexicographic ordering as well as the
fact that all diagrams commute since the underlying directed graph of the category is merely the
total order. The most interesting examples of $n$-fold monoidal categories are those whose objects
can be represented by Ferrer or Young diagrams  (the underlying shapes of Young tableaux.) These
exhibit products with the geometrical interpretation ``combining stacks of boxes.'' Managers of
warehouses or quarries perhaps may already be well acquainted with the three dimensional version of
the main example of iterated monoidal categories we introduce here. Imagine that floor space in the
quarry or warehouse is at a premium and that therefore you are charged with combining several
stacks of crates or stone blocks by restacking some together vertically and shifting others
together horizontally. It turns out that the best result in terms of gained floor space is always
to be achieved most efficiently by doing the restacking and shifting in  a very particular
order--horizontally first, then vertically.

The main new contribution is the theory of operads within, or enriched in, iterated monoidal
categories.  This theory is based upon the fact that the natural setting of operads turns out to be
in a category with lax interchange between multiple operations, as opposed to the full strength of
a braiding or symmetry as is classically assumed.  Batanin's definition of $n$-operad also relies
on this insight \cite{bat}. In that paper he notes that an iterated monoidal category ${\cal V}$
would be an example of a globular monoidal category with a single object, and a single arrow in
each dimension up to $n$, in which last dimension the arrows would actually be the objects of
${\cal V}.$ Of course the invertibility of the interchange would also have to be dropped from his
definition. In that case the $n$-fold operads defined here would correspond to Batanin's
$n$-operads. The advantages of seeing them in a single categorical dimension are in the way that
doing so generalizes the fact that operads in a symmetric monoidal category inherit its symmetric
structure.  We  investigate the somewhat flexible structure of the iterated monoidal 2-category
that $n$-fold operads comprise. Flexibility arises from the difference between $n$ and $k$, where
one is investigating $n$-fold operads in a $k$-fold monoidal category ${\cal V}$, where $n < k-1$.
It turns out that choosing $n$ much smaller than $k$ allows multiple interchanging products to be
defined on the category of operads, whereas choosing $n$ closer to $k$ allows the operads to take
on multiple operad structures at once with respect to the products in ${\cal V}.$ Examples of
combinatorial operads living in the previously introduced combinatorially defined categories are
utilized to demonstrate the sharpness of several of the resulting theorems, i.e. to provide
counterexamples. The examples start to take on a life of their own, however, as theorems and open
questions about the classification of operads   in combinatorial $n$-fold monoidal categories
arise.  The definition of operad in the categories with morphisms given by ordering leads to
descriptions of interesting kinds of growth.  We give a complete description of the simple example
of 2-fold operads in the natural numbers.  We then give the elementary results for operads in the
category of Young diagrams. In the basic examples linear and logarithmic growth characterize
respective dimensions in a single sequence of Young diagrams.  These phenomena hint towards a
theory of operadic growth. Full investigation and further classification must await a sequel to
this paper. Applications might be found in scientific fields such as the theory of small world
networks, where the diameter of a network is the logarithm of the number of nodes.

First, however, we look at  some of the  history and philosophy of
the two major players here, operads and iterated monoidal categories.
Operads in a category of topological spaces are the crystallization
of several approaches to the recognition problem for iterated loop
spaces. Beginning with Stasheff's associahedra and Boardman and
Vogt's little $n$-cubes, and continuing with more general $A_{\infty}$,
$E_n$ and $E_{\infty}$ operads described by May and others, that
problem has largely been solved \cite{Sta}, \cite{BV1}, \cite{May}.
Loop spaces are characterized by admitting an operad action of the
appropriate kind. More lately Batanin's approach to higher categories
through internal and higher operads  promises to shed further light
on the combinatorics of $E_n$ spaces \cite{bat2}, \cite{bat3}.

Recently there has also been growing interest in the application
of higher dimensional structured categories to the characterization
of loop spaces. The program being advanced by many categorical
homotopy theorists seeks to model the coherence laws governing
homotopy types with the coherence axioms of structured  $n$-categories.
By modeling we mean a connection that will  be in the form of a
functorial equivalence between categories of special categories and
categories of special spaces.  The largest challenges currently are
to find the most natural and efficient definition of (weak)
$n$-category, and to determine the nature of the functor from
categories to spaces.  The latter will almost certainly be analogous
to the nerve functor on 1-categories, which preserves homotopy
equivalence.  In \cite{StAlg} Street defines the nerve of a strict
$n$-category.  Recently Duskin in \cite{Dusk} has worked out the
description of the nerve of a bicategory.  A second part of the
latter paper promises the full description  of the functor including
how it takes morphisms of bicategories to continuous maps.

One major recent advance is the discovery of Balteanu, Fiedorowicz,
Schw\"anzl and Vogt in \cite{Balt} that the nerve functor
on categories gives a direct connection between iterated monoidal
categories and iterated loop spaces.  Stasheff \cite{Sta} and Mac
Lane \cite{Mac} showed that monoidal categories are precisely
analogous to 1-fold loop spaces. There is a similar connection
between symmetric monoidal categories and infinite loop spaces. The
first step in filling in the gap between 1 and infinity was made
in \cite{ZF} where it is shown that the group completion of the
nerve of a braided monoidal category is a 2-fold loop space.   In
\cite{Balt} the authors finish this process by, in their words,
``pursuing an analogy to the tautology that an $n$-fold loop space
is a loop space in the category of $(n-1)$-fold loop spaces.'' The
first thing they focus on is the fact that a braided category is a
special case of a carefully defined 2-fold monoidal category. Based
on their observation of the  correspondence between loop spaces and
monoidal categories, they iteratively define the notion of $n$-fold
monoidal category as a monoid in the category of $(n-1)$-fold
monoidal categories.  In \cite{Balt} a symmetric category is seen
as a category that is $n$-fold monoidal for all $n$.  The main
result in that paper is that the group completion of the nerve of
an $n$-fold monoidal category is an $n$-fold loop space.  It is
still open whether this is a complete characterization, that is,
whether every $n$-fold loop space arises as the nerve of an $n$-fold
monoidal category. Much progress towards the answer to this question
was made by the original authors in their sequel paper, but the
desired result was later shown to remain unproven. One of the future
goals of the program begun here is to use weakenings or deformations
of the examples of $n$-fold monoidal categories introduced here to
model specific loop spaces in a direct way.

The connection between the $n$-fold monoidal categories of Fiedorowicz
and the theory of higher categories is through the  periodic table
as laid out in \cite{Baez1}.  Here Baez organizes the $k$-tuply
monoidal $n$-categories, by which terminology he refers to
$(n+k)$-categories  that are trivial below dimension $k.$ The
triviality of lower cells allows the higher ones to compose freely,
and thus these special cases of $(n+k)$-categories are viewed as
$n$-categories with $k$ multiplications.  Of course a $k$-tuply
monoidal $n$-category is a special $k$-fold monoidal $n$-category.
The specialization results from the definition(s) of $n$-category,
all of which seem to include the axiom that the interchange
transformation between two ways of composing four higher morphisms
along two different lower dimensions is required to be an isomorphism.
As will be mentioned in the next section the  property of having
iterated loop space nerves held by the $k$-fold monoidal categories
relies on interchange transformations that are not isomorphisms.
If those transformations are indeed isomorphisms then the $k$-fold
monoidal 1-categories do reduce to the braided and symmetric
1-categories of the periodic table. Whether this continues for
higher dimensions, yielding for example the sylleptic monoidal
2-categories of the periodic table as 3-fold monoidal 2-categories
with interchange isomorphisms, is yet to be determined.

A further refinement of higher categories is to require all morphisms
to have inverses. These special cases are referred to as $n$-groupoids,
and since their nerves are simpler to describe it has long been
suggested that they model homotopy $n$-types through a construction
of a fundamental $n$-groupoid. This has in fact been shown to hold
in Tamsamani's definition of weak $n$-category \cite{tam}, and in
a recent paper by Cisinski to hold in the definition of Batanin as
found in \cite{bat}.  A homotopy $n$-type is a topological space
$X$ for which $\pi_k(X)$ is trivial for all $k>n.$
 It has been
suggested that a key requirement for any useful definition of
$n$-category is that a $k$-tuply monoidal $n$-groupoid be associated
functorially (by a nerve) to a topological space which is a  homotopy
$n$-type and a $k$-fold loop space \cite{Baez1}.  The loop degree
will be precise for $k<n+1,$ but for $k>n$ the associated homotopy
$n$-type will be an infinite loop space.  This last statement is a
consequence of the stabilization hypothesis , which states that
there should be a left adjoint to forgetting monoidal structure
that is an equivalence of $(n+k+2)$-categories between $k$-tuply
monoidal $n$-categories and $(k+1)$-tuply monoidal $n$-categories
for $k>n+1.$ This hypothesis has been shown by Simpson to hold in
the case of Tamsamani's definition \cite{sim}. For the case of $n=1$
if the interchange transformations are isomorphic then a $k$-fold
monoidal 1-category is equivalent to a symmetric category for $k>2.$
With these facts in mind it is possible that if we wish to precisely
model homotopy $n$-type $k$-fold loop spaces for $k>n$ then we may
need to consider $k$-fold as well as $k$-tuply monoidal $n$-categories.
This paper is part of an embryonic  effort in that direction.

Since a loop space can be efficiently described as an operad algebra,
it is not surprising that there are several existing definitions
of $n$-category that utilize operad actions.  These definitions
fall into two main classes: those that define an $n$-category as
an algebra of a higher order operad, and those that achieve an
inductive definition using classical operads in symmetric monoidal
categories to parameterize iterated enrichment. The first class of
definitions is typified by Batanin and Leinster \cite{bat},\cite{lst}.

The former author defines monoidal globular categories in which
interchange transformations are isomorphisms and which thus resemble
free strict $n$-categories.  Globular operads live in these, and
take all sorts of pasting diagrams as input types, as opposed to
just a string of objects as in the case of classical operads.  The
binary composition in an $n$-category derives from the action of a
certain one  of these globular operads.  Leinster expands this
concept to describe $n$-categories with unbiased composition of any
number of cells.  The second class of definitions is typified by
the works of Trimble and May \cite{trimble}, \cite{may2}.

The former parameterizes iterated enrichment with a series of operads
in $(n-1)$-Cat achieved by taking the fundamental $(n-1)$-groupoid
of the $k\th$ component of the topological path composition operad
$E.$ The latter begins with an $A_{\infty}$ operad in a symmetric
monoidal category ${\cal V}$ and requires his enriched categories
to be tensored over ${\cal V}$ so that the iterated enrichment
always refers to the same original operad.

Iterated enrichment over $n$-fold categories is described in
\cite{forcey1} and \cite{forcey2}.  It seems worthwhile to define
$n$-fold operads in $n$-fold monoidal categories in a way that is
consistent with the spirit of Batanin's globular operads.  Their
potential value may include using them to weaken enrichment over
$n$-fold monoidal categories in a way that is in the spirit of May
and Trimble.  This program carries with it the promise of
characterizing $k$-fold loop spaces with homotopy $n$-type for all
$n,k$ by describing  the categories with exactly those spaces as
nerves.  As a candidate for the type of category with such a nerve
we suggest a weak $n$-category with $k$ multiplications that
interchange only in the lax sense. In this paper ``lax'' will
indicate that the morphisms involved in a definition are not
necessarily isomorphisms. Lax interchangers will obey coherence
axioms, as seen in the next section.

In this paper we follow May by  defining $n$-fold operads in terms
of monoids in a certain category of collections.  A more abstract
approach for future consideration would begin by finding an equivalent
definition in the language of Weber, where an operad lives within
a monoidal pseudo algebra of a 2-monad \cite{web}.  This latter is
a general notion of operad which includes as instances both the
classical operads and the higher operads of Batanin.

\section{\texorpdfstring{$k$}{K}-fold monoidal categories}
This sort of category was developed and defined in \cite{Balt}. The authors describe its structure
as arising recursively from its description as a monoid ${\cal V}, \otimes_k, \lambda$ in the
category of $(k-1)$-fold monoidal categories and lax monoidal functors, with the cartesian product.
Actually, the multiplication $\otimes_k$ is a lax monoidal functor, but the unit $\lambda$ is a
strict monoidal functor. Here we present that definition (in an expanded form) altered only
slightly to make visible the coherent associators as in \cite{forcey1}. That latter paper describes
its structure in terms of tensor objects in the category of $(k-1)$-fold monoidal categories. Our
variation has the effect of making visible the associators $\alpha^{i}_{ABC}.$ It is desirable to
do so for several reasons. One is that the associators are included in definitions of other related
objects such as enriched categories. Another reason is that this inclusion makes easier a direct
comparison with Batanin's definition of monoidal globular categories as in \cite{bat}. A monoidal
globular category can be seen as a quite special case of an iterated monoidal category, with source
and target maps that take objects to those in a category with one less product, and with
interchanges that are isomorphisms.

A third reason is that in this paper we will consider a category of collections in an iterated
monoidal category which will be (iterated) monoidal only up to natural associators.  That being
said, in much of the remainder of this paper we will consider examples with strict associativity,
where each $\alpha$ is the identity, and in interest of clarity will often hide associators.
Another expansion beyond \cite{forcey1} in the following definition is that the products are
allowed to have distinct units.

\begin{definition}\label{iterated} A (strong) $n${\it -fold monoidal category}
with distinct units is a category ${\cal V}$ with the following structure.
\begin{enumerate}
\item There are $n$ multiplications
\[\otimes_1,\otimes_2,\dots, \otimes_n:{\cal V}\times{\cal V}\to{\cal V}\]
each equipped with an associator $\alpha_{UVW}$, a natural isomorphism
which satisfies the pentagon equation:
\end{enumerate}

\noindent
\begin{center}
\resizebox{5.5in}{!}{
$$
\xymatrix@R=45pt@C=-35pt{
&((U\otimes_i V)\otimes_i W)\otimes_i X \text{ }\text{ }
\ar[rr]^{ \alpha^{i}_{UVW}\otimes_i 1_{X}}
\ar[dl]^{ \alpha^{i}_{(U\otimes_i V)WX}}
&\text{ }\text{ }\text{ }\text{ }\text{ }&\text{ }\text{ }(U\otimes_i (V\otimes_i W))\otimes_i X
\ar[dr]^{ \alpha^{i}_{U(V\otimes_i W)X}}&\\
(U\otimes_i V)\otimes_i (W\otimes_i X)
\ar[drr]^{ \alpha^{i}_{UV(W\otimes_i X)}}
&&&&U\otimes_i ((V\otimes_i W)\otimes_i X)
\ar[dll]^{ 1_{U}\otimes_i \alpha^{i}_{VWX}}
\\&&U\otimes_i (V\otimes_i (W\otimes_i X))&&&
}
$$
}
\end{center}
\begin{enumerate}
\item[(2)] ${\cal V}$ has objects $I^i$, $i=1\dots n$, which are strict units
for multiplications: $A\otimes_i I^i = A = I^i \otimes_i A$. Since $\otimes_j$ is a (lax) monoidal
functor (with strict units) with respect to $\otimes_i$ for $1\le i<j\le n$ these units obey
$I^i\otimes_j I^i = I^i.$ Since the unit map (from the single object category to ${\cal V}$ seen as
a monoid with multiplication $\otimes_j$) is itself a strict monoidal functor these units also obey
$I^j\otimes_i I^j = I^j.$

\item[(3)] For each pair $(i,j)$ such that $1\le i<j\le n$ there is an interchanger
natural transformation
\[\eta^{ij}_{ABCD}: (A\otimes_j B)\otimes_i(C\otimes_j D)
    \to (A\otimes_i C)\otimes_j(B\otimes_i D).\]
These natural transformations $\eta^{ij}$ are subject to the following
conditions:
\begin{enumerate}
\item Internal unit condition:
    $\eta^{ij}_{ABI^iI^i}=\eta^{ij}_{I^iI^iAB}=1_{A\otimes_j B}$
\item External unit condition:
    $\eta^{ij}_{AI^jBI^j}=\eta^{ij}_{I^jAI^jB}=1_{A\otimes_i B}$
\item Internal associativity condition: The following diagram commutes.
\end{enumerate}
\end{enumerate}

\begin{noindent}
\begin{center}
\resizebox{6.5in}{!}{
$$
\diagram
((U\otimes_j V)\otimes_i (W\otimes_j X))\otimes_i (Y\otimes_j Z)
\xto[rrr]^{\eta^{ij}_{UVWX}\otimes_i 1_{Y\otimes_j Z}}
\ar[d]^{\alpha^i}
&&&\bigl((U\otimes_i W)\otimes_j(V\otimes_i X)\bigr)\otimes_i (Y\otimes_j Z)
\dto^{\eta^{ij}_{(U\otimes_i W)(V\otimes_i X)YZ}}\\
(U\otimes_j V)\otimes_i ((W\otimes_j X)\otimes_i (Y\otimes_j Z))
\dto^{1_{U\otimes_j V}\otimes_i \eta^{ij}_{WXYZ}}
&&&((U\otimes_i W)\otimes_i Y)\otimes_j((V\otimes_i X)\otimes_i Z)
\ar[d]^{\alpha^i \otimes_j \alpha^i}
\\
(U\otimes_j V)\otimes_i \bigl((W\otimes_i Y)\otimes_j(X\otimes_i Z)\bigr)
\xto[rrr]^{\eta^{ij}_{UV(W\otimes_i Y)(X\otimes_i Z)}}
&&& (U\otimes_i (W\otimes_i Y))\otimes_j(V\otimes_i (X\otimes_i Z))
\enddiagram
$$
}
\end{center}
\end{noindent}
\begin{enumerate}
\item[{}]
\begin{enumerate}
\item[(d)] External associativity condition: The following diagram commutes.
\end{enumerate}
\end{enumerate}
\noindent
\begin{center}
\resizebox{6.5in}{!}{
$$
\diagram
((U\otimes_j V)\otimes_j W)\otimes_i ((X\otimes_j Y)\otimes_j Z)
\xto[rrr]^{\eta^{ij}_{(U\otimes_j V)W(X\otimes_j Y)Z}}
\ar[d]^{\alpha^j \otimes_i \alpha^j}
&&& \bigl((U\otimes_j V)\otimes_i (X\otimes_j Y)\bigr)\otimes_j(W\otimes_i Z)
\dto^{\eta^{ij}_{UVXY}\otimes_j 1_{W\otimes_i Z}}\\
(U\otimes_j (V\otimes_j W))\otimes_i (X\otimes_j (Y\otimes_j Z))
\dto^{\eta^{ij}_{U(V\otimes_j W)X(Y\otimes_j Z)}}
&&&((U\otimes_i X)\otimes_j(V\otimes_i Y))\otimes_j(W\otimes_i Z)
\ar[d]^{\alpha^j}
\\
(U\otimes_i X)\otimes_j\bigl((V\otimes_j W)\otimes_i (Y\otimes_j Z)\bigr)
\xto[rrr]^{1_{U\otimes_i X}\otimes_j\eta^{ij}_{VWYZ}}
&&& (U\otimes_i X)\otimes_j((V\otimes_i Y)\otimes_j(W\otimes_i Z))
\enddiagram
$$
}
\end{center}
\begin{enumerate}
\item[{}]
\begin{enumerate}
\item[(e)] Finally it is required  for each triple $(i,j,k)$ satisfying
$1\le i<j<k\le n$ that the giant hexagonal interchange diagram commutes.
\end{enumerate}
\end{enumerate}
\noindent
\begin{center}
\resizebox{6.5in}{!}{
$$
\xymatrix@C=-118pt{
&((A\otimes_k A')\otimes_j (B\otimes_k B'))\otimes_i((C\otimes_k C')\otimes_j (D\otimes_k D'))
\ar[ddl]|{\eta^{jk}_{AA'BB'}\otimes_i \eta^{jk}_{CC'DD'}}
\ar[ddr]|{\eta^{ij}_{(A\otimes_k A')(B\otimes_k B')(C\otimes_k C')(D\otimes_k D')}}
\\\\
((A\otimes_j B)\otimes_k (A'\otimes_j B'))\otimes_i((C\otimes_j D)\otimes_k (C'\otimes_j D'))
\ar[dd]|{\eta^{ik}_{(A\otimes_j B)(A'\otimes_j B')(C\otimes_j D)(C'\otimes_j D')}}
&&((A\otimes_k A')\otimes_i (C\otimes_k C'))\otimes_j((B\otimes_k B')\otimes_i (D\otimes_k D'))
\ar[dd]|{\eta^{ik}_{AA'CC'}\otimes_j \eta^{ik}_{BB'DD'}}
\\\\
((A\otimes_j B)\otimes_i (C\otimes_j D))\otimes_k((A'\otimes_j B')\otimes_i (C'\otimes_j D'))
\ar[ddr]|{\eta^{ij}_{ABCD}\otimes_k \eta^{ij}_{A'B'C'D'}}
&&((A\otimes_i C)\otimes_k (A'\otimes_i C'))\otimes_j((B\otimes_i D)\otimes_k (B'\otimes_i D'))
\ar[ddl]|{\eta^{jk}_{(A\otimes_i C)(A'\otimes_i C')(B\otimes_i D)(B'\otimes_i D')}}
\\\\
&((A\otimes_i C)\otimes_j (B\otimes_i D))\otimes_k((A'\otimes_i C')\otimes_j (B'\otimes_i D'))
}
$$
}
\end{center}
\end{definition}

As noted in the introduction, the terminology for the case in which the interchangers are
isomorphisms is $k$-tuply monoidal. If the associators $\alpha^i$ are identities then we call the
category strict $n$-fold monoidal; if the units $I^i$ are identical then we say the category has a
common unit $I$; and to follow \cite{Balt} if there is are no quantifiers then we refer to the
strict category with a common unit. If the associators are merely natural transformations then we
call the category lax $n$-fold monoidal (with strict units). The units will always be strict unless
specified, and so the parenthetical specification will be omitted.

Note that in the case of a common unit $I$, for $q>p$ we have natural transformations
\[
\eta^{pq}_{AIIB}:A\otimes_p B\to A\otimes_q B\qquad\text{ and }\qquad
\eta^{pq}_{IABI}:A\otimes_p B\to B\otimes_q A.
\]

Joyal and Street \cite{JS} considered a similar concept to Balteanu,
Fiedorowicz, Schw\"anzl and Vogt's idea of 2--fold
monoidal category.  The former pair required the natural transformation
$\eta_{ABCD}$ to be an isomorphism and showed that the resulting
category is naturally equivalent to a braided monoidal category.
As explained in \cite{Balt}, given such a category one obtains an
equivalent braided monoidal category by discarding one of the two
operations, say $\otimes_2$, and defining the commutativity isomorphism
for the remaining operation $\otimes_1$ to be the composite
\[
\diagram
A\otimes_1 B\rrto^{\eta_{IABI}}
&& B\otimes_2 A\rrto^{\eta_{BIIA}^{-1}}
&& B\otimes_1 A.
\enddiagram
\]

The authors of \cite{Balt} remark that a symmetric monoidal category
is $n$-fold monoidal for all $n$.  This they demonstrate by letting
\[\otimes_1=\otimes_2=\dots=\otimes_n=\otimes\]
and defining
\[\eta^{ij}_{ABCD}=\alpha^{-1}\circ (1_A\otimes \alpha)
    \circ(1_A\otimes (c_{BC}\otimes 1_D))
    \circ(1_A\otimes \alpha^{-1})\circ \alpha \]
for all $i<j$. Here $c_{BC}: B\otimes C \to C\otimes B$ is the symmetry natural transformation.

Joyal and Street \cite{JS} require that the interchange natural
transformations $\eta^{ij}_{ABCD}$ be isomorphisms and observed
that for $n\ge3$ the resulting sort of category is equivalent to a
symmetric monoidal category. Thus as Balteanu, Fiedorowicz, Schw\"anzl
and Vogt point out, the nerves of such categories have group
completions which are infinite loop spaces rather than only $n$--fold
loop spaces.

Because of the nature of the definition of iterated monoidal category, there are multiple forgetful
functors implied. Specifically, letting $n<k$, from the category of $k$-fold monoidal categories to
the category of $n$-fold monoidal categories there are $\binom k n$ forgetful functors which forget
all but the chosen set of products.

The coherence theorem for strict iterated monoidal categories with a common unit states  that any
diagram composed solely of interchange transformations commutes; i.e.  if two compositions of
various interchange transformations (legs of a diagram) have the same source and target then they
describe the same morphism. Furthermore we can easily determine when a composition of interchanges
exists between objects. Here are the necessary definitions and Theorem as given in \cite{Balt}.

\begin{definition} \cite{Balt}
Let ${\cal F}_n(S)$ be the free strict $n$-fold monoidal category on the finite set $S.$ Its
objects are all finite expressions generated by the elements of $S$ using the products $\otimes_k,
k=1..n.$ Its common unit is the empty expression. Its morphisms are finite composites of finite
expressions generated by the symbols $\eta^{ij}_{ABCD}$ with $1\le i < j \le n$, and $A,B,C,D$
objects of ${\cal F}_n(S)$, using the associative operations $\otimes_k, k=1..n.$ Two morphisms are
identified if they can be converted into one another by use of functoriality, naturality and
associativity axioms. By ${\cal M}_n(S)$ we denote the full sub-category of ${\cal F}_n(S)$ whose
objects are expressions in which each element of $S$ occurs exactly once.
\end{definition}

If $S \subset T$ then there is a restriction functor ${\cal M}_n(T)
\to {\cal M}_n(S)$, induced by the functor ${\cal F}_n(T) \to {\cal
F}_n(S)$, which sends $T - S$ to the empty expression 0.

\begin{definition}
Let $A$ be an object of ${\cal M}_n(S).$ For $a,b \in S$ we say
that $a\otimes_ib$ $in$ $A$ if the restriction functor ${\cal M}_n(S)
\to {\cal M}_n({a,b})$ sends $A$ to $a\otimes_ib.$
\end{definition}

\begin{theorem}\cite{Balt}
Let $A$ and $B$ be objects of ${\cal M}_n(S).$ Then
\begin{enumerate}
\item There is at most one morphism $A\to B.$
\item Moreover, there exists a morphism $A\to B$ if and only if, for
every $a,b \in S$, $a\otimes_i b$ in $A$ implies that either $a\otimes_j b$ is in $B$ for some
$j\ge i$ or $b\otimes_j a$ is in $B$ for some $j> i.$
\end{enumerate}
\end{theorem}

Now the coherence theorem can be used to check for commutativity of diagrams in an $n$-fold
monoidal ${\cal V.}$ If two legs of a diagram involving $k$ possibly indistinct operands at its
vertices are formed solely of instances of the interchangers, then they are equal by considering a
morphism from ${\cal M}_n(\{1,\dots ,k\})$ to ${\cal V.}$ Note that if ${\cal V}$ is lax or strong
(associators not identities) then the coherence theorem has not yet been shown to hold. Also note
that if ${\cal V}$ has distinct units then the ``if'' portion of the second part of the coherence
theorem does not imply anything about ${\cal V.}$ However in this case the coherence theorem can
still be used to show commutativity in ${\cal V}$  by considering a morphism from the full
subcategory of ${\cal M}_n(\{1,\dots ,k\})$ formed by discarding the empty expression.

\section{\texorpdfstring{$n$}{N}-fold operads}
The two principle components of an operad are a collection,
historically a sequence, of objects in a monoidal category and a
family of composition maps. Operads are often described as
parameterizations  of $n$-ary operations.  Peter May's original
definition of operad in a symmetric (or braided) monoidal category
\cite{May} has a composition $\gamma$ that takes the tensor product
of the $n\th$ object ($n$-ary operation) and $n$ others (of various
arity) to a resultant that sums the arities of those others.  The
$n\th$ object or $n$-ary operation is often pictured as a tree with
$n$ leaves, and the composition appears like this:

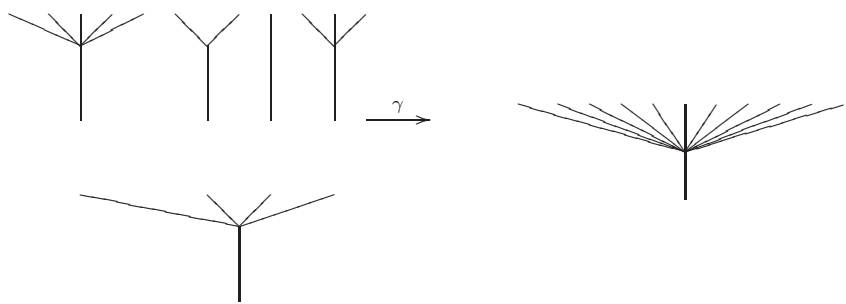
\begin{tabular}{ll}
$$
\xymatrix@W=0pc @H=2.2pc @R=1pc @C=1pc{
*=0{}\ar@{-}[drr] &*=0{}\ar@{-}[dr] &*=0{}\ar@{-}[d] &*=0{}\ar@{-}[dl] &*=0{}\ar@{-}[dll] &*=0{}\ar@{-}[dr] &*=0{} &*=0{}\ar@{-}[dl] &*=0{}\ar@{-}[d] &*=0{}\ar@{-}[dr] &*=0{}\ar@{-}[d] &*=0{}\ar@{-}[dl]\\
*=0{} &*=0{} &*=0{}\ar@{-}[d] &*=0{} &*=0{} &*=0{} &*=0{}\ar@{-}[d] &*=0{} &*=0{}\ar@{-}[d] &*=0{}  &*=0{}\ar@{-}[d] &*=0{} &*=0{}&*=0{}&*=0{}&*=0{}&*=0{}&*=0{}\\&*=0{}&*=0{}&*=0{}&*=0{}&*=0{}&*=0{}&*=0{}&*=0{}&*=0{}&*=0{}&*=0{}\ar[rr]^{\gamma}&*=0{}&*=0{}&*=0{}&*=0{}\\
*=0{} &*=0{} &*=0{}\ar@{-}[drrrrr] &*=0{} &*=0{} &*=0{}&*=0{}\ar@{-}[dr] &*=0{} &*=0{}\ar@{-}[dl] &*=0{} &*=0{}\ar@{-}[dlll] &*=0{} &*=0{}&*=0{}&*=0{}&*=0{}&*=0{}&*=0{}&*=0{}&*=0{}&*=0{}\\
*=0{} &*=0{} &*=0{} &*=0{}      &*=0{} &*=0{} &*=0{}&*=0{}\ar@{-}[d] &*=0{} &*=0{}\\ &*=0{} &*=0{} &*=0{} &*=0{} &*=0{} &*=0{} &*=0{} &*=0{}&*=0{}&*=0{}&*=0{}
}
$$ & \hspace{-6pc}
$$
\xymatrix@W=0pc @H=2.2pc @R=1.5pc @C=1pc{\\
*=0{}\ar@{-}[drrrrr]&*=0{}\ar@{-}[drrrr]&*=0{}\ar@{-}[drrr]&*=0{}\ar@{-}[drr]&*=0{}\ar@{-}[dr]&*=0{}\ar@{-}[d]&*=0{}\ar@{-}[dl]&*=0{}\ar@{-}[dll]&*=0{}\ar@{-}[dlll]&*=0{}\ar@{-}[dllll]&*=0{}\ar@{-}[dlllll]\\
*=0{}&*=0{}&*=0{}&*=0{}&*=0{}&*=0{}\ar@{-}[d]\\
*=0{}*=0{}*=0{}*=0{}*=0{}*=0{}*=0{}&*=0{}&*=0{}&*=0{}&*=0{}&*=0{}&*=0{}&*=0{}&*=0{}&*=0{}\\
}
$$
\end{tabular}

By requiring this composition to be associative we mean that it obeys this sort of pictured commuting diagram:

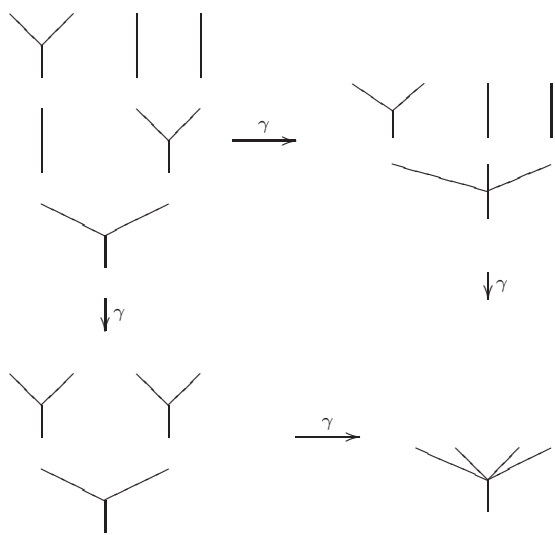
\begin{center}
\begin{tabular}{ll}
$$
\xymatrix@W=0pc @H=2.2pc @R=1pc @C=1pc{
*=0{}\ar@{-}[dr] &*=0{} &*=0{}\ar@{-}[dl] &*=0{} &*=0{}\ar@{-}[d] &*=0{} &*=0{}\ar@{-}[d] \\
*=0{} &*=0{}\ar@{-}[d] &*=0{} &*=0{} &*=0{}\ar@{-}[d] &*=0{} &*=0{}\ar@{-}[d] \\
*=0{}&*=0{}&*=0{}&*=0{}&*=0{}&*=0{}&*=0{}&*=0{}&*=0{}\\
*=0{} &*=0{}\ar@{-}[d] &*=0{} &*=0{} &*=0{}\ar@{-}[dr] &*=0{} &*=0{}\ar@{-}[dl] \\
*=0{} &*=0{}\ar@{-}[d] &*=0{} &*=0{} &*=0{} &*=0{}\ar@{-}[d] &*=0{}&*=0{}\ar[rr]^{\gamma}&*=0{}&*=0{}&*=0{}&*=0{} \\
*=0{}&*=0{}&*=0{}&*=0{}&*=0{}&*=0{}&*=0{}&*=0{}&*=0{}\\
*=0{} &*=0{}\ar@{-}[drr]  &*=0{} &*=0{} &*=0{} &*=0{}\ar@{-}[dll]\\
*=0{} &*=0{} &*=0{} &*=0{}\ar@{-}[d] \\
*=0{} &*=0{} &*=0{} &*=0{} \\
*=0{} &*=0{} &*=0{} &*=0{}\ar[d]^{\gamma}&*=0{}\\
*=0{} &*=0{} &*=0{} &*=0{} &*=0{}
}
$$ & \hspace{-4pc}
$$
\xymatrix@W=0pc @H=2.2pc @R=.85pc @C=1pc{\\
*=0{}\ar@{-}[dr] &*=0{} &*=0{}\ar@{-}[dl] &*=0{} &*=0{}\ar@{-}[d] &*=0{} &*=0{}\ar@{-}[d] \\
*=0{} &*=0{}\ar@{-}[d] &*=0{} &*=0{} &*=0{}\ar@{-}[d] &*=0{} &*=0{}\ar@{-}[d] \\
*=0{}&*=0{}&*=0{}&*=0{}&*=0{}&*=0{}&*=0{}&*=0{}&*=0{}\\
*=0{} &*=0{}\ar@{-}[drrr] &*=0{} &*=0{} &*=0{}\ar@{-}[d] &*=0{} &*=0{}\ar@{-}[dll] \\
*=0{} &*=0{} &*=0{} &*=0{} &*=0{}\ar@{-}[d] &*=0{} &*=0{} \\
*=0{}&*=0{}&*=0{}&*=0{}&*=0{}&*=0{}&*=0{}\\
*=0{} \\
*=0{} &*=0{} &*=0{} &*=0{} &*=0{}\ar[d]^{\gamma}&*=0{}\\
*=0{} &*=0{} &*=0{} &*=0{} &*=0{} &*=0{}
}
$$ \\
$$
\xymatrix@W=0pc @H=2.2pc @R=1pc @C=1pc{
*=0{}\ar@{-}[dr] &*=0{} &*=0{}\ar@{-}[dl] &*=0{} &*=0{}\ar@{-}[dr] &*=0{} &*=0{}\ar@{-}[dl] \\
*=0{} &*=0{}\ar@{-}[d] &*=0{} &*=0{} &*=0{} &*=0{}\ar@{-}[d] &*=0{} \\
*=0{}&*=0{}&*=0{}&*=0{}&*=0{}&*=0{}&*=0{}&*=0{}&*=0{}&*=0{}\ar[rr]^{\gamma}&*=0{}&*=0{}&*=0{}&*=0{} \\
*=0{} &*=0{}\ar@{-}[drr]  &*=0{} &*=0{} &*=0{} &*=0{}\ar@{-}[dll]\\
*=0{} &*=0{} &*=0{} &*=0{}\ar@{-}[d] \\
*=0{} &*=0{} &*=0{} &*=0{} \\
}
$$
& \hspace{-2pc}
$$
\xymatrix@W=0pc @H=2.2pc @R=1pc @C=1pc{\\
*=0{}\ar@{-}[drr] &*=0{}\ar@{-}[dr] &*=0{} &*=0{}\ar@{-}[dl] &*=0{}\ar@{-}[dll]\\
*=0{} &*=0{} &*=0{}\ar@{-}[d] &*=0{} \\
*=0{} &*=0{} &*=0{} &*=0{}
}
$$
\end{tabular}
\end{center}

In the above pictures the tensor products are shown just by
juxtaposition, but now we would like to think about the products
more explicitly. If the monoidal category is not strict, then there
is actually required another leg of the associativity diagram, where
the tensoring is reconfigured  so that the composition can operate
in an alternate order.  Here is how that rearranging looks in a
symmetric (braided) category, where the shuffling is accomplished
by use of the symmetry (braiding):

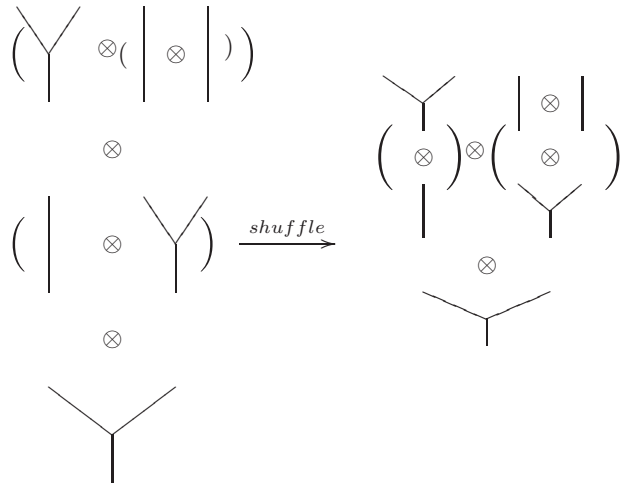
\begin{center}
\begin{tabular}{ll}
$$
\xymatrix@W=0pc @H=2.2pc @R=1.5pc @C=1pc{
*=0{}\ar@{-}[dr] &*=0{} &*=0{}\ar@{-}[dl] &*=0{} &*=0{}\ar@{-}[d] &*=0{} &*=0{}\ar@{-}[d] \\
*=0{\txt{\huge(}} &*=0{}\ar@{-}[d] &*=0{} &*=0{\otimes\txt{(}} &*=0{}\ar@{-}[d] &*=0{\otimes} &*=0{}\ar@{-}[d] &*=0{\left.\right)\txt{\huge)}}\\
*=0{}&*=0{}&*=0{}&*=0{}&*=0{}&*=0{}&*=0{}&*=0{}&*=0{}\\
*=0{}&*=0{}&*=0{}&*=0{\otimes}&*=0{}&*=0{}&*=0{}&*=0{}&*=0{}\\
*=0{} &*=0{}\ar@{-}[d] &*=0{} &*=0{} &*=0{}\ar@{-}[dr] &*=0{} &*=0{}\ar@{-}[dl] \\
*=0{\txt{\huge(}} &*=0{}\ar@{-}[d] &*=0{} &*=0{\otimes} &*=0{} &*=0{}\ar@{-}[d] &*=0{\txt{\huge)}}&*=0{}\ar[rrr]^{shuffle}&*=0{}&*=0{}&*=0{}&*=0{} \\
*=0{}&*=0{}&*=0{}&*=0{}&*=0{}&*=0{}&*=0{}&*=0{}&*=0{}\\
*=0{}&*=0{}&*=0{}&*=0{\otimes}&*=0{}&*=0{}&*=0{}&*=0{}&*=0{}\\
*=0{} &*=0{}\ar@{-}[drr]  &*=0{} &*=0{} &*=0{} &*=0{}\ar@{-}[dll]\\
*=0{} &*=0{} &*=0{} &*=0{}\ar@{-}[d] \\
*=0{} &*=0{} &*=0{} &*=0{} \\
}
$$ & \hspace{-1.25pc}
$$
\xymatrix@W=0pc @H=2.2pc @R=.85pc @C=1pc{\\
*=0{}\ar@{-}[dr] &*=0{} &*=0{}\ar@{-}[dl] &*=0{} &*=0{}\ar@{-}[d] &*=0{} &*=0{}\ar@{-}[d] \\
*=0{} &*=0{}\ar@{-}[d] &*=0{} &*=0{} &*=0{}\ar@{-}[d] &*=0{\otimes} &*=0{}\ar@{-}[d] &*=0{}\\
*=0{}&*=0{}&*=0{}&*=0{}&*=0{}&*=0{}&*=0{}&*=0{}&*=0{}\\
*=0{\txt{\Huge(}}&*=0{\otimes}&*=0{\txt{\Huge)}}&*=0{\otimes\txt{\Huge(}}&*=0{}&*=0{\otimes}&*=0{}&*=0{\txt{\Huge)}}&*=0{}\\
*=0{} &*=0{}\ar@{-}[d] &*=0{} &*=0{} &*=0{}\ar@{-}[dr] &*=0{} &*=0{}\ar@{-}[dl] \\
*=0{} &*=0{}\ar@{-}[d] &*=0{} &*=0{} &*=0{} &*=0{}\ar@{-}[d] &*=0{}&*=0{} \\
*=0{}&*=0{}&*=0{}&*=0{}&*=0{}&*=0{}&*=0{}&*=0{}&*=0{}\\
*=0{}&*=0{}&*=0{}&*=0{\otimes}&*=0{}&*=0{}&*=0{}&*=0{}&*=0{}\\
*=0{} &*=0{}\ar@{-}[drr]  &*=0{} &*=0{} &*=0{} &*=0{}\ar@{-}[dll]\\
*=0{} &*=0{} &*=0{} &*=0{}\ar@{-}[d] \\
*=0{} &*=0{} &*=0{} &*=0{} \\
}
$$
\end{tabular}
\end{center}

We now foreshadow our definition of operads in an iterated monoidal
category with the same picture as above but using two tensor products,
$\otimes_1$ and $\otimes_2.$ It becomes clear that the true nature
of the shuffle is in fact that of an interchange transformation.

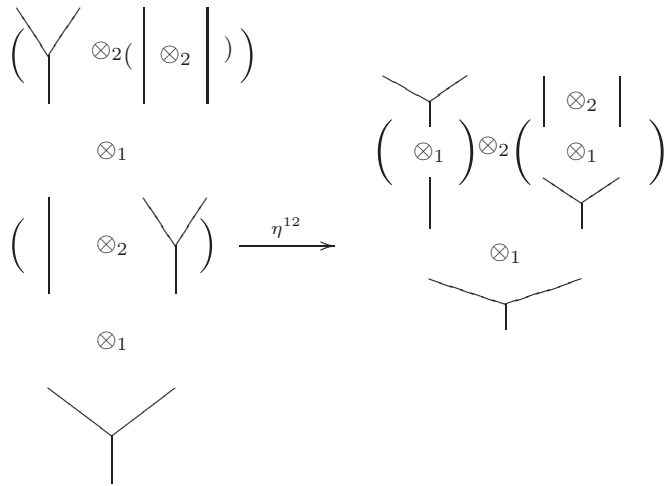
\begin{center}
\begin{tabular}{ll}
$$
\xymatrix@W=0pc @H=2.2pc @R=1.5pc @C=1pc{
*=0{}\ar@{-}[dr] &*=0{} &*=0{}\ar@{-}[dl] &*=0{} &*=0{}\ar@{-}[d] &*=0{} &*=0{}\ar@{-}[d] \\
*=0{\txt{\huge(}} &*=0{}\ar@{-}[d] &*=0{} &*=0{\otimes_2\txt{(}} &*=0{}\ar@{-}[d] &*=0{\otimes_2} &*=0{}\ar@{-}[d] &*=0{\left.\right)\txt{\huge)}}\\
*=0{}&*=0{}&*=0{}&*=0{}&*=0{}&*=0{}&*=0{}&*=0{}&*=0{}\\
*=0{}&*=0{}&*=0{}&*=0{\otimes_1}&*=0{}&*=0{}&*=0{}&*=0{}&*=0{}\\
*=0{} &*=0{}\ar@{-}[d] &*=0{} &*=0{} &*=0{}\ar@{-}[dr] &*=0{} &*=0{}\ar@{-}[dl] \\
*=0{\txt{\huge(}} &*=0{}\ar@{-}[d] &*=0{} &*=0{\otimes_2} &*=0{} &*=0{}\ar@{-}[d] &*=0{\txt{\huge)}}&*=0{}\ar[rrr]^{\eta^{12}}&*=0{}&*=0{}&*=0{}&*=0{} \\
*=0{}&*=0{}&*=0{}&*=0{}&*=0{}&*=0{}&*=0{}&*=0{}&*=0{}\\
*=0{}&*=0{}&*=0{}&*=0{\otimes_1}&*=0{}&*=0{}&*=0{}&*=0{}&*=0{}\\
*=0{} &*=0{}\ar@{-}[drr]  &*=0{} &*=0{} &*=0{} &*=0{}\ar@{-}[dll]\\
*=0{} &*=0{} &*=0{} &*=0{}\ar@{-}[d] \\
*=0{} &*=0{} &*=0{} &*=0{} \\
}
$$ & \hspace{-1.25pc}
$$
\xymatrix@W=0pc @H=2.2pc @R=.8pc @C=1.2pc{\\
*=0{}\ar@{-}[dr] &*=0{} &*=0{}\ar@{-}[dl] &*=0{} &*=0{}\ar@{-}[d] &*=0{} &*=0{}\ar@{-}[d] \\
*=0{} &*=0{}\ar@{-}[d] &*=0{} &*=0{} &*=0{}\ar@{-}[d] &*=0{\otimes_2} &*=0{}\ar@{-}[d] &*=0{}\\
*=0{}&*=0{}&*=0{}&*=0{}&*=0{}&*=0{}&*=0{}&*=0{}&*=0{}\\
*=0{\txt{\Huge(}}&*=0{\otimes_1}&*=0{\txt{\Huge)}}&*=0{\otimes_2\txt{\Huge(}}&*=0{}&*=0{\otimes_1}&*=0{}&*=0{\txt{\Huge)}}&*=0{}\\
*=0{} &*=0{}\ar@{-}[d] &*=0{} &*=0{} &*=0{}\ar@{-}[dr] &*=0{} &*=0{}\ar@{-}[dl] \\
*=0{} &*=0{}\ar@{-}[d] &*=0{} &*=0{} &*=0{} &*=0{}\ar@{-}[d] &*=0{}&*=0{} \\
*=0{}&*=0{}&*=0{}&*=0{}&*=0{}&*=0{}&*=0{}&*=0{}&*=0{}\\
*=0{}&*=0{}&*=0{}&*=0{\otimes_1}&*=0{}&*=0{}&*=0{}&*=0{}&*=0{}\\
*=0{} &*=0{}\ar@{-}[drr]  &*=0{} &*=0{} &*=0{} &*=0{}\ar@{-}[dll]\\
*=0{} &*=0{} &*=0{} &*=0{}\ar@{-}[d] \\
*=0{} &*=0{} &*=0{} &*=0{} \\
}
$$
\end{tabular}
\end{center}

To see this just focus on the actual domain and range of $\eta^{12}$
which are the upper two levels of trees in the pictures, with the
tensor product $\left({\mathbf|}\otimes_2{\mathbf|}\right)$ considered
as a single object.

Now we are ready to give the technical definitions. We begin with
the definition of 2-fold operad in an $n$-fold monoidal category,
as in the above picture, and then mention how it generalizes the
case of operad in a braided category. Because of this generalization
of the well known case, and since there are easily described examples
of 2-fold monoidal categories based on a braided category as in
\cite{forcey3}, it seems worthwhile to work out the theory for the
2-fold operads in its entirety before moving on to $n$-fold operads.

\begin{definition}
Let ${\cal V}$ be a strict $n$-fold monoidal category as defined in
Section 2. A 2-fold operad ${\cal C}$ in ${\cal V}$ consists of
objects ${\cal C}(j)$, $j\ge 0$, a unit map ${\cal J}:I\to {\cal
C}(1)$, and composition maps in ${\cal V}$
\[
\gamma^{12}:{\cal C}(k) \otimes_1 ({\cal C}(j_1) \otimes_2 \dots \otimes_2 {\cal C}(j_k))\to {\cal C}(j)
\]
for $k\ge 1$, $j_s\ge0$ for $s=1\dots k$ and
$\smash{\sum\limits_{n=1}^k j_n = j}$. The composition maps obey
the following axioms:
\begin{enumerate}
\item Associativity: The following diagram is required to commute
for all  $k\ge 1$, $j_s\ge 0$ and $i_t\ge 0$, and where
$\sum\limits_{s=1}^k j_s = j$ and $\sum\limits_{t=1}^j i_t = i.$ Let
$g_s= \sum\limits_{u=1}^s j_u$ and let
$h_s=\sum\limits_{u=1+g_{s-1}}^{g_s} i_u$.  The $\eta^{12}$ labeling
the leftmost arrow actually stands for a variety of equivalent maps
which factor into instances of the $12$ interchange.
\[
\xymatrix{
{\cal C}(k)\otimes_1\left(\bigotimes\limits_{s=1}^k{}_2 {\cal C}(j_s)\right)\otimes_1
\left(\bigotimes\limits_{t=1}^j{}_2 {\cal C}(i_t)\right)
\ar[rr]^>>>>>>>>>>>>{\gamma^{12} \otimes_1 \text{id}}
\ar[dd]_{\text{id} \otimes_1 \eta^{12}}
&&{\cal C}(j)\otimes_1 \left(\bigotimes\limits_{t=1}^j{}_2 {\cal C}(i_t)\right)
\ar[d]^{\gamma^{12}}
\\
&&{\cal C}(i)
\\
{\cal C}(k)\otimes_1 \left(\bigotimes\limits_{s=1}^k{}_2 \left({\cal
C}(j_s)\otimes_1 \left(\bigotimes\limits_{u=1}^{j_s}{}_2 {\cal
C}(i_{u+g_{s-1}})\right)\right)\right) \ar[rr]_>>>>>>>>>>{\text{id}
\otimes_1(\otimes^k_2\gamma^{12})} &&{\cal
C}(k)\otimes_1\left(\bigotimes\limits_{s=1}^k{}_2 {\cal
C}(h_s)\right) \ar[u]_{\gamma^{12}} }
\]

\item Respect of units is required just as in the symmetric case.
The following unit diagrams commute.
\[
\xymatrix{
{\cal C}(k)\otimes_1 (\otimes_2^k I)
\ar[d]_{1\otimes_1(\otimes_2^k {\cal J})}
\ar@{=}[r]^{}
&{\cal C}(k)\\
{\cal C}(k)\otimes_1(\otimes_2^k {\cal C}(1))
\ar[ur]_{\gamma^{12}}
}
\xymatrix{
I\otimes_1 {\cal C}(k)
\ar[d]_{{\cal J}\otimes_1 1}
\ar@{=}[r]^{}
&{\cal C}(k)\\
{\cal C}(1)\otimes_1 {\cal C}(k)
\ar[ur]^{\gamma^{12}}
}
\]
\end{enumerate}
\end{definition}

Note that operads in a braided monoidal category are examples of
2-fold operads. This is true based on the arguments of Joyal and
Street \cite{JS}, who showed that braided categories arise as 2-fold
monoidal categories where the interchanges are isomorphisms. Also
note that given such a perspective on a braided category, the two
products are equivalent and the use of the braiding to shuffle in
the operad associativity requirement can be rewritten as the use
of the interchange.

It is immediately clear that we can define operads using more than
just the first two products in an $n$-fold monoidal category. The
best way of going about this is to use the theory of monoids, (and
more generally enriched categories), in iterated monoidal
categories. We continue by first describing this procedure for
2-fold operads. Operads in a symmetric (braided) monoidal category
are often efficiently defined as the monoids of a category of
collections. For a braided category $({\cal V}, \otimes)$ with
coproducts that are preserved by both functors $ (\dash\otimes A)$
and $(A\otimes\dash)$ the objects of $\col({\cal V})$ are functors
from the discrete category of nonnegative integers to ${\cal V}.$ In
other words the data for a collection ${\cal C}$ is a sequence of
objects ${\cal C}(j).$ Morphisms in $\col({\cal V})$ are natural
transformations. The tensor product in $\col({\cal V})$ is given by
\[
({\cal B}\otimes{\cal C})(j) =
    \coprod\limits_{\substack{k\ge 0\\ j_1+\dots +j_k =j}}
    {\cal B}(k) \otimes({\cal C}(j_1) \otimes \dots \otimes {\cal C}(j_k))
\]
where $j_i \ge 0.$ This product is associative by use of the symmetry
or braiding, and due to the hypothesis that the tensor product
preserves the coproduct.  The unit is the collection
$(\emptyset, I, \emptyset, \dots)$ where $\emptyset$ is an initial
object in ${\cal V}.$

Now recall how the interchange transformations generalize braiding.
For ${\cal V}$ a $2$-fold monoidal category with all coproducts in
which both $\otimes_1$ and $\otimes_2$ strictly preserve the
coproduct, define the objects and morphisms of      $ \col(\cal V)$
in precisely the same way as in the braided case, but define the
product to be
\[
({\cal B}{\otimes}^{12}{\cal C})(j) =
    \coprod\limits_{\substack{k\ge 0\\ j_1+\dots +j_k =j}}
    {\cal B}(k) \otimes_1({\cal C}(j_1) \otimes_2
    \dots \otimes_2 {\cal C}(j_k))
\]
In general the interchangers will not be isomorphisms, so this
product can not be that of a monoidal category with the usual strong
associativity.  However the interchangers can be used to make the
product in question obey lax associativity, where the associator is
a coherent natural transformation: it obeys the usual pentagon axiom
but is not required to be an isomorphism. This lax associativity is
seen by inspection of the two 3-fold products $({\cal
B}{\otimes}^{12}{\cal C}){\otimes}^{12}{\cal D}$ and ${\cal
B}{\otimes}^{12}({\cal C}{\otimes}^{12}{\cal D})$. In the braided
case mentioned above, the two large coproducts in question are seen
to be composed of the same terms up to a braiding between them. Here
the terms of the two coproducts are related by instances of the
interchange transformation $\eta^{12}$ from the term in $(({\cal
B}{\otimes}^{12}{\cal C}){\otimes}^{12}{\cal D})(j)$ to the
corresponding term in $({\cal B}{\otimes}^{12}({\cal
C}{\otimes}^{12}{\cal D}))(j).$ For example upon expansion of the
two three-fold products we see that in the coproduct which is
$(({\cal B}{\otimes}^{12}{\cal C}){\otimes}^{12}{\cal D})(2)$ we
have the term
\[{\cal B}(2)\otimes_1({\cal C}(1)\otimes_2{\cal C}(1))\otimes_1({\cal D}(1)\otimes_2{\cal D}(1))\]
while in
$({\cal B}{\otimes}^{12}({\cal C}{\otimes}^{12}{\cal D}))(2)$
we have the term
\[{\cal B}(2)\otimes_1({\cal C}(1)\otimes_1{\cal D}(1))\otimes_2({\cal C}(1)\otimes_1{\cal D}(1)).\]
Note that the first of these terms appears courtesy of the fact that
strict preservation of coproducts by the tensor product means
precisely that there is a distributive law $(\coprod B_n) \otimes A
= \coprod(B_n\otimes A)$.

Commutativity of the pentagon axiom for the associators is implied by functoriality of the products
and by the coherence theorem for strict $n$-fold monoidal categories,  since the legs of the
diagram are made of distributions over the coproduct (identities) and of  compositions of
interchangers $\eta^{12}$ in ${\cal V}$. Some remarks about the non-invertibility of $\alpha$ are
in order. Note that  Mac Lane proves his coherence theorem in two steps \cite{MacLane}. First it is
shown that every diagram involving only $\alpha$ (no $\alpha^{-1}$) commutes. Then it is noted that
this suffices to make every diagram of both $\alpha$ and $\alpha^{-1}$ commute since for every
binary word there exists a path of just instances of $\alpha$ from that word to the word
parenthesized all to the right. (Here we are taking the domain of $\alpha$ to be $(A\otimes
B)\otimes C.$) Thus when $\alpha$ is not invertible we still have that every diagram commutes.
There are still canonical maps from every binary word to the word parenthesized all to the right.
However there are necessarily fewer diagrams. For instance if $({\cal V}, \otimes)$ is lax monoidal
there is no canonical map between the two objects $(B\otimes B)\otimes(B\otimes B)$ and $(B\otimes
(B\otimes B))\otimes B$. This affects the statement of the general associativity theorem for
monoids in a lax monoidal category. Only the specific case of the general associativity theorem as
stated by Mac Lane  holds, as follows.

\begin{theorem}
Let $(A,\mu)$ be a monoid in a (lax) monoidal category. Let $A^{n}$
be the product given by $A^2 = A \otimes A, A^{n+1} = A\otimes A^n$,
i.e. parenthesized to the right. Define the composition $\mu^{(n)}$
by $\mu^{(2)} = \mu, \mu^{(n+1)} = \mu \circ (1 \otimes \mu^{(n)})$.
Then
\[ \mu^{(n)} \circ (\mu^{(k_1)}\otimes \dots \otimes \mu^{(k_n)}) = \mu^{(k_1 + \dots +k_n)}\circ \alpha' \]
for all $n, k_i \ge 2$ where $\alpha'$ stands for the canonical map
to $A^{k_1 + \dots +k_n}.$
\end{theorem}

\begin{proof}
This is just the special case of the general associative law for
monoids shown by Mac Lane, which only depends on the existence of
the canonical map $\alpha'$ \cite{MacLane}.
\end{proof}
Now we have a condensed way of defining 2-fold operads.
\begin{theorem}
2-fold operads in $2$-fold monoidal ${\cal V}$ are monoids in
$\col(\cal V).$
\end{theorem}
\begin{proof}
A monoid in $ \col(\cal V)$ is an object ${\cal C}$ in $\col(\cal
V)$ with multiplication and unit morphisms. Since morphisms of
$\col(\cal V)$ are natural transformations the multiplication and
unit consist of families of maps in ${\cal V}$ indexed by the
natural numbers, with source and target exactly as required for
operad composition and unit. The operad axioms are equivalent to the
associativity and unit requirements of monoids.
\end{proof}

This brings us back to the question of defining operads in an
$n$-fold monoidal ${\cal V}$ using the higher products and
interchanges. This idea will correspond to a series of higher
products, denoted by $\otimes^{pq},$ in the category of collections.
These are defined just as for the first case $\otimes^{12}$ above.
Associators are as described above for the first product, using
$\eta^{pq}$ for the associator $\alpha\colon{\cal
A}\otimes^{pq}({\cal B}\otimes^{pq}{\cal C})
    \to ({\cal A}\otimes^{pq}{\cal B})\otimes^{pq}{\cal C}$.
The unit for each is the collection  $(\emptyset, I, \emptyset,
\dots)$ where $\emptyset$ is an initial object in ${\cal V}.$ Notice
that these products do not interchange; i.e they are not functorial
with respect to each other.  Notice also that the associators in
these categories of collections are not isomorphisms unless we are
considering the special cases of braiding or symmetry. Instead the
category of collections with substitution product $\otimes^{pq}$ is
lax monoidal, by which we will mean that the associator is merely a
natural transformation, which obeys the pentagon coherence condition
by the same argument as for $\otimes^{12}$ above.

Now we will focus on the products $\otimes^{(m-1)m}$ in the category
of collections in  $n$-fold monoidal ${\cal V}$, for $m\le n,$ since
these will be seen to suffice for defining all operad compositions.
Before defining $m$-fold operads as monoids with respect to
$\otimes^{(m-1)m},$ we note that there is also fibrewise monoidal
structure.  This will be important in the description of the monoidal
structure of the category of operads. In fact, we have the following

\begin{definition}
Let $ ({\cal V},\otimes_1,\dots, \otimes_n)$  be a strict $n$-fold
monoidal category with coproducts and an initial object $\emptyset.$
For $i=1\dots n$ and $A$ an object of ${\cal V}$ let each of the
functors $(\dash\otimes_i A)$ and $(A\otimes_i\dash)$ preserve
coproducts and let $\emptyset\otimes_i A = \emptyset$ and
$A\otimes_i\emptyset = \emptyset $. Let $({\cal
V},\coprod,\otimes_1,\dots, \otimes_n)$ be a strict $(n+1)$-fold
monoidal category with distinct units for which forgetting the first
tensor product (given by the coproduct) recovers $ {\cal V}.$ The
unit for $\coprod$ is the initial object $\emptyset.$  Let $n\ge m
\ge 2.$ We denote by $\col_m(\cal V)$ the category of collections in
${\cal V}$ with the following products:

\begin{align*}
({\cal B}\hat{\otimes}_1{\cal C})(j) &=
    \coprod\limits_{\substack {k\ge 0\\ j_1+\dots +j_k =j}}
    {\cal B}(k) \otimes_{m-1}
    ({\cal C}(j_1) \otimes_m \dots \otimes_m {\cal C}(j_k))
\intertext{and}
({\cal B}\hat{\otimes}_2{\cal C})(j) &= {\cal B}(j){\otimes}_{m+1}{\cal C}(j) \\
&\hspace{0.5em}\vdots\\
({\cal B}\hat{\otimes}_{n-m+1}{\cal C})(j) &= {\cal
B}(j){\otimes}_{n}{\cal C}(j)
\end{align*}
\end{definition}

\begin{theorem}\label{foo}
 $\col_m(\cal V)$ is an $(n-m+1)$-fold lax monoidal category with (two) distinct strict units.
\end{theorem}

\begin{proof}
The first tensor product is $\hat{\otimes}_1 = \otimes^{(m-1)m}$ and
the others are the higher fibrewise products starting with fibrewise
$\otimes_{m+1}.$ Thus the  unit for $\hat{\otimes}_1$ is ${\cal I} =
(\emptyset, I, \emptyset, \dots)$ and the unit for all the other
products is ${\cal M} = (I, I, \dots)$. It is  not hard to check the
unit conditions which are required for the fibrewise products to be
the multiplication for a monoid in the category of monoidal
categories.  The extra requirement of the two sorts of unit is that
${\cal M} \hat{\otimes}_1 {\cal M} = {\cal M} $ and that ${\cal I}
\hat{\otimes}_k {\cal I} = {\cal I}$ for $k >1.$ These equations do
indeed hold.

 Now we must check that
there are interchangers, natural transformations
\[\xi^{1j}: ({\cal A}\hat{\otimes}_j{\cal B})\hat{\otimes}_1({\cal C}\hat{\otimes}_j{\cal D})
    \to ({\cal A}\hat{\otimes}_1{\cal C})\hat{\otimes}_j({\cal B}\hat{\otimes}_1{\cal D}).\]
These utilize the $\eta^{ij}$ of ${\cal V}$ and thus exist by inspection of the terms of the
compound products.  For example, in the product in $\col_2{\cal V}:$
\[ (({\cal A}\hat{\otimes}_2 {\cal B})\hat{\otimes}_1({\cal C}\hat{\otimes}_2{\cal D}))(2) \]
we find the term
\[ \left(({\cal A}(1)\otimes_3{\cal B}(1))\otimes_1({\cal C}(2)\otimes_3{\cal D}(2))\right) \coprod \left(({\cal A}(2)\otimes_3 {\cal B}(2))\otimes_1(({\cal C}(1)\otimes_3 {\cal D}(1))
    \otimes_2({\cal C}(1)\otimes_3 {\cal D}(1)))\right) \]
while in
\[ (({\cal A}\hat{\otimes}_1{\cal C})\hat{\otimes}_2({\cal B}\hat{\otimes}_1{\cal D}))(2) \]
we find the term
\[\left( ({\cal A}(1)\otimes_1{\cal C}(2))\coprod ({\cal A}(2)\otimes_1 ({\cal C}(1)\otimes_2 {\cal C}(1)))\right)\otimes_3
    \left(({\cal B}(1)\otimes_1{\cal D}(2))\coprod({\cal B}(2)\otimes_1({\cal D}(1)\otimes_2 {\cal D}(1)))\right) \]
 The map
$\xi^{12}$ thus uses first $\eta^{23}$, then $\eta^{13}$ and finally the hypothesis that $( {\cal
V},\coprod,\otimes_1,\dots, \otimes_n)$ is an $(n+1)$-fold monoidal category; specifically
instances of the interchange between $\coprod$ and $\otimes_3.$  Thus the external associativity
axiom for $\xi^{1j}$   and the giant hexagon axiom for $i,j,k= 1,j,k$ hold due to the coherence
theorem for strict iterated monoidal categories applied to the $(n+1)$-fold monoidal category $(
{\cal V},\coprod,\otimes_{1},\dots, \otimes_n).$ The internal associativity axiom holds due to
functoriality of the products (since the associators for $\otimes^{pq}$ use an equality followed by
an interchange) as well as coherence. Note that in this category there are distinct strict units,
$\emptyset$ and $I$. Therefore the second part of the coherence theorem which describes existence
of maps does not hold. However since we verify the existence of maps by inspection, the first part
of the coherence theorem applies, which states that any two morphisms (sharing a common source and
target and both being formed of instances of the interchanges which obey the axioms of \cite{Balt})
are the same morphism.

The unit conditions for the interchangers $\xi^{ij}$ are seen to
hold based on the unit conditions for the $(n+1)$-fold monoidal
category $( {\cal V},\coprod,\otimes_{1},\dots, \otimes_n)$ and on
the hypothesis that $\emptyset\otimes_i A = \emptyset$ and
$A\otimes_i\emptyset = \emptyset .$
 Thus the first product
together with any of the fibrewise products are those of a 2-fold
lax monoidal category.

For the products $\hat{\otimes}_{2}$ and higher the associators and
interchange transformations are fibrewise and the axioms hold since
they hold for each fiber.
\end{proof}

\begin{remark}
In the context of \cite{batainf} the  lax functoriality of the
tensor product with respect to the coproduct is due to the
hypothesis that the symmetric monoidal category ${\cal V}$ is closed
(from the right) with respect to the tensor product.  This
guarantees that that product preserves colimits on the first
operand, since the functor $(\dash\otimes B)$ has as a right adjoint
the internal hom, denoted by  $[B,\dash].$ Applied to the coproduct
this fact in turn implies that there is a canonical map in ${\cal
V}$ from $(A \otimes B)\coprod(C\otimes D)$ to $(A \coprod C)\otimes
(B\coprod D).$ From the universal properties of the coproduct it can
be  checked that this map satisfies the the middle interchange law
that is required of a monoidal functor.  Also in \cite{batainf}
Batanin points out that a fibrewise product is a monoidal functor
with respect to the collection product. In that paper the existence
of the transformation $\xi$ depends on the symmetry (braiding) and
the lax functoriality of the tensor product with respect to the
coproduct. In this paper we chose to simply include the necessary
iterated monoidal structure as a hypothesis, rather than the
hypothesis of closedness, in the interest of generality.

Theorem~\ref{foo} is quite useful for describing $n$-fold operads
and their higher-categorical structure, especially when coupled with
two other facts. The first is that monoids are equivalently defined
as single object enriched categories, and the second is the
following result from \cite{forcey1} and \cite{forcey2}, where the
concept of $n$-fold monoidal 2-category is discussed. In those
sources the quantifier lax is sometimes left off, but the proofs in
question nowhere require the associator to be an isomorphism.
\end{remark}

\begin{theorem}\label{enrich}
For ${\cal V}$ $n$-fold (lax) monoidal with distinct strict units
the category of enriched categories over $({\cal V}, \otimes_1)$ is
an $(n-1)$-fold (lax) monoidal 2-category with distinct strict
units.
\end{theorem}
\begin{proof}
The proof in \cite{forcey1} holds exactly as it is stated but with
the modification that given distinct strict units $I^i$ for the
products $\otimes_i$ of ${\cal V}$, we define unit enriched
categories ${\bcal I}^i$ in ${\cal V}$-Cat to each have the single
object $0$ and to have ${\bcal I}^i(0,0) = I^{i+1}$. Then the unit
axioms for distinct units play the same role in the new proof as did
the standard unit axioms in the original.
\end{proof}

For our purposes we translate the theorem about enriched categories
into its single object corollary about the category $\mon(\cal V)$
of monoids in ${\cal V}$.

\begin{corollary}\label{mon}
For ${\cal V}$ $n$-fold (lax) monoidal with distinct strict units,
the category $\mon(\cal V)$ is an $(n-1)$-fold (lax) monoidal
2-category with distinct strict units.
\end{corollary}
\begin{proof}
The product of enriched categories always has as its object set the
cartesian product of the object sets of its components. Thus one
object enriched categories have products with one object as well.
\end{proof}

\begin{definition}\label{itop}
If an $n$-fold monoidal category $ {\cal V}$ has coproducts and
$({\cal V},\coprod,\otimes_3,\dots, \otimes_n)$ is an $(n-1)$-fold
monoidal category in which each of the functors $(\dash\otimes_i A)$
and $(A\otimes_i\dash), i=1\dots n$ preserves coproducts we define
the category of $m$-fold operads $\oper_m(\cal V)$ to be the
category of monoids in the category of collections $( \col_m(\cal
V), \hat{\otimes}_1)$ for $n\ge m \ge 2.$
\end{definition}

\begin{corollary}\label{op}
$\oper_m(\cal V)$ is an $(n-m)$-fold monoidal 2-category.
\end{corollary}

\begin{proof}
Rather than starting with monoids in an $n$-fold monoidal ${\cal V}$
as in the previous corollary we are actually beginning with monoids
in $(n-m+1)$-fold monoidal $\col_m(\cal V)$.  Note that in
\cite{forcey1} the products in ${\cal V}$ are assumed to have a
common unit. To generalize to our situation here, where the unit for
the first product in the category of collections is distinct from
the other units, we need to add slightly to the definitions in
\cite{forcey1}. When enriching (or more specifically taking monoids)
we are doing so with respect to the first available product. Thus
the unit morphism for enriched categories has its domain the unit
for that first product, ${\cal I}.$ However the unit enriched
category ${\bcal I}$ has one object, denoted $0$, and ${\bcal
I}(0,0) = {\cal M}$.
\end{proof}

\begin{remark}
This last corollary justifies our focus on the first $m$ products of
${\cal V}$ as opposed to any subset of the $n$ products. Our choice
of focus is due to the way in which this focus allows us to describe
the resulting structure on the category of $m$-fold operads. Of
course, we can use the forgetful functors mentioned in Section 2 to
pass from $n$-fold monoidal ${\cal V}$ to ${\cal V}$ with any of the
subsets of products.  The $m$-fold operads do behave as expected
under this forgetting, retaining all but the structure which depends
on the forgotten products. This will be seen more clearly upon
inspection of the unpacked definition to follow.  In short, we will
see that an $m$-fold operad is also an $(m-1)$-fold operad.
\end{remark}

\begin{remark}
We note that since a symmetric monoidal category is $n$-fold
monoidal for all $n$, then operads in a symmetric monoidal category
are $n$-fold monoidal for all $n$ as well.  More generally, if $n\ge
3$ and the interchanges are isomorphisms, then by the Eckmann-Hilton
argument the products collapse into one and the result is a
symmetric monoidal category, and so operads in it are again $n$-fold
monoidal for all $n.$ Here we are always discussing ordinary
``non-symmetric,'' (``non-braided'') operads. The possible faithful
actions of symmetry or braid groups can be considered after the
definition, which we leave for a later paper. We do point out that
the proper direction in which to expand this work is seen in Weber's
paper \cite{web}. He generalizes by making a distinction between the
binary and $k$-ary products in the domain of the composition map
$\gamma\colon{\cal C}(k)\otimes({\cal C}(j_1)\otimes
    \dots \otimes{\cal C}(j_k))
    \to {\cal C}(j)$.
The binary tensor product is seen formally as a pseudo-monoid
structure and the $k$-ary product as a pseudo-algebra structure for
a 2-monad which can contain the information needed to describe
actions of braid or symmetry groups. The two structures are defined
using strong monoidal morphisms, and so the products coincide and
give rise to the braiding which is used to describe the associativity
of composition. To encompass the definitions in this paper we would
move to operads in lax-monoidal pseudo algebras, where instead of
pseudo monoids and strong monoidal morphisms in a pseudo algebra
we would consider the same picture but with lax monoidal morphisms.
\end{remark}

The fact that monoids are single object enriched categories also
leads to an efficient expanded definition of $m$-fold operads in an
$n$-fold monoidal category.  Let ${\cal V}$ be an $n$-fold monoidal
category.

\begin{definition}\label{expandop}
For $2\le m\le n$ an $m$-fold operad ${\cal C}$ in ${\cal V}$ consists of
objects ${\cal C}(j)$, $j\ge 0$,
a unit map ${\cal J}:I\to {\cal C}(1)$,
and composition maps in ${\cal V}$
\[
\gamma^{pq}:{\cal C}(k) \otimes_p ({\cal C}(j_1) \otimes_q \dots \otimes_q {\cal C}(j_k))\to {\cal C}(j)
\]
for $m\ge q>p \ge 1$, $k\ge 1$, $j_s\ge0$ for $s=1\dots k$ and $\sum\limits_{n=1}^k j_n = j$. The composition maps obey the following axioms:
\begin{enumerate}
\item Associativity: The following diagram is required to commute for all $m\ge q>p \ge 1$, $k\ge 1$, $j_s\ge 0$ and $i_t\ge 0$, and
where $\sum\limits_{s=1}^k j_s = j$ and $\sum\limits_{t=1}^j i_t = i.$ Let $g_s= \sum\limits_{u=1}^s j_u$ and
let $h_s=\sum\limits_{u=1+g_{s-1}}^{g_s} i_u$.

The $\eta^{pq}$ labeling the leftmost arrow actually stands for a
variety of equivalent maps which factor into instances of the $pq$
interchange.

\[
\xymatrix{
{\cal C}(k)\otimes_p\left(\bigotimes\limits_{s=1}^k{}_q {\cal C}(j_s)\right)\otimes_p
\left(\bigotimes\limits_{t=1}^j{}_q {\cal C}(i_t)\right)
\ar[rr]^>>>>>>>>>>>>{\gamma^{pq} \otimes_p \text{id}}
\ar[dd]_{\text{id} \otimes_p \eta^{pq}}
&&{\cal C}(j)\otimes_p \left(\bigotimes\limits_{t=1}^j{}_q {\cal C}(i_t)\right)
\ar[d]^{\gamma^{pq}}
\\
&&{\cal C}(i)
\\
{\cal C}(k)\otimes_p \left(\bigotimes\limits_{s=1}^k{}_q \left({\cal
C}(j_s)\otimes_p \left(\bigotimes\limits_{u=1}^{j_s}{}_q {\cal
C}(i_{u+g_{s-1}})\right)\right)\right) \ar[rr]_>>>>>>>>>>{\text{id}
\otimes_p(\otimes^k_q\gamma^{pq})} &&{\cal
C}(k)\otimes_p\left(\bigotimes\limits_{s=1}^k{}_q {\cal
C}(h_s)\right) \ar[u]_{\gamma^{pq}} }
\]

\item Respect of units is required just as in the symmetric case.
The following unit diagrams commute for all $m\ge q>p \ge 1$.

\[
\xymatrix{
{\cal C}(k)\otimes_p (\otimes_q^k I)
\ar[d]_{1\otimes_p(\otimes_q^k {\cal J})}
\ar@{=}[r]^{}
&{\cal C}(k)\\
{\cal C}(k)\otimes_p(\otimes_q^k {\cal C}(1))
\ar[ur]_{\gamma^{pq}}
}
\xymatrix{
I\otimes_p {\cal C}(k)
\ar[d]_{{\cal J}\otimes_p 1}
\ar@{=}[r]^{}
&{\cal C}(k)\\
{\cal C}(1)\otimes_p {\cal C}(k)
\ar[ur]^{\gamma^{pq}}
}
\]

\end{enumerate}
\end{definition}

\begin{theorem}
The description of $m$-fold operad in Definition~\ref{expandop} is
equivalent to that given in Definition~\ref{itop}.
\end{theorem}

\begin{proof}
If a collection has an operad composition $\gamma^{q,q+1}$ using
$\otimes_q$ and $\otimes_{q+1}$ then it automatically has an operad
composition for any pair of products $\otimes_p$ and $\otimes_{s}$
for $p< s \le q+1$. This  follows from the fact that for $p < q$
we have natural transformations
$ \eta^{pq}_{AIIB}:A\otimes_p B\to A\otimes_q B$,
as described at the end of Definition~\ref{iterated}.
Thus if we have $\gamma^{q,q+1}$ then we can form
$\gamma^{ps} = \gamma^{q,q+1} \circ (\eta^{pq}
    \circ (1 \otimes_q\eta^{s,q+1})).$
The new $\gamma^{ps} $ is associative based on the old $\gamma$'s
associativity, the naturality of $\eta,$ and the coherence of
$\eta.$ Thus follows our claim that generally operads are preserved
as such by the forgetful functors mentioned in Section 2 and
specifically that an $m$-fold operad is also an $(m-1)$-fold operad.
The converse of this latter statement is not true, as we will see by
counterexample in the final section. It will demonstrate the
existence of $m$-fold operads which are not $(m+1)$-fold operads.
\end{proof}

It is also worth while to expand the definition of the tensor
products of $m$-fold operads that is implicit in their depiction
as monoids in the category of collections in an $n$-fold monoidal
category. Here is the expanded version of the definition:

\begin{definition}\label{tensor}
Let ${\cal C}$,${\cal D}$ be  $m$-fold operads.  For $1 \le i \le (n-m)$
and using a $\otimes'_k$ to denote the product of two $m$-fold
operads, we define that product to be given by:
\[({\cal C}\otimes'_i {\cal D})(j) = {\cal C}(j) \otimes_{i+m} {\cal D}(j).\]
\end{definition}

We note that the new $\gamma$ is in terms of the two old ones, for
$m\ge q>p \ge 1$:
\[
\gamma^{pq}_{{\cal C}\otimes'_i {\cal D}} =
(\gamma^{pq}_{{\cal C}}\otimes_{i+m}\gamma^{pq}_{{\cal D}})\circ \eta^{p,i+m} \circ (1 \otimes_p \eta^{q,i+m})
\]
where the subscripts denote the $n$-fold operad the $\gamma$ belongs
to and the $\eta$'s actually stand for any of the equivalent maps
which factor into instances of the indicated interchange. Note that
this expansion also helps make clear why it is that the
monoidalness, or number of products, of $m$-fold operads must
decrease by the same number $m.$ From the condensed version this is
expected due to the iterated enrichment. From the expanded view this
allows us to define the new composition since in order for the
products of operads to be closed, $\gamma$ for the $i^{th}$ product
utilizes an interchange with superscript $i+m.$ Defined this way $i$
can only be allowed to be as large as $n-m.$ We demonstrate in the
last section in fact a counterexample which shows that the degree of
monoidalness for the category of $m$-fold operads in an $n$-fold
monoidal category is in general no greater than $n-m.$

\section{Examples of iterated monoidal categories}
\label{section:categoryexamples}
\begin{lemma}\label{one}
Given a totally ordered set $S$ with a least element $e \in S$, then
the elements of $S$ with morphisms given by the ordering make up the
objects of a strict monoidal category.
\end{lemma}
The category will also be denoted $S.$  Its morphisms are given by
the ordering means that there is only an arrow $a \to b$ if $a\le
b.$ The product is max and the 2-sided unit is the least element
$e$. We must check that the product is functorial since this defines
monoidal structure on morphisms.  Here it is so since if $a\le b$
and $a'\le b'$ then $\max(a,a')\le\max(b,b').$ Also the identity is
clearly preserved.

\begin{example}\label{uno}
The basic example is the nonnegative integers $\nat$ with their
ordering~$\le.$
\end{example}

\begin{lemma}\label{semi}
Any ordered monoid with its identity element $e$ also its least
element forms the object set of a 2-fold monoidal category.
\end{lemma}

\begin{proof} Morphisms are again given by the ordering. The products
are given by $\max$ and the monoid operation: $a\otimes_1b =
\max(a,b)$ and $a\otimes_2b = ab$.  The shared two-sided unit for
these products is the identity element $e.$ The products are both
strictly associative and functorial since if $a\le b$ and $a'\le b'$
then $aa'\le bb'$ and $\max(a,a')\le\max(b,b').$ The interchange
natural transformations  exist since
$\max(ab,cd)\le\max(a,c)\max(b,d).$ That is because
\begin{align*}
a\le\max(a,c)\quad&\text{and}\quad b\le\max(b,d) \\
\intertext{so}
ab\le\max(a,c)\max(b,d)\quad&\text{and}\quad cd\le\max(a,c)\max(b,d)
\end{align*}

The internal and external unit and associativity conditions of
Definition~\ref{iterated} are all satisfied due to the fact that
there is only one morphism between two objects. More generally,
given any ordered $n$-fold monoidal category with $I$ the least
object we can potentially form an $(n+1)$-fold monoidal category
with  morphisms ordering, and the new $\otimes_1 = \max.$
\end{proof}

\begin{example}\label{nat}
Again we have in mind $\nat$ with its ordering and addition.
\end{example}
Other examples of such  monoids as in Lemma~\ref{semi} are the pure
braids on $n$ strands with only right-handed crossings \cite{rolf}.
Notice that braid composition is a non-symmetric example.  Further
examples are found in the papers on semirings and idempotent
mathematics, such as \cite{LitSob} and its references as well as on
the related concept of tropical geometry, such as \cite{Sturm} and
its references. Semirings that arise in these two areas of study
would require some  translation of the lemmas we have stated thus
far, since the idempotent operation is usually $\min$ and its unit
$\infty.$ Also, since the operation given by addition has unit 0, we
have to consider distinct units. Recall that the additional
requirement is that the two distinct units obey each other's
operations: i.e $I_1\otimes_2I_1 = I_1$ and $I_2\otimes_1I_2 = I_2$.
For example, $\min(0,0)= 0$ and $\infty+\infty=\infty.$

\begin{example}\label{seq}
If $S$ is an ordered set with least element $e$ then by $\seq(S)$ we
denote the infinite sequences $X = \{X_n\}_{n\ge 0}$ of elements of
$S$ for which there exists a natural number  $l(X)$ called the
length such that $k> l(X)$ implies $X_k = e$ and $X_{l(X)}\ne e.$
Under lexicographic ordering $\seq(S)$ is in turn a totally ordered
set with a least element. The latter is the sequence 0 where $0_n =
e$ for all $n.$ We let $l(0) = 0.$ The lexicographic order means
that $A \le B$ if either $A_k = B_k$ for all $k$ or there is a
natural number $n=n_{AB}$ such that $A_k = B_k$ for all $k < n$, and
such that $A_{n} < B_{n}.$

The ordering is easily shown to be reflexive, transitive, and
antisymmetric. See for instance \cite{Schrod} where the case of
lexicographic ordering of $n$-tuples of natural numbers is discussed.
In our case we will need to modify the proof given in that source
by always making comparisons of $\max(l(A),l(B))$-tuples.

As a category $\seq(S)$ is 2-fold monoidal since we can demonstrate
two interchanging products. They are max using the lexicographic
order: $A\otimes_1 B = \max(A,B);$ and concatenation of sequences:
\[ (A\otimes_2 B)_n =
\begin{cases}
A_n,    &n \le l(A) \\
B_{n-l(A)},    &n >   l(A)
\end{cases} \]
Concatenation clearly preserves the ordering, and so
Lemma~\ref{semi} applies.
\end{example}

\begin{example}
Letting $S$ be the set with a single element recovers Example~\ref{nat}
as $\seq(S)$.
\end{example}

\begin{lemma}
If we have an ordered monoid $(M,+)$ as in Lemma~\ref{semi} and
reconsider $\seq(M)$ as in Example~\ref{seq} then we can describe a
3-fold monoidal category $\seq(M,+)$ (with $\seq(M)$ the image of
forgetting the third product of pointwise addition) if and only if
the monoid operation $+$ is such that $0<a<b$ and $c\le d$ imply
both $a+c < b+d$ and $c+a<d+b$ strictly.
\end{lemma}
\begin{proof}
The first two products are again lexicographic max and concatenation
of sequences. The third product $\otimes_3$ is pointwise application
of $+$, $(A\otimes_3B)_n = A_n+B_n$. The last condition that the
monoid operation $+$ strictly respect strict ordering is necessary
to guarantee that the third product both respect the lexicographic
ordering and interchange correctly  with concatenation. To see the
former let sequences $A\le B, C\le D.$ Note that if  $A=B, C=D$ then
$A \otimes_3C = B\otimes_3D.$ Otherwise let $k = \min\{j~|~A_j < B_j
\text{ or } C_j < D_j\}.$ Then $(A\otimes_3C)_k < (B\otimes_3D)_k$
and $(A\otimes_3C)_i = (B\otimes_3D)_i$ for $i<k.$

 To
see that $\otimes_3$ respects the lexicographic ordering only if
addition strictly respects the order, consider a case where $0<a<b$
and $c\le d$ but $a+c=b+d$. Then the sequences $A=(a,a)$, $B=(b,0)$,
$C=(c,0)$, $D=(d,0)$ are such that lexicographically $A<B$ and $C\le
D$ but $A\otimes_3C = (a+c,a) > B\otimes_3D = (b+d,0).$

To see the interchange $(A\otimes_3B)\otimes_2(C\otimes_3D) \le
(A\otimes_2C)\otimes_3(B\otimes_2D) $ notice that we can assume that
$l(A) > l(B).$ Then

\[ \concat(A+B,C+D)\le \concat(A,C)+ \concat(B,D) \]

due to the fact that if $D$ has a first non-zero term, it will be
added to an earlier term of the concatenation of $A$ and $C$ in the
second four-fold product.
\end{proof}

\begin{example}
$\seq(\nat,+)$ plays an important role in  Example~\ref{bee}.
Example interchanges in $\seq(\nat,+)$ are as follow. Let $A =
(0,1,2,0,\dots),~B=(1,1,0,\dots),~C=(2,1,3,0,\dots) \text{ and }
D=(0,2,0,\dots).$ Then
\begin{align*}
(A\otimes_2 B) \otimes_1 (C \otimes_2 D) = (2,1,3,0,2,0,\dots)
\\
(A\otimes_1 C) \otimes_2 (B \otimes_1 D) = (2,1,3,1,1,0,\dots)
\\
(A\otimes_3 B) \otimes_2 (C \otimes_3 D) = (1,2,2,2,3,3,0,\dots)
\\
(A\otimes_2 C) \otimes_3 (B \otimes_2 D) = (1,2,2,4,1,3,0,\dots)
\\
(A\otimes_3 B) \otimes_1 (C \otimes_3 D) = (2,3,3,0,\dots)
\\
(A\otimes_1 C) \otimes_3 (B \otimes_1 D) =(3,2,3,0,\dots)
\end{align*}
\end{example}

\begin{remark}
A non-example is seen if we begin with the monoid of Lemma~\ref{one},
an ordered set with a least element where the product is max. Here
max does not strictly preserve strict ordering, and so pointwise
max does not respect lexicographic ordering. Neither do concatenation
and pointwise max interchange.
\end{remark}

\begin{corollary}\label{duh}
Given any ordered $n$-fold monoidal category $C$ with $I$ the least
object and $\otimes_1$ the max, and whose higher products strictly
respect strict ordering, we can form an $(n+1)$-fold monoidal
category $\seq(C)$.
\end{corollary}
\begin{proof}
The new products of $\seq(C)$ are the lexicographic max, the
concatenation, and the pointwise application of $\otimes_i$ for $i
= 2\dots n.$ The pointwise application of the original products to
the sequences directly inherits the interchange properties. For
instance, if $A, B, C, D \in $ $\seq(C)$ then
$(A_n\otimes_2B_n)\otimes_1(C_n\otimes_2D_n) \le
(A_n\otimes_1C_n)\otimes_2(B_n\otimes_1D_n)$ for all $n$, which
certainly implies that the pointwise 4-fold products are ordered
lexicographically.
\end{proof}

\begin{example}\label{2castles}
Even more symmetrical structure is found in examples with a natural
geometric representation which allows the use of addition in each
product. One such category is that whose objects are Young diagrams,
by which we mean the underlying shapes or diagrams of Young
tableaux. These can be presented by a decreasing sequence of
nonnegative integers in two ways: the sequence that gives the
heights of the columns or the sequence that gives the lengths of the
rows. We let $\otimes_3$ be the product which adds the heights of
columns of two diagrams, $\otimes_2$ adds the length of rows. We
often refer to these as vertical and horizontal stacking
respectively.  If
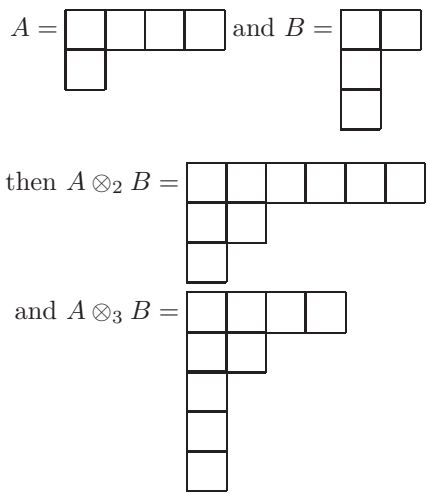
\[
A = \xymatrix@W=.75pc @H=.75pc @R=0pc @C=0pc @*[F-]{~&~&~&~\\~}
\text{ and } B = \xymatrix@W=.75pc @H=.75pc @R=0pc @C=0pc @*[F-]{~&~\\~\\~}
\]
\begin{align*}
\text{then } A\otimes_2 B &= \xymatrix@W=.75pc @H=.75pc @R=0pc @C=0pc @*[F-]{~&~&~&~&~&~\\~&~\\~} \\
\text{ and } A\otimes_3 B &= \xymatrix@W=.75pc @H=.75pc @R=0pc
@C=0pc @*[F-]{~&~&~&~\\~&~\\~\\~\\~}
\end{align*}

 We can take as
morphisms the totally ordered structure of the Young diagrams given
by lexicographic ordering applied to the sequences of column
heights.
Thus we may retain the lexicographic max as $\otimes_1$, and will
refer to the entire category simply as the category of Young
diagrams.

By previous discussion of sequences the Young diagrams with
$\otimes_1$ the lexicographic max and  $\otimes_3$ the piecewise
addition (thought of here as vertical stacking) form a subcategory
of the 3-fold monoidal category called $\seq(\nat,+)$.  To see that
with the additional $\otimes_2$ of horizontal stacking this becomes
a valid 3-fold monoidal category  we look at that operation from
another point of view. Note that the horizontal product of Young
diagrams $A$ and $B$ can be described as a reorganization of all the
columns of both $A$ and $B$ into a new Young diagram made up of
those columns in descending order of height. Rather than (but
equivalent to) the addition of rows, we see horizontal stacking as
the concatenation of monotone decreasing sequences (of columns)
followed by sorting greatest to least.  We call this operation
\emph{merging}.
\end{example}

\begin{lemma}\label{sort}
Let $(S, \le, +)$ be an ordered monoid and consider the sequences
$\seq(S,+)$ with lexicographic ordering, piecewise addition $+$ and
the function of sorting denoted by
\[ s\colon\seq(S,+)\to \seq(S,+) \]
Then the triangle inequality holds for two sequences:
$s(A+B)\le s(A)+s(B).$
\end{lemma}

\begin{proof}
Consider $s(A+B)$, where we start with the two sequences and add
them piecewise before sorting. We can metamorphose this into
$s(A)+s(B)$ in stages by using an algorithm to sort $A$ and $B$.
Note that if $A$ and $B$ are already sorted, the inequality becomes
an equality.  For our algorithm we choose parallel bubble sorting.
This consists of a series of passes through the sequences comparing
$A_n$ and $A_{n+1}$ and  comparing $B_n$ and $B_{n+1}$
simultaneously. If the two elements of a given sequence are not
already in strictly decreasing order we switch their places.  We
claim that switching consecutive sequence elements into order always
results in a lexicographically larger sequence after adding
piecewise and sorting. If both the elements of $A$ and of $B$ are
switched, or if neither, then the result is unaltered. Therefore
without loss of generality we assume that $A_n < A_{n+1}$ and that
$B_{n+1} < B_n.$ Then we compare the original result of sorting
after adding and the same but after the switching of $A_n$ and
$A_{n+1}.$ It is simplest to note that the new result includes
$A_{n+1} + B_n,$ which is larger than both $A_n + B_n$ and $A_{n+1}
+ B_{n+1}.$ So after adding and sorting the new result is indeed
larger lexicographically.  Thus since each move of the parallel
bubble sort results in a larger expression after first adding and
then sorting, and after all the moves the result of adding and then
sorting the now pre-sorted sequences is the same as first sorting
then adding, the triangle inequality follows.
\end{proof}

\begin{theorem}
The category of
Young diagrams forms a 3-fold monoidal
category.
\end{theorem}
\begin{proof}
  The
products on Young diagrams are $\otimes_1 = $ lexicographic max,
$\otimes_2 = $ horizontal stacking and $\otimes_3 =$ vertical
stacking. We need to check first  that horizontal stacking, or
merging, is functorial with respect to morphisms (defined as the
$\le$ relations of the lexicographic ordering.)  The cases where
$A=B$ or $C=D$ are easy. For example let $A_k = B_k$ for all $k$
and $C_k = D_k$ for all $k < n_{CD},$ where $n_{CD}$ is as defined
in Example~\ref{seq}. Thus the columns from the copies of, for
instance $A$ in $A\otimes_1 C$ and $A\otimes_1 D$ fall into the
same final spot under the sortings right up to the critical location,
so if $C \le D$, then $A\otimes_1 C \le A\otimes_1 D.$ Similarly,
it is clear that $A \le B$ implies $(A \otimes_1 D) \le (B \otimes_1 D).$
Hence if $A\le B$ and $C\le D$, then
$A\otimes_1 C \le A\otimes_1 D \le B\otimes_1 D$
which by transitivity gives us our desired property.

Next we check that our interchange transformations will always
exist. $\eta^{1j}$ exists by the proof of Lemma~\ref{semi} for
$j=2,3$ since the higher products both respect morphisms(ordering)
and are thus ordered monoid operations. We need to check for existence
of $\eta^{23},$ i.e. we need to show that
$(A \otimes_3 B)\otimes_2(C \otimes_3 D)
    \le (A \otimes_2 C)\otimes_3 (B \otimes_2 D).$
This inequality follows from Lemma~\ref{sort} on the triangle
inequality for sorting.  To prove the new inequality we consider
the special case of two sequences formed by letting $A'$ be $A$
followed by $C$ and letting $B'$ be $B$ followed by $D$. By ``followed
by'' we mean padded by zeroes so that $l(A') = \max(l(A),l(B)) + l(C) $
and $l(B') = \max(l(A),l(B)) + l(D).$ Thus piecewise addition
of $A'$ and $B'$ results in piecewise addition of $A$ and $B$, and
respectively $C$ and $D.$ Then to our new sequences $A'$ and $B'$
we apply the result of Lemma~\ref{sort} and have our desired result.
\end{proof}

Here is an example of the inequality we have just shown to always hold.
Let four Young diagrams be as follow:
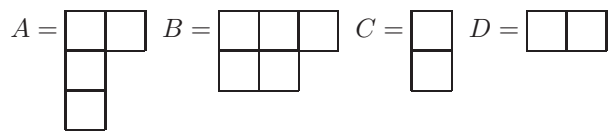
\[
A = \xymatrix@W=.75pc @H=.75pc @R=0pc @C=0pc @*[F-]{~&~\\~\\~}~~
B = \xymatrix@W=.75pc @H=.75pc @R=0pc @C=0pc @*[F-]{~&~&~\\~&~}~~
C = \xymatrix@W=.75pc @H=.75pc @R=0pc @C=0pc @*[F-]{~\\~}~~
D = \xymatrix@W=.75pc @H=.75pc @R=0pc @C=0pc @*[F-]{~&~}
\]
Then the fact that
$(A \otimes_3 B)\otimes_2(C \otimes_3 D)
    \le (A \otimes_2 C)\otimes_3 (B \otimes_2 D)$
appears as follows:
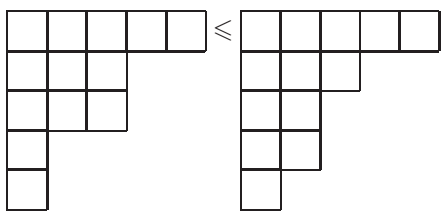
\[
\xymatrix@W=.75pc @H=.75pc @R=0pc @C=0pc @*[F-]{~&~&~&~&~\\~&~&~\\~&~&~\\~\\~} \le \xymatrix@W=.75pc @H=.75pc @R=0pc @C=0pc @*[F-]{~&~&~&~&~\\~&~&~\\~&~\\~&~\\~}
\]

\begin{remark}\label{hpre}
Alternatively we can create a category equivalent to the non-negative
integers in Example~\ref{uno} by pre-ordering the Young diagrams
by height.  Here the height $h(A)$ of the Young diagram is the
number of boxes in its leftmost column, and  we say $A \le B$ if
$h(A)\le h(B)$.  Two Young diagrams with the same height are
isomorphic objects, and the one-column stacks form both a full
subcategory and a  skeleton of the height preordered category.
Everything works as for the previous example of natural numbers
since $h(A\otimes_2 B) = h(A) +h(B)$ and $h(A\otimes_1 B) = \max(h(A),
h(B)).$ There is also a max product; the new max with respect to
the height preordering is defined as
\[ \max(A,B) = \begin{cases}
    A,& \text{if $B\le A$} \\
    B,& \text{otherwise.}
\end{cases} \]
In the height preordered category this latter product is equivalent
to horizontal stacking,~$\otimes_1$.
\end{remark}

\begin{remark}
Notice that we can start with any totally ordered monoids $\{M,\le,
+\}$ such that the identity 0 is less than any other element and
such that $0<a<b$ and $c\le d$ implies both $a+c < b+d$ and $c+a<d+b$
for all $a,b,c \in G.$ We create a 3-fold monoidal category
$\modseq(M,+)$ with objects monotone decreasing finitely non-zero
sequences of elements of $M$ and morphisms given by the lexicographic
ordering.  The products are as described for the category of Young
diagrams $\modseq(\nat,+)$ in the previous example.  The common
unit is the zero sequence.   The proofs we have given in the previous
example for $M=\nat$  are all still valid.
\end{remark}

By Corollary~\ref{duh} we can also consider 4-fold monoidal categories
such as $\seq(\modseq(M))$ and other combinations of Seq and ModSeq.
For instance if $\modseq(\nat,+)$ is our category of Young diagrams
then $\modseq(\modseq(\nat,+))$ has objects monotone decreasing
sequences of Young diagrams, which we can visualize along the
$z$-axis. Here the lexicographic-lexicographic max is  $\otimes_1$,
lexicographic merging is $\otimes_2$, pointwise merging (pointwise
horizontal or $y$-axis stacking) is $\otimes_3$ and pointwise-pointwise
addition (pointwise $x$-axis stacking) is $\otimes_4$. For example,
if:
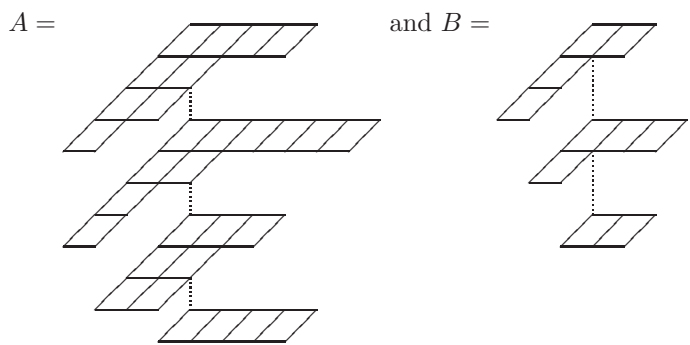
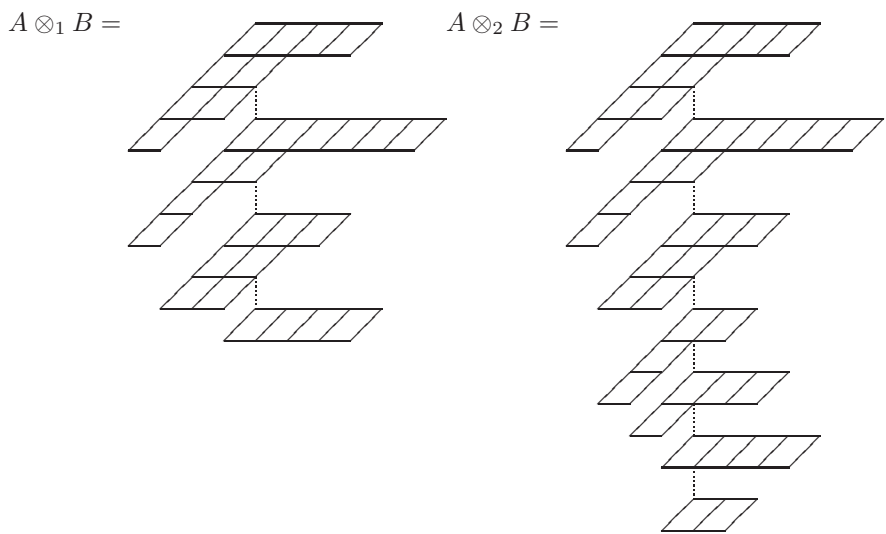
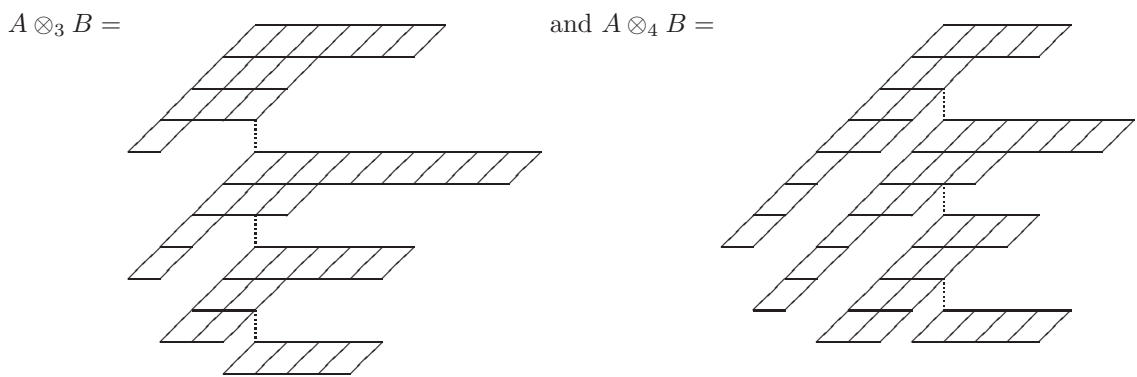
\[
A=\xymatrix @W=1.0pc @H=1.0pc @R=1.0pc @C=1.0pc
{*=0{}&*=0{}&*=0{}&*=0{}&*=0{}\ar@{-}[r]\ar@{-}[dl]&*=0{}\ar@{-}[r]\ar@{-}[dl]&*=0{}\ar@{-}[r]\ar@{-}[dl]&*=0{}\ar@{-}[r]\ar@{-}[dl]&*=0{}\ar@{-}[dl]
\\*=0{}&*=0{}&*=0{}&*=0{}\ar@{-}[r]\ar@{-}[dl]&*=0{}\ar@{-}[r]\ar@{-}[dl]&*=0{}\ar@{-}[r]\ar@{-}[dl]&*=0{}\ar@{-}[r]&*=0{}
\\*=0{}&*=0{}&*=0{}\ar@{-}[r]\ar@{-}[dl]&*=0{}\ar@{-}[r]\ar@{-}[dl]&*=0{}\ar@{-}[dl]\ar@{..}[d]
\\*=0{}&*=0{}\ar@{-}[r]\ar@{-}[dl]&*=0{}\ar@{-}[r]\ar@{-}[dl]&*=0{}&*=0{}\ar@{-}[r]\ar@{-}[dl]&*=0{}\ar@{-}[r]\ar@{-}[dl]&*=0{}\ar@{-}[r]\ar@{-}[dl]&*=0{}\ar@{-}[r]\ar@{-}[dl]&*=0{}\ar@{-}[r]\ar@{-}[dl]&*=0{}\ar@{-}[r]\ar@{-}[dl]&*=0{}\ar@{-}[dl]
\\*=0{}\ar@{-}[r]&*=0{}&*=0{}&*=0{}\ar@{-}[r]\ar@{-}[dl]&*=0{}\ar@{-}[r]\ar@{-}[dl]&*=0{}\ar@{-}[r]\ar@{-}[dl]&*=0{}\ar@{-}[r]&*=0{}\ar@{-}[r]&*=0{}\ar@{-}[r]&*=0{}
\\*=0{}&*=0{}&*=0{}\ar@{-}[r]\ar@{-}[dl]&*=0{}\ar@{-}[r]\ar@{-}[dl]&*=0{}\ar@{..}[d]
\\*=0{}&*=0{}\ar@{-}[r]\ar@{-}[dl]&*=0{}\ar@{-}[dl]&*=0{}&*=0{}\ar@{-}[r]\ar@{-}[dl]&*=0{}\ar@{-}[r]\ar@{-}[dl]&*=0{}\ar@{-}[r]\ar@{-}[dl]&*=0{}\ar@{-}[dl]
\\*=0{}\ar@{-}[r]&*=0{}&*=0{}&*=0{}\ar@{-}[r]\ar@{-}[dl]&*=0{}\ar@{-}[r]\ar@{-}[dl]&*=0{}\ar@{-}[r]\ar@{-}[dl]&*=0{}
\\*=0{}&*=0{}&*=0{}\ar@{-}[r]\ar@{-}[dl]&*=0{}\ar@{-}[r]\ar@{-}[dl]&*=0{}\ar@{-}[dl]\ar@{..}[d]
\\*=0{}&*=0{}\ar@{-}[r]&*=0{}\ar@{-}[r]&*=0{}&*=0{}\ar@{-}[r]\ar@{-}[dl]&*=0{}\ar@{-}[r]\ar@{-}[dl]&*=0{}\ar@{-}[r]\ar@{-}[dl]&*=0{}\ar@{-}[r]\ar@{-}[dl]&*=0{}\ar@{-}[dl]
\\*=0{}&*=0{}&*=0{}&*=0{}\ar@{-}[r]&*=0{}\ar@{-}[r]&*=0{}\ar@{-}[r]&*=0{}\ar@{-}[r]&*=0{}
\\*=0{}}
\text{ and }B=\xymatrix @W=1.0pc @H=1.0pc @R=1.0pc @C=1.0pc
{*=0{}&*=0{}&*=0{}&*=0{}\ar@{-}[r]\ar@{-}[dl]&*=0{}\ar@{-}[r]\ar@{-}[dl]&*=0{}\ar@{-}[dl]
\\*=0{}&*=0{}&*=0{}\ar@{-}[r]\ar@{-}[dl]&*=0{}\ar@{-}[r]\ar@{-}[dl]\ar@{..}[dd]&*=0{}
\\*=0{}&*=0{}\ar@{-}[r]\ar@{-}[dl]&*=0{}\ar@{-}[dl]
\\*=0{}\ar@{-}[r]&*=0{}&*=0{}&*=0{}\ar@{-}[r]\ar@{-}[dl]&*=0{}\ar@{-}[r]\ar@{-}[dl]&*=0{}\ar@{-}[r]\ar@{-}[dl]&*=0{}\ar@{-}[dl]
\\*=0{}&*=0{}&*=0{}\ar@{-}[r]\ar@{-}[dl]&*=0{}\ar@{-}[r]\ar@{-}[dl]\ar@{..}[dd]&*=0{}\ar@{-}[r]&*=0{}
\\*=0{}&*=0{}\ar@{-}[r]&*=0{}
\\
*=0{}&*=0{}&*=0{}&*=0{}\ar@{-}[r]\ar@{-}[dl]&*=0{}\ar@{-}[r]\ar@{-}[dl]&*=0{}\ar@{-}[dl]
\\*=0{}&*=0{}&*=0{}\ar@{-}[r]&*=0{}\ar@{-}[r]&*=0{}
\\*=0{}}
\]
then
\[
A\otimes_1 B =\xymatrix @W=1.0pc @H=1.0pc @R=1.0pc @C=1.0pc
{*=0{}&*=0{}&*=0{}&*=0{}&*=0{}\ar@{-}[r]\ar@{-}[dl]&*=0{}\ar@{-}[r]\ar@{-}[dl]&*=0{}\ar@{-}[r]\ar@{-}[dl]&*=0{}\ar@{-}[r]\ar@{-}[dl]&*=0{}\ar@{-}[dl]
\\*=0{}&*=0{}&*=0{}&*=0{}\ar@{-}[r]\ar@{-}[dl]&*=0{}\ar@{-}[r]\ar@{-}[dl]&*=0{}\ar@{-}[r]\ar@{-}[dl]&*=0{}\ar@{-}[r]&*=0{}
\\*=0{}&*=0{}&*=0{}\ar@{-}[r]\ar@{-}[dl]&*=0{}\ar@{-}[r]\ar@{-}[dl]&*=0{}\ar@{-}[dl]\ar@{..}[d]
\\*=0{}&*=0{}\ar@{-}[r]\ar@{-}[dl]&*=0{}\ar@{-}[r]\ar@{-}[dl]&*=0{}&*=0{}\ar@{-}[r]\ar@{-}[dl]&*=0{}\ar@{-}[r]\ar@{-}[dl]&*=0{}\ar@{-}[r]\ar@{-}[dl]&*=0{}\ar@{-}[r]\ar@{-}[dl]&*=0{}\ar@{-}[r]\ar@{-}[dl]&*=0{}\ar@{-}[r]\ar@{-}[dl]&*=0{}\ar@{-}[dl]
\\*=0{}\ar@{-}[r]&*=0{}&*=0{}&*=0{}\ar@{-}[r]\ar@{-}[dl]&*=0{}\ar@{-}[r]\ar@{-}[dl]&*=0{}\ar@{-}[r]\ar@{-}[dl]&*=0{}\ar@{-}[r]&*=0{}\ar@{-}[r]&*=0{}\ar@{-}[r]&*=0{}
\\*=0{}&*=0{}&*=0{}\ar@{-}[r]\ar@{-}[dl]&*=0{}\ar@{-}[r]\ar@{-}[dl]&*=0{}\ar@{..}[d]
\\*=0{}&*=0{}\ar@{-}[r]\ar@{-}[dl]&*=0{}\ar@{-}[dl]&*=0{}&*=0{}\ar@{-}[r]\ar@{-}[dl]&*=0{}\ar@{-}[r]\ar@{-}[dl]&*=0{}\ar@{-}[r]\ar@{-}[dl]&*=0{}\ar@{-}[dl]
\\*=0{}\ar@{-}[r]&*=0{}&*=0{}&*=0{}\ar@{-}[r]\ar@{-}[dl]&*=0{}\ar@{-}[r]\ar@{-}[dl]&*=0{}\ar@{-}[r]\ar@{-}[dl]&*=0{}
\\*=0{}&*=0{}&*=0{}\ar@{-}[r]\ar@{-}[dl]&*=0{}\ar@{-}[r]\ar@{-}[dl]&*=0{}\ar@{-}[dl]\ar@{..}[d]
\\*=0{}&*=0{}\ar@{-}[r]&*=0{}\ar@{-}[r]&*=0{}&*=0{}\ar@{-}[r]\ar@{-}[dl]&*=0{}\ar@{-}[r]\ar@{-}[dl]&*=0{}\ar@{-}[r]\ar@{-}[dl]&*=0{}\ar@{-}[r]\ar@{-}[dl]&*=0{}\ar@{-}[dl]
\\*=0{}&*=0{}&*=0{}&*=0{}\ar@{-}[r]&*=0{}\ar@{-}[r]&*=0{}\ar@{-}[r]&*=0{}\ar@{-}[r]&*=0{}
\\*=0{}}
A \otimes_2 B=\xymatrix @W=1.0pc @H=1.0pc @R=1.0pc @C=1.0pc
{*=0{}&*=0{}&*=0{}&*=0{}&*=0{}\ar@{-}[r]\ar@{-}[dl]&*=0{}\ar@{-}[r]\ar@{-}[dl]&*=0{}\ar@{-}[r]\ar@{-}[dl]&*=0{}\ar@{-}[r]\ar@{-}[dl]&*=0{}\ar@{-}[dl]
\\*=0{}&*=0{}&*=0{}&*=0{}\ar@{-}[r]\ar@{-}[dl]&*=0{}\ar@{-}[r]\ar@{-}[dl]&*=0{}\ar@{-}[r]\ar@{-}[dl]&*=0{}\ar@{-}[r]&*=0{}
\\*=0{}&*=0{}&*=0{}\ar@{-}[r]\ar@{-}[dl]&*=0{}\ar@{-}[r]\ar@{-}[dl]&*=0{}\ar@{-}[dl]\ar@{..}[d]
\\*=0{}&*=0{}\ar@{-}[r]\ar@{-}[dl]&*=0{}\ar@{-}[r]\ar@{-}[dl]&*=0{}&*=0{}\ar@{-}[r]\ar@{-}[dl]&*=0{}\ar@{-}[r]\ar@{-}[dl]&*=0{}\ar@{-}[r]\ar@{-}[dl]&*=0{}\ar@{-}[r]\ar@{-}[dl]&*=0{}\ar@{-}[r]\ar@{-}[dl]&*=0{}\ar@{-}[r]\ar@{-}[dl]&*=0{}\ar@{-}[dl]
\\*=0{}\ar@{-}[r]&*=0{}&*=0{}&*=0{}\ar@{-}[r]\ar@{-}[dl]&*=0{}\ar@{-}[r]\ar@{-}[dl]&*=0{}\ar@{-}[r]\ar@{-}[dl]&*=0{}\ar@{-}[r]&*=0{}\ar@{-}[r]&*=0{}\ar@{-}[r]&*=0{}
\\*=0{}&*=0{}&*=0{}\ar@{-}[r]\ar@{-}[dl]&*=0{}\ar@{-}[r]\ar@{-}[dl]&*=0{}\ar@{..}[d]
\\*=0{}&*=0{}\ar@{-}[r]\ar@{-}[dl]&*=0{}\ar@{-}[dl]&*=0{}&*=0{}\ar@{-}[r]\ar@{-}[dl]&*=0{}\ar@{-}[r]\ar@{-}[dl]&*=0{}\ar@{-}[r]\ar@{-}[dl]&*=0{}\ar@{-}[dl]
\\*=0{}\ar@{-}[r]&*=0{}&*=0{}&*=0{}\ar@{-}[r]\ar@{-}[dl]&*=0{}\ar@{-}[r]\ar@{-}[dl]&*=0{}\ar@{-}[r]\ar@{-}[dl]&*=0{}
\\*=0{}&*=0{}&*=0{}\ar@{-}[r]\ar@{-}[dl]&*=0{}\ar@{-}[r]\ar@{-}[dl]&*=0{}\ar@{-}[dl]\ar@{..}[d]
\\*=0{}&*=0{}\ar@{-}[r]&*=0{}\ar@{-}[r]&*=0{}&*=0{}\ar@{-}[r]\ar@{-}[dl]&*=0{}\ar@{-}[r]\ar@{-}[dl]&*=0{}\ar@{-}[dl]
\\*=0{}&*=0{}&*=0{}&*=0{}\ar@{-}[r]\ar@{-}[dl]&*=0{}\ar@{-}[r]\ar@{-}[dl]\ar@{..}[d]&*=0{}
\\*=0{}&*=0{}&*=0{}\ar@{-}[r]\ar@{-}[dl]&*=0{}\ar@{-}[dl]&*=0{}\ar@{-}[r]\ar@{-}[dl]&*=0{}\ar@{-}[r]\ar@{-}[dl]&*=0{}\ar@{-}[r]\ar@{-}[dl]&*=0{}\ar@{-}[dl]
\\*=0{}&*=0{}\ar@{-}[r]&*=0{}&*=0{}\ar@{-}[r]\ar@{-}[dl]&*=0{}\ar@{-}[r]\ar@{-}[dl]\ar@{..}[d]&*=0{}\ar@{-}[r]&*=0{}
\\*=0{}&*=0{}&*=0{}\ar@{-}[r]&*=0{}&*=0{}\ar@{-}[r]\ar@{-}[dl]&*=0{}\ar@{-}[r]\ar@{-}[dl]&*=0{}\ar@{-}[r]\ar@{-}[dl]&*=0{}\ar@{-}[r]\ar@{-}[dl]&*=0{}\ar@{-}[dl]
\\*=0{}&*=0{}&*=0{}&*=0{}\ar@{-}[r]&*=0{}\ar@{-}[r]\ar@{..}[d]&*=0{}\ar@{-}[r]&*=0{}\ar@{-}[r]&*=0{}
\\
*=0{}&*=0{}&*=0{}&*=0{}&*=0{}\ar@{-}[r]\ar@{-}[dl]&*=0{}\ar@{-}[r]\ar@{-}[dl]&*=0{}\ar@{-}[dl]
\\*=0{}&*=0{}&*=0{}&*=0{}\ar@{-}[r]&*=0{}\ar@{-}[r]&*=0{}
\\*=0{}}
\]
\[
A \otimes_3 B=\xymatrix @W=1.0pc @H=1.0pc @R=1.0pc @C=1.0pc
{*=0{}&*=0{}&*=0{}&*=0{}&*=0{}\ar@{-}[r]\ar@{-}[dl]&*=0{}\ar@{-}[r]\ar@{-}[dl]&*=0{}\ar@{-}[r]\ar@{-}[dl]&*=0{}\ar@{-}[r]\ar@{-}[dl]&*=0{}\ar@{-}[r]\ar@{-}[dl]&*=0{}\ar@{-}[r]\ar@{-}[dl]&*=0{}\ar@{-}[dl]
\\*=0{}&*=0{}&*=0{}&*=0{}\ar@{-}[r]\ar@{-}[dl]&*=0{}\ar@{-}[r]\ar@{-}[dl]&*=0{}\ar@{-}[r]\ar@{-}[dl]&*=0{}\ar@{-}[r]\ar@{-}[dl]&*=0{}\ar@{-}[r]&*=0{}\ar@{-}[r]&*=0{}
\\*=0{}&*=0{}&*=0{}\ar@{-}[r]\ar@{-}[dl]&*=0{}\ar@{-}[r]\ar@{-}[dl]&*=0{}\ar@{-}[r]\ar@{-}[dl]&*=0{}\ar@{-}[dl]
\\*=0{}&*=0{}\ar@{-}[r]\ar@{-}[dl]&*=0{}\ar@{-}[r]\ar@{-}[dl]&*=0{}\ar@{-}[r]&*=0{}\ar@{..}[d]
\\*=0{}\ar@{-}[r]&*=0{}&*=0{}&*=0{}&*=0{}\ar@{-}[r]\ar@{-}[dl]&*=0{}\ar@{-}[r]\ar@{-}[dl]&*=0{}\ar@{-}[r]\ar@{-}[dl]&*=0{}\ar@{-}[r]\ar@{-}[dl]&*=0{}\ar@{-}[r]\ar@{-}[dl]&*=0{}\ar@{-}[r]\ar@{-}[dl]&*=0{}\ar@{-}[r]\ar@{-}[dl]&*=0{}\ar@{-}[r]\ar@{-}[dl]&*=0{}\ar@{-}[r]\ar@{-}[dl]&*=0{}\ar@{-}[dl]
\\*=0{}&*=0{}&*=0{}&*=0{}\ar@{-}[r]\ar@{-}[dl]&*=0{}\ar@{-}[r]\ar@{-}[dl]&*=0{}\ar@{-}[r]\ar@{-}[dl]&*=0{}\ar@{-}[r]\ar@{-}[dl]&*=0{}\ar@{-}[r]&*=0{}\ar@{-}[r]&*=0{}\ar@{-}[r]&*=0{}\ar@{-}[r]&*=0{}\ar@{-}[r]&*=0{}
\\*=0{}&*=0{}&*=0{}\ar@{-}[r]\ar@{-}[dl]&*=0{}\ar@{-}[r]\ar@{-}[dl]&*=0{}\ar@{-}[r]\ar@{..}[d]&*=0{}
\\*=0{}&*=0{}\ar@{-}[r]\ar@{-}[dl]&*=0{}\ar@{-}[dl]&*=0{}&*=0{}\ar@{-}[r]\ar@{-}[dl]&*=0{}\ar@{-}[r]\ar@{-}[dl]&*=0{}\ar@{-}[r]\ar@{-}[dl]&*=0{}\ar@{-}[r]\ar@{-}[dl]&*=0{}\ar@{-}[r]\ar@{-}[dl]&*=0{}\ar@{-}[dl]
\\*=0{}\ar@{-}[r]&*=0{}&*=0{}&*=0{}\ar@{-}[r]\ar@{-}[dl]&*=0{}\ar@{-}[r]\ar@{-}[dl]&*=0{}\ar@{-}[r]\ar@{-}[dl]&*=0{}\ar@{-}[r]&*=0{}\ar@{-}[r]&*=0{}
\\*=0{}&*=0{}&*=0{}\ar@{-}[r]\ar@{-}[dl]&*=0{}\ar@{-}[r]\ar@{-}[dl]&*=0{}\ar@{-}[dl]\ar@{..}[d]
\\*=0{}&*=0{}\ar@{-}[r]&*=0{}\ar@{-}[r]&*=0{}&*=0{}\ar@{-}[r]\ar@{-}[dl]&*=0{}\ar@{-}[r]\ar@{-}[dl]&*=0{}\ar@{-}[r]\ar@{-}[dl]&*=0{}\ar@{-}[r]\ar@{-}[dl]&*=0{}\ar@{-}[dl]
\\*=0{}&*=0{}&*=0{}&*=0{}\ar@{-}[r]&*=0{}\ar@{-}[r]&*=0{}\ar@{-}[r]&*=0{}\ar@{-}[r]&*=0{}
\\*=0{}}
\text{ and }
A \otimes_4 B=\xymatrix @W=1.0pc @H=1.0pc @R=1.0pc @C=1.0pc
{*=0{}&*=0{}&*=0{}&*=0{}&*=0{}&*=0{}&*=0{}&*=0{}\ar@{-}[r]\ar@{-}[dl]&*=0{}\ar@{-}[r]\ar@{-}[dl]&*=0{}\ar@{-}[r]\ar@{-}[dl]&*=0{}\ar@{-}[r]\ar@{-}[dl]&*=0{}\ar@{-}[dl]
\\*=0{}&*=0{}&*=0{}&*=0{}&*=0{}&*=0{}&*=0{}\ar@{-}[r]\ar@{-}[dl]&*=0{}\ar@{-}[r]\ar@{-}[dl]&*=0{}\ar@{-}[r]\ar@{-}[dl]&*=0{}\ar@{-}[r]&*=0{}
\\*=0{}&*=0{}&*=0{}&*=0{}&*=0{}&*=0{}\ar@{-}[r]\ar@{-}[dl]&*=0{}\ar@{-}[r]\ar@{-}[dl]&*=0{}\ar@{-}[dl]\ar@{..}[d]
\\*=0{}&*=0{}&*=0{}&*=0{}&*=0{}\ar@{-}[r]\ar@{-}[dl]&*=0{}\ar@{-}[r]\ar@{-}[dl]&*=0{}\ar@{-}[dl]&*=0{}\ar@{-}[r]\ar@{-}[dl]&*=0{}\ar@{-}[r]\ar@{-}[dl]&*=0{}\ar@{-}[r]\ar@{-}[dl]&*=0{}\ar@{-}[r]\ar@{-}[dl]&*=0{}\ar@{-}[r]\ar@{-}[dl]&*=0{}\ar@{-}[r]\ar@{-}[dl]&*=0{}\ar@{-}[dl]
\\*=0{}&*=0{}&*=0{}&*=0{}\ar@{-}[r]\ar@{-}[dl]&*=0{}\ar@{-}[r]\ar@{-}[dl]&*=0{}&*=0{}\ar@{-}[r]\ar@{-}[dl]&*=0{}\ar@{-}[r]\ar@{-}[dl]&*=0{}\ar@{-}[r]\ar@{-}[dl]&*=0{}\ar@{-}[r]\ar@{-}[dl]&*=0{}\ar@{-}[r]&*=0{}\ar@{-}[r]&*=0{}
\\*=0{}&*=0{}&*=0{}\ar@{-}[r]\ar@{-}[dl]&*=0{}\ar@{-}[dl]&*=0{}&*=0{}\ar@{-}[r]\ar@{-}[dl]&*=0{}\ar@{-}[r]\ar@{-}[dl]&*=0{}\ar@{-}[r]\ar@{-}[dl]\ar@{..}[d]&*=0{}
\\*=0{}&*=0{}\ar@{-}[r]\ar@{-}[dl]&*=0{}\ar@{-}[dl]&*=0{}&*=0{}\ar@{-}[r]\ar@{-}[dl]&*=0{}\ar@{-}[r]\ar@{-}[dl]&*=0{}&*=0{}\ar@{-}[r]\ar@{-}[dl]&*=0{}\ar@{-}[r]\ar@{-}[dl]&*=0{}\ar@{-}[r]\ar@{-}[dl]&*=0{}\ar@{-}[dl]
\\*=0{}\ar@{-}[r]&*=0{}&*=0{}&*=0{}\ar@{-}[r]\ar@{-}[dl]&*=0{}\ar@{-}[dl]&*=0{}&*=0{}\ar@{-}[r]\ar@{-}[dl]&*=0{}\ar@{-}[r]\ar@{-}[dl]&*=0{}\ar@{-}[r]\ar@{-}[dl]&*=0{}
\\*=0{}&*=0{}&*=0{}\ar@{-}[r]\ar@{-}[dl]&*=0{}\ar@{-}[dl]&*=0{}&*=0{}\ar@{-}[r]\ar@{-}[dl]&*=0{}\ar@{-}[r]\ar@{-}[dl]&*=0{}\ar@{-}[dl]\ar@{..}[d]
\\*=0{}&*=0{}\ar@{-}[r]&*=0{}&*=0{}&*=0{}\ar@{-}[r]\ar@{-}[dl]&*=0{}\ar@{-}[r]\ar@{-}[dl]&*=0{}\ar@{-}[dl]&*=0{}\ar@{-}[r]\ar@{-}[dl]&*=0{}\ar@{-}[r]\ar@{-}[dl]&*=0{}\ar@{-}[r]\ar@{-}[dl]&*=0{}\ar@{-}[r]\ar@{-}[dl]&*=0{}\ar@{-}[dl]
\\*=0{}&*=0{}&*=0{}&*=0{}\ar@{-}[r]&*=0{}\ar@{-}[r]&*=0{}&*=0{}\ar@{-}[r]&*=0{}\ar@{-}[r]&*=0{}\ar@{-}[r]&*=0{}\ar@{-}[r]&*=0{}
\\*=0{}}
\]

\begin{example}\label{3-castles}
It might be nice to retain the geometric picture of the products of
Young diagrams in terms of vertical and horizontal stacking, and
stacking in other directions as dimension increases. This is not
found in the just illustrated category, which relies on the merging
viewpoint.  The ``diagram stacking'' point of view is restored if we
restrict to 3-d  Young diagrams. We can represent these objects as
infinite matrices with finitely many nonzero natural number entries,
and with monotone decreasing columns and rows.  We  require that
${A_n}_k$ be decreasing in $n$ for constant $k$, and decreasing in
$k$ for constant $n.$ We choose the sequence of rows to represent
the sequence of sequences, i.e. each row represents a Young diagram
which we draw as being parallel to the $xy$ plane.  This choice is
important because it determines the total ordering of matrices and
thus the morphisms of the category.  Thus $y$-axis stacking is
horizontal concatenation (disregarding  trailing zeroes) of matrices
followed by sorting the new longer rows (row merging).  $x$-axis
stacking is addition of matrices.  Now we define $z$-axis stacking
as vertical concatenation  of matrices followed by sorting the new
long columns (column merging).

\end{example}
Here is a visual example of the three new products, beginning with
$z$-axis stacking, labeled $\otimes_1$: if
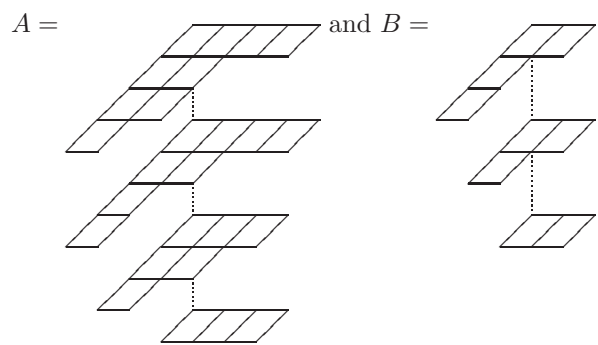
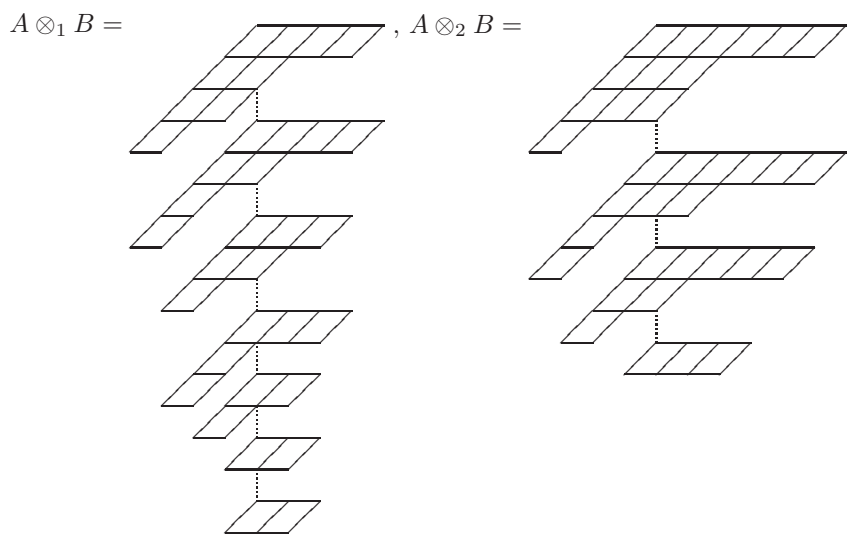
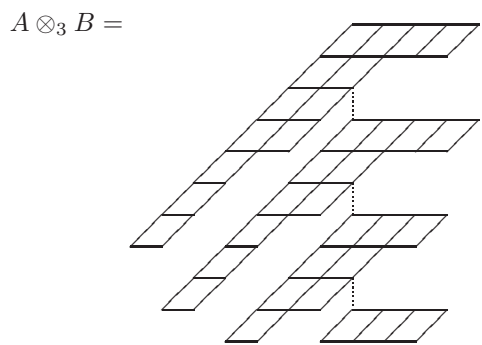
\[
A=\xymatrix @W=1.0pc @H=1.0pc @R=1.0pc @C=1.0pc
{*=0{}&*=0{}&*=0{}&*=0{}&*=0{}\ar@{-}[r]\ar@{-}[dl]&*=0{}\ar@{-}[r]\ar@{-}[dl]&*=0{}\ar@{-}[r]\ar@{-}[dl]&*=0{}\ar@{-}[r]\ar@{-}[dl]&*=0{}\ar@{-}[dl]
\\*=0{}&*=0{}&*=0{}&*=0{}\ar@{-}[r]\ar@{-}[dl]&*=0{}\ar@{-}[r]\ar@{-}[dl]&*=0{}\ar@{-}[r]\ar@{-}[dl]&*=0{}\ar@{-}[r]&*=0{}
\\*=0{}&*=0{}&*=0{}\ar@{-}[r]\ar@{-}[dl]&*=0{}\ar@{-}[r]\ar@{-}[dl]&*=0{}\ar@{-}[dl]\ar@{..}[d]
\\*=0{}&*=0{}\ar@{-}[r]\ar@{-}[dl]&*=0{}\ar@{-}[r]\ar@{-}[dl]&*=0{}  &*=0{}\ar@{-}[r]\ar@{-}[dl]&*=0{}\ar@{-}[r]\ar@{-}[dl]&*=0{}\ar@{-}[r]\ar@{-}[dl]&*=0{}\ar@{-}[r]\ar@{-}[dl]&*=0{}\ar@{-}[dl]
\\*=0{}\ar@{-}[r]&*=0{} &*=0{}&*=0{}\ar@{-}[r]\ar@{-}[dl]&*=0{}\ar@{-}[r]\ar@{-}[dl]&*=0{}\ar@{-}[r]\ar@{-}[dl]&*=0{}\ar@{-}[r]&*=0{}
\\*=0{}&*=0{}&*=0{}\ar@{-}[r]\ar@{-}[dl]&*=0{}\ar@{-}[r]\ar@{-}[dl]&*=0{}\ar@{..}[d]
\\*=0{}&*=0{}\ar@{-}[r]\ar@{-}[dl]&*=0{}\ar@{-}[dl]&*=0{}&*=0{}\ar@{-}[r]\ar@{-}[dl]&*=0{}\ar@{-}[r]\ar@{-}[dl]&*=0{}\ar@{-}[r]\ar@{-}[dl]&*=0{}\ar@{-}[dl]
\\*=0{}\ar@{-}[r]&*=0{}&*=0{}&*=0{}\ar@{-}[r]\ar@{-}[dl]&*=0{}\ar@{-}[r]\ar@{-}[dl]&*=0{}\ar@{-}[r]\ar@{-}[dl]&*=0{}
\\*=0{}&*=0{}&*=0{}\ar@{-}[r]\ar@{-}[dl]&*=0{}\ar@{-}[r]\ar@{-}[dl]&*=0{}\ar@{..}[d]
\\*=0{}&*=0{}\ar@{-}[r]&*=0{}&*=0{}&*=0{}\ar@{-}[r]\ar@{-}[dl]&*=0{}\ar@{-}[r]\ar@{-}[dl]&*=0{}\ar@{-}[r]\ar@{-}[dl]&*=0{}\ar@{-}[dl]&*=0{}
\\*=0{}&*=0{}&*=0{}&*=0{}\ar@{-}[r]&*=0{}\ar@{-}[r]&*=0{}\ar@{-}[r]&*=0{}
\\*=0{}}
\text{ and }B=\xymatrix @W=1.0pc @H=1.0pc @R=1.0pc @C=1.0pc
{*=0{}&*=0{}&*=0{}&*=0{}\ar@{-}[r]\ar@{-}[dl]&*=0{}\ar@{-}[r]\ar@{-}[dl]&*=0{}\ar@{-}[dl]
\\*=0{}&*=0{}&*=0{}\ar@{-}[r]\ar@{-}[dl]&*=0{}\ar@{-}[r]\ar@{-}[dl]\ar@{..}[dd]&*=0{}
\\*=0{}&*=0{}\ar@{-}[r]\ar@{-}[dl]&*=0{}\ar@{-}[dl]
\\*=0{}\ar@{-}[r]&*=0{}&*=0{}&*=0{}\ar@{-}[r]\ar@{-}[dl]&*=0{}\ar@{-}[r]\ar@{-}[dl]&*=0{}\ar@{-}[dl]&*=0{}
\\*=0{}&*=0{}&*=0{}\ar@{-}[r]\ar@{-}[dl]&*=0{}\ar@{-}[r]\ar@{-}[dl]\ar@{..}[dd]&*=0{}
\\*=0{}&*=0{}\ar@{-}[r]&*=0{}
\\
*=0{}&*=0{}&*=0{}&*=0{}\ar@{-}[r]\ar@{-}[dl]&*=0{}\ar@{-}[r]\ar@{-}[dl]&*=0{}\ar@{-}[dl]
\\*=0{}&*=0{}&*=0{}\ar@{-}[r]&*=0{}\ar@{-}[r]&*=0{}
\\*=0{}}
\]
then we let
\[
A \otimes_1 B=\xymatrix @W=1.0pc @H=1.0pc @R=1.0pc @C=1.0pc
{*=0{}&*=0{}&*=0{}&*=0{}&*=0{}\ar@{-}[r]\ar@{-}[dl]&*=0{}\ar@{-}[r]\ar@{-}[dl]&*=0{}\ar@{-}[r]\ar@{-}[dl]&*=0{}\ar@{-}[r]\ar@{-}[dl]&*=0{}\ar@{-}[dl]
\\*=0{}&*=0{}&*=0{}&*=0{}\ar@{-}[r]\ar@{-}[dl]&*=0{}\ar@{-}[r]\ar@{-}[dl]&*=0{}\ar@{-}[r]\ar@{-}[dl]&*=0{}\ar@{-}[r]&*=0{}
\\*=0{}&*=0{}&*=0{}\ar@{-}[r]\ar@{-}[dl]&*=0{}\ar@{-}[r]\ar@{-}[dl]&*=0{}\ar@{-}[dl]\ar@{..}[d]
\\*=0{}&*=0{}\ar@{-}[r]\ar@{-}[dl]&*=0{}\ar@{-}[r]\ar@{-}[dl]&*=0{}&*=0{}\ar@{-}[r]\ar@{-}[dl]&*=0{}\ar@{-}[r]\ar@{-}[dl]&*=0{}\ar@{-}[r]\ar@{-}[dl]&*=0{}\ar@{-}[r]\ar@{-}[dl]&*=0{}\ar@{-}[dl]
\\*=0{}\ar@{-}[r]&*=0{}&*=0{}&*=0{}\ar@{-}[r]\ar@{-}[dl]&*=0{}\ar@{-}[r]\ar@{-}[dl]&*=0{}\ar@{-}[r]\ar@{-}[dl]&*=0{}\ar@{-}[r]&*=0{}
\\*=0{}&*=0{}&*=0{}\ar@{-}[r]\ar@{-}[dl]&*=0{}\ar@{-}[r]\ar@{-}[dl]&*=0{}\ar@{..}[d]
\\*=0{}&*=0{}\ar@{-}[r]\ar@{-}[dl]&*=0{}\ar@{-}[dl]&*=0{}&*=0{}\ar@{-}[r]\ar@{-}[dl]&*=0{}\ar@{-}[r]\ar@{-}[dl]&*=0{}\ar@{-}[r]\ar@{-}[dl]&*=0{}\ar@{-}[dl]
\\*=0{}\ar@{-}[r]&*=0{}&*=0{}&*=0{}\ar@{-}[r]\ar@{-}[dl]&*=0{}\ar@{-}[r]\ar@{-}[dl]&*=0{}\ar@{-}[r]\ar@{-}[dl]&*=0{}
\\*=0{}&*=0{}&*=0{}\ar@{-}[r]\ar@{-}[dl]&*=0{}\ar@{-}[r]\ar@{-}[dl]&*=0{}\ar@{..}[d]
\\*=0{}&*=0{}\ar@{-}[r]&*=0{}&*=0{}&*=0{}\ar@{-}[r]\ar@{-}[dl]&*=0{}\ar@{-}[r]\ar@{-}[dl]&*=0{}\ar@{-}[r]\ar@{-}[dl]&*=0{}\ar@{-}[dl]&*=0{}
\\*=0{}&*=0{}&*=0{}&*=0{}\ar@{-}[r]\ar@{-}[dl]&*=0{}\ar@{-}[r]\ar@{-}[dl]\ar@{..}[d]&*=0{}\ar@{-}[r]&*=0{}&*=0{}
\\*=0{}&*=0{}&*=0{}\ar@{-}[r]\ar@{-}[dl]&*=0{}\ar@{-}[dl]&*=0{}\ar@{-}[r]\ar@{-}[dl]&*=0{}\ar@{-}[r]\ar@{-}[dl]&*=0{}\ar@{-}[dl]&*=0{}
\\*=0{}&*=0{}\ar@{-}[r]&*=0{}&*=0{}\ar@{-}[r]\ar@{-}[dl]&*=0{}\ar@{-}[r]\ar@{-}[dl]\ar@{..}[d]&*=0{}
\\*=0{}&*=0{}&*=0{}\ar@{-}[r]&*=0{}&*=0{}\ar@{-}[r]\ar@{-}[dl]&*=0{}\ar@{-}[r]\ar@{-}[dl]&*=0{}\ar@{-}[dl]
\\*=0{}&*=0{}&*=0{}&*=0{}\ar@{-}[r]&*=0{}\ar@{-}[r]\ar@{..}[d]&*=0{}&*=0{}
\\
*=0{}&*=0{}&*=0{}&*=0{}&*=0{}\ar@{-}[r]\ar@{-}[dl]&*=0{}\ar@{-}[r]\ar@{-}[dl]&*=0{}\ar@{-}[dl]
\\*=0{}&*=0{}&*=0{}&*=0{}\ar@{-}[r]&*=0{}\ar@{-}[r]&*=0{}
\\*=0{}}
\text{ , }
A \otimes_2 B=\xymatrix @W=1.0pc @H=1.0pc @R=1.0pc @C=1.0pc
{*=0{}&*=0{}&*=0{}&*=0{}&*=0{}\ar@{-}[r]\ar@{-}[dl]&*=0{}\ar@{-}[r]\ar@{-}[dl]&*=0{}\ar@{-}[r]\ar@{-}[dl]&*=0{}\ar@{-}[r]\ar@{-}[dl] &*=0{}\ar@{-}[r]\ar@{-}[dl] &*=0{}\ar@{-}[r]\ar@{-}[dl] &*=0{}\ar@{-}[dl]
\\*=0{}&*=0{}&*=0{}&*=0{}\ar@{-}[r]\ar@{-}[dl]&*=0{}\ar@{-}[r]\ar@{-}[dl]&*=0{}\ar@{-}[r]\ar@{-}[dl] &*=0{}\ar@{-}[r]\ar@{-}[dl] &*=0{}\ar@{-}[r] &*=0{}\ar@{-}[r]&*=0{}
\\*=0{}&*=0{}&*=0{}\ar@{-}[r]\ar@{-}[dl]&*=0{}\ar@{-}[r]\ar@{-}[dl]&*=0{}\ar@{-}[r]\ar@{-}[dl]&*=0{}\ar@{-}[dl]
\\*=0{}&*=0{}\ar@{-}[r]\ar@{-}[dl]&*=0{}\ar@{-}[r]\ar@{-}[dl]&*=0{}\ar@{-}[r]&*=0{}\ar@{..}[d]
\\*=0{}\ar@{-}[r]&*=0{}&*=0{}&*=0{}&*=0{}\ar@{-}[r]\ar@{-}[dl]&*=0{}\ar@{-}[r]\ar@{-}[dl]&*=0{}\ar@{-}[r]\ar@{-}[dl]&*=0{}\ar@{-}[r]\ar@{-}[dl]&*=0{}\ar@{-}[r]\ar@{-}[dl]&*=0{}\ar@{-}[r]\ar@{-}[dl]&*=0{}\ar@{-}[dl]
\\*=0{}&*=0{}&*=0{}&*=0{}\ar@{-}[r]\ar@{-}[dl]&*=0{}\ar@{-}[r]\ar@{-}[dl]&*=0{}\ar@{-}[r]\ar@{-}[dl]&*=0{}\ar@{-}[r]\ar@{-}[dl]&*=0{}\ar@{-}[r]&*=0{}\ar@{-}[r]&*=0{}
\\*=0{}&*=0{}&*=0{}\ar@{-}[r]\ar@{-}[dl]&*=0{}\ar@{-}[r]\ar@{-}[dl]&*=0{}\ar@{-}[r]\ar@{..}[d]&*=0{}
\\*=0{}&*=0{}\ar@{-}[r]\ar@{-}[dl]&*=0{}\ar@{-}[dl]&*=0{}&*=0{}\ar@{-}[r]\ar@{-}[dl]&*=0{}\ar@{-}[r]\ar@{-}[dl]&*=0{}\ar@{-}[r]\ar@{-}[dl]&*=0{}\ar@{-}[r]\ar@{-}[dl]&*=0{}\ar@{-}[r]\ar@{-}[dl]&*=0{}\ar@{-}[dl]
\\*=0{}\ar@{-}[r]&*=0{}&*=0{}&*=0{}\ar@{-}[r]\ar@{-}[dl]&*=0{}\ar@{-}[r]\ar@{-}[dl]&*=0{}\ar@{-}[r]\ar@{-}[dl]&*=0{}\ar@{-}[r]&*=0{}\ar@{-}[r]&*=0{}
\\*=0{}&*=0{}&*=0{}\ar@{-}[r]\ar@{-}[dl]&*=0{}\ar@{-}[r]\ar@{-}[dl]&*=0{}\ar@{..}[d]
\\*=0{}&*=0{}\ar@{-}[r]&*=0{}&*=0{}&*=0{}\ar@{-}[r]\ar@{-}[dl]&*=0{}\ar@{-}[r]\ar@{-}[dl]&*=0{}\ar@{-}[r]\ar@{-}[dl]&*=0{}\ar@{-}[dl]&*=0{}
\\*=0{}&*=0{}&*=0{}&*=0{}\ar@{-}[r]&*=0{}\ar@{-}[r]&*=0{}\ar@{-}[r]&*=0{}
\\*=0{}}
\]
and
\[
A \otimes_3 B=\xymatrix @W=1.0pc @H=1.0pc @R=1.0pc @C=1.0pc
{*=0{}&*=0{}&*=0{}&*=0{}&*=0{}&*=0{}&*=0{}&*=0{}\ar@{-}[r]\ar@{-}[dl]&*=0{}\ar@{-}[r]\ar@{-}[dl]&*=0{}\ar@{-}[r]\ar@{-}[dl]&*=0{}\ar@{-}[r]\ar@{-}[dl]&*=0{}\ar@{-}[dl]
\\*=0{}&*=0{}&*=0{}&*=0{}&*=0{}&*=0{}&*=0{}\ar@{-}[r]\ar@{-}[dl]&*=0{}\ar@{-}[r]\ar@{-}[dl]&*=0{}\ar@{-}[r]\ar@{-}[dl]&*=0{}\ar@{-}[r]&*=0{}&*=0{}&*=0{}
\\*=0{}&*=0{}&*=0{}&*=0{}&*=0{}&*=0{}\ar@{-}[r]\ar@{-}[dl]&*=0{}\ar@{-}[r]\ar@{-}[dl]&*=0{}\ar@{-}[dl]\ar@{..}[d]&*=0{}&*=0{}&*=0{}&*=0{}
\\*=0{}&*=0{}&*=0{}&*=0{}&*=0{}\ar@{-}[r]\ar@{-}[dl]&*=0{}\ar@{-}[r]\ar@{-}[dl]&*=0{}\ar@{-}[dl]&*=0{}\ar@{-}[r]\ar@{-}[dl]&*=0{}\ar@{-}[r]\ar@{-}[dl]&*=0{}\ar@{-}[r]\ar@{-}[dl]&*=0{}\ar@{-}[r]\ar@{-}[dl]&*=0{}\ar@{-}[dl]&*=0{}&*=0{}&*=0{}&*=0{}
\\*=0{}&*=0{}&*=0{}&*=0{}\ar@{-}[r]\ar@{-}[dl]&*=0{}\ar@{-}[r]\ar@{-}[dl]&*=0{}&*=0{}\ar@{-}[r]\ar@{-}[dl]&*=0{}\ar@{-}[r]\ar@{-}[dl]&*=0{}\ar@{-}[r]\ar@{-}[dl]&*=0{}\ar@{-}[r]&*=0{}&*=0{}&*=0{}&*=0{}&*=0{}
\\*=0{}&*=0{}&*=0{}\ar@{-}[r]\ar@{-}[dl]&*=0{}\ar@{-}[dl]&*=0{}&*=0{}\ar@{-}[r]\ar@{-}[dl]&*=0{}\ar@{-}[r]\ar@{-}[dl]&*=0{}\ar@{-}[dl]\ar@{..}[d]&*=0{}&*=0{}&*=0{}&*=0{}
\\*=0{}&*=0{}\ar@{-}[r]\ar@{-}[dl]&*=0{}\ar@{-}[dl]&*=0{}&*=0{}\ar@{-}[r]\ar@{-}[dl]&*=0{}\ar@{-}[r]\ar@{-}[dl]&*=0{}&*=0{}\ar@{-}[r]\ar@{-}[dl]&*=0{}\ar@{-}[r]\ar@{-}[dl]&*=0{}\ar@{-}[r]\ar@{-}[dl]&*=0{}\ar@{-}[dl]&*=0{}&*=0{}&*=0{}&*=0{}
\\*=0{}\ar@{-}[r]&*=0{}&*=0{}&*=0{}\ar@{-}[r]\ar@{-}[dl]&*=0{}\ar@{-}[dl]&*=0{}&*=0{}\ar@{-}[r]\ar@{-}[dl]&*=0{}\ar@{-}[r]\ar@{-}[dl]&*=0{}\ar@{-}[r]\ar@{-}[dl]&*=0{}&*=0{}&*=0{}&*=0{}&*=0{}
\\*=0{}&*=0{}&*=0{}\ar@{-}[r]\ar@{-}[dl]&*=0{}\ar@{-}[dl]&*=0{}&*=0{}\ar@{-}[r]\ar@{-}[dl]&*=0{}\ar@{-}[r]\ar@{-}[dl]&*=0{}\ar@{-}[dl]\ar@{..}[d]&*=0{}&*=0{}&*=0{}&*=0{}
\\*=0{}&*=0{}\ar@{-}[r]&*=0{}&*=0{}&*=0{}\ar@{-}[r]\ar@{-}[dl]&*=0{}\ar@{-}[r]\ar@{-}[dl]&*=0{}&*=0{}\ar@{-}[r]\ar@{-}[dl]&*=0{}\ar@{-}[r]\ar@{-}[dl]&*=0{}\ar@{-}[r]\ar@{-}[dl]&*=0{}\ar@{-}[dl]&*=0{}&*=0{}&*=0{}&*=0{}&*=0{}
\\*=0{}&*=0{}&*=0{}&*=0{}\ar@{-}[r]&*=0{}&*=0{}&*=0{}\ar@{-}[r]&*=0{}\ar@{-}[r]&*=0{}\ar@{-}[r]&*=0{}&*=0{}&*=0{}&*=0{}&*=0{}&*=0{}
\\*=0{}&*=0{}&*=0{}&*=0{}&*=0{}&*=0{}&*=0{}&*=0{}&*=0{}}
\]
Note that in this restricted setting of decreasing matrices the
lexicographic merging of sequences (rows) of two matrices does not
preserve the total decreasing property (decreasing in rows and
columns).

These three products just shown preserve the total sum of the
entries in both matrices, and do interact via interchanges to form
the structure of a 3-fold monoidal category. Renumbered, they are:
$\otimes_1$ ($z$-axis stacking) is the vertical concatenation of
matrices followed by sorting the new longer columns, $\otimes_2$
($y$-axis stacking) is horizontal concatenation of matrices followed
by sorting the new longer rows and $\otimes_3$  ($x$-axis stacking)
is the addition of matrices.  For comparison, here is the same
example of the products as just given above shown by matrices. Only
the non-zero entries of the matrices are shown.
\[
A = \left[\begin{array}{cccc}4&3&1&1 \\ 4&2&1&1 \\ 3&2&1 \\ 1&1&1\end{array} \right] B= \left[\begin{array}{cc}3&1 \\ 2&1 \\ 1&1 \end{array} \right]
\]
\[
A \otimes_1 B = \left[\begin{array}{cccc} 4&3&1&1 \\ 4&2&1&1 \\ 3&2&1 \\ 3&1&1 \\ 2&1 \\ 1&1 \\ 1&1\end{array} \right]
A\otimes_2 B = \left[\begin{array}{cccccc} 4&3&3&1&1&1 \\ 4&2&2&1&1&1 \\ 3&2&1&1&1 \\ 1&1&1 \end{array} \right]
\text{ and } A\otimes_3 B = \left[\begin{array}{cccc} 7&4&1&1 \\ 6&3&1&1 \\ 4&3&1 \\ 1&1&1\end{array} \right]
\]

\begin{theorem}\label{3dY}
The category of 3-d  Young diagrams with lexicographic ordering and
the products just described possesses the structure  of a 3-fold
monoidal category.
\end{theorem}

The proof will require the following two lemmas.
\begin{lemma}\label{minmax}
For two sequences of $n$ elements each, the first given by $a_1
\dots a_n$ and the second by $b_1 \dots b_n$, then considering pairs
of elements $a_{\sigma(i)}$ and $b_{\tau(i)}$ for permutations
$\sigma, \tau \in S_n$, we have the following inequality:
\[
\max(\min(a_{\sigma(1)},b_{\tau(1)}),\dots,\min(a_{\sigma(n)},b_{\tau(n)}))
\le \min(\max(a_1,\dots,a_n), \max(b_1, \dots, b_n)).
\]
\end{lemma}

\begin{proof}
This is true since for $i=1\dots n$ we have $a_i \le
\max(a_1,\dots,a_n)$ and $b_i \le \max(b_1, \dots, b_n).$ Therefore
$\min(a_{\sigma(i)}, b_{\tau(i)})
    \le \min(\max(a_1,\dots,a_n), \max(b_1,  \dots, b_n))$
and the inequality follows.
\end{proof}

\begin{lemma}\label{matrixsort}
For a given finite matrix $M$ with $n$ rows, we claim that first
sorting each row (greater to lesser) and then sorting each resulting
column gives a final result that is lexicographically less than or
equal to the final result of sorting each column of $M$ and then
each row. The lexicographic ordering here is applied to the
sequences of entries read from the matrices by rows.
\end{lemma}

\begin{proof}
This is seen by a chain of inequalities that each correspond to a
single step in a parallel bubble sorting of the rows of $M.$
Consider the final result of sorting each column vertically and then
each row. We gradually evolve this into the reverse procedure by
performing a series of steps, each of which begins by comparing two
adjacent columns in the current stage of the evolution. The step
consists of  switches that insure each horizontal pair in the
columns is in order, i.e. switching the positions of the two
elements in each row only if the one in the left column is smaller
than the one in the right.  We call this a parallel switch, or just
a switch.  The result of taking the switched matrix and vertically
sorting its columns and then horizontally sorting its rows will be
shown to be lexicographically less than or equal to the result of
vertically sorting columns and then horizontally sorting rows before
the parallel switch. The entire series of steps together constitute
sorting each row of $M.$ Since after vertically sorting a matrix
which began with sorted rows the new rows still remain sorted, then
at the end of the evolution we are indeed doing the reverse
procedure; that is sorting horizontally first and then vertically.

For a single step in the parallel bubble sort, we claim that after
the parallel switch and then vertical sorting of the two adjacent
columns the pairs in each resulting row will be either all identical
to those in the result of vertically sorting the unswitched columns,
or there will be a first row $k$ in which the pair in the switched
version of the columns consists of one element equal to one element
of the corresponding pair in the unswitched version and one element
less than the other element in the unswitched version.

Since no other columns are changed at this step, then this will
imply that after vertically sorting the other columns and then all
the rows in both matrices, the two resulting matrices will be
identical  or just identical up to the  $k^{th}$ row, where the
switched matrix will be lexicographically less than the unswitched.

The claim for two columns follows from repeated application of
Lemma~\ref{minmax}. Let the two columns be $a_1 \dots a_n$ and $b_1
\dots b_n$ After the parallel switching, the left column holds the
max of each pair and the right the min.  Vertical sorting moves the
max of each column to the top row, and leaves all the new rows (of
two elements each) still sorted  left to right. Located in the left
position of the new top row is
\[ \max(\max(a_1,b_1), \dots ,\max(a_n,b_n))
    = \max(\max(a_1, \dots ,a_n),\max(b_1, \dots ,b_n))\]
the latter of which is the in the top row of the vertically sorted
unswitched columns. The right position in the top row of the
switched columns is \[r=\max(\min(a_1,b_1),\dots,\min(a_n,b_n)),\]
which is less than or equal to the other element in the top row of
the vertically sorted unswitched columns
\[s=\min(\max(a_1,\dots,a_n),\max(b_1,\dots,b_n)),\]
by the preceding Lemma~\ref{minmax} (with trivial permutations). If
$r<s$ then we are done. If $r=s$ then we note that the remaining
rows $2\dots n$ contain the same collection of elements $a_i$ and
$b_i$ in both the switched and unswitched columns, i.e. we may
assume that in vertically sorting either version we moved $a_j$ and
$b_l$ to the top row.  Note that since the rows in the switched
version are sorted, $\max(a_l,b_l) \ge \min(a_j,b_j)$ and
$\max(a_j,b_j) \ge \min(a_l,b_l).$ Thus the $\max(a_l,b_j)$ will
always be in the left column and $\min(a_l,b_j)$ in the right.

Then the second row of the vertically sorted switched pair of
columns is
\[\max(\max(a_1,b_1),\dots,\widehat{\max(a_j,b_j)},\dots,
    \widehat{\max(a_l,b_l)},\dots,\max(a_n,b_n),\max(a_l,b_j))\]
in the first position and
\[\max(\min(a_1,b_1),\dots,\widehat{\min(a_j,b_j)},\dots,
    \widehat{\min(a_l,b_l)},\dots,\min(a_n,b_n),\min(a_l,b_j))\]
in the second position, where the hats indicate missing elements.
Whereas the second row of the vertically sorted unswitched columns
is made up of
\[\max(\max(a_1,\dots,\hat{a_j},\dots,a_n))\text{ and }
    \max(\max(b_1,\dots,\hat{b_l},\dots,b_n)).\]
Thus the left position in the second row of the switched version is
the same value as one of the elements in the second row of the
unswitched vertically sorted columns.  By Lemma~\ref{minmax} with
the evident permutations, the right position in the second row is
less than or equal to the other element in the second row of the
unswitched vertically sorted columns.  If less than, then we are
done, if equal then the process continues. If the $1^{st}$ through
$(n-1)^{st}$ rows of the switched and unswitched columns contain the
same values after vertical sorting, then so do the $n\th$ rows. This
completes the proof of the lemma.
\end{proof}

Now we can proceed to the proof of Theorem~\ref{3dY}.
\begin{proof} (of Theorem~\ref{3dY})
We already have existence of $\eta^{23}$ by the argument about
pointwise application of two interchanging products in the proof of
Corollary~\ref{duh}. Here the two products are merging and vertical
stacking applied pointwise to the sequence of rows seen as a
sequence of Young diagrams. To show existence of
$\eta^{13}:(A\otimes_3B)\otimes_1(C\otimes_3D)\to
(A\otimes_1C)\otimes_3(B\otimes_1D)$ we need to check that sorting
each of the columns of two pairs of vertically concatenated matrices
before pointwise adding gives a larger lexicographic result with
respect to rows than adding first and then sorting columns. This
follows from Lemma~\ref{sort}, applied to each pair of sequences
which are the $n^{th}$ columns in the two new matrices formed by
vertically concatenating  $A$ and $C$ and respectively $B$ and $D$,
padded with zeroes so that adding the new matrices results in adding
$A$ and $B$ and respectively $C$ and $D$. From  the lemma then we
have that $(A\otimes_1C)\otimes_3(B\otimes_1D)$ gives a result whose
$n^{th}$ column is lexicographically greater than or equal to the
$n^{th}$ column of $(A\otimes_3B)\otimes_1(C\otimes_3D).$ This
implies that either the pairs of respective columns are each equal
sequences or that there is some least row $i$ and column $j$ such
that all the pairs of columns are identical  in rows less than $i$
and that the two rows $i$ are identical in columns less than $j$,
but that the $i,j$ position in $(A\otimes_3B)\otimes_1(C\otimes_3D)$
is less than the corresponding position in
$(A\otimes_1C)\otimes_3(B\otimes_1D).$  Thus the existence of the
required inequality is shown.

The existence of $\eta^{12}$ is due to the fact
that we are ordering the matrices by giving precedence to the rows.
The two four-fold products can be seen as two alternate operations
on a single large matrix $M$. This matrix is constructed by arranging
$A,B,C,D$ with added zeroes so that $(A\otimes_1C)\otimes_2(B\otimes_1D)$
is the result of first sorting each column vertically, greater
values at the top, and then each row horizontally, greater values
to the left, while $(A\otimes_2B)\otimes_1(C\otimes_2D)$ is achieved
by sorting horizontally first and then vertically.  Recall that in
the ordering rows are given precedence over columns.  Here is an
illustration of the inequality, showing the process of constructing
the large matrix.
\[
A = \left[\begin{array}{ccc} 3&3&2 \\ 1&1 \end{array} \right] B = \left[\begin{array}{c} 9 \\9 \\9 \end{array} \right]
C = \left[\begin{array}{c} 2 \\ 1 \end{array} \right] D= \left[\begin{array}{c} 5 \end{array} \right]
\]
\[
M = \left[\begin{array}{cccc} 3&3&2&9 \\ 1&1&0&9 \\ 0&0&0&9 \\ 2&0&0&5 \\ 1&0&0&0 \end{array} \right]
\]
\[
(A\otimes_2B)\otimes_1(C\otimes_2D) = \left[\begin{array}{cccc} 9&3&3&2 \\ 9&2&1 \\ 9&1 \\ 5 \\ 1 \end{array} \right]
<  \left[\begin{array}{cccc} 9&3&3&2 \\ 9&2&1 \\ 9&1 \\ 5&1 \end{array} \right]= (A\otimes_1C)\otimes_2(B\otimes_1D)
\]
The proof that this inequality always holds uses
Lemma~\ref{matrixsort}.
 By applying that lemma to the large
matrix $M$ constructed of the four matrices $A,B,C,D$ as described
above, we have the proof of the theorem.
\end{proof}

Now we define the general $n$-fold monoidal category of $n$-dimensional
Young diagrams.  The proof of the theorem for three dimensions plays
an important role in the general theorem,  since each interchanger
involves two products.  Once we have decided to represent Young
diagrams of higher dimension by arrays of natural numbers which
decrease in each index, it is clear that each interchanger will
either involve directly two of the indices of the array or one index
as well as pointwise addition.

\begin{definition} The category of $n$-dimensional Young diagrams
consists of
\begin{enumerate}
\item Objects $A_{i_1 i_2\dots i_{n-1}}$, finitely nonzero
$n$-dimensional arrays of nonnegative integers which are monotone
decreasing in each index, and
\item Morphisms the order relations in the lexicographic ordering
with precedence given to lesser indices.
\end{enumerate}
\end{definition}

There are $n$ ways to take a product of two $n$-dimensional Young
diagrams, which we visualize as arrays of natural numbers in $n-1$
dimensions. The products correspond to merging, i.e. concatenating
and then sorting, in each of the $n-1$ possible directions, as well
as pointwise addition as $\otimes_n$. The order of products is the
reverse of the order of the indices. That is, for $k=1\dots n-1,$
$\otimes_k$ is merging in the direction of the index $i_{n-k}.$

\begin{theorem}
The category of $n$-dimensional Young diagrams with the above
products constitutes an $n$-fold monoidal category.
\end{theorem}

\begin{proof}
We must show the existence of the interchangers $\eta^{jk}$ as
inequalities for $1\le j<k\le n.$ First we demonstrate the existence
of the required inequality when $k < n.$ For $A,B,C,D$ $n$-dimensional
Young diagrams seen as $(n-1)$-dimensional arrays, we let
$M_{i_1 i_2\dots i_{n-1}}$ be a large array made by concatenating $A$
and $B$ in the direction of the index $i_k$, concatenating $C$ and $D$
in the direction of the index $i_k$, and then concatenating those
two results in the direction of the index $i_j$. Zeros are added
(see above for the two dimensional array example) so that the
products
$(A\otimes_k B)\otimes_j(C\otimes_k D)$ and
$(A\otimes_j C)\otimes_k(B\otimes_j D)$
can then both be described as sorting $M_{i_1 i_2\dots i_{n-1}}$
in two directions; first $i_k$ then $i_j$ or vice versa respectively.
That the inequality holds is seen as we compare the results position
by position in the lexicographic order, i.e. reading lower indices
first.  The first differing value we come upon, say in location
$(i_1 i_2\dots i_{n-1}),$ then will necessarily be the first
difference in  the  sub-array of two dimensions in the directions
$i_j$ and $i_k$ determined by the location $(i_1 i_2\dots i_{n-1}).$
Thus by the proof of Lemma~\ref{matrixsort}, the value in
$(A\otimes_k B)\otimes_j(C\otimes_k D)$
is less than the corresponding value in
$(A\otimes_j C)\otimes_k(B\otimes_j D).$

Secondly we check the cases that have $k=n.$ We can see the
four-fold products as operations on two arrays, one made by
concatenating $A$ and $C$ in the $i_j$ direction, and another made
by concatenating $B$ and $D$ in the $i_j$ direction, padded with
zeroes so that adding the two pointwise results in pointwise
addition of $A$ with $B$, and of $C$ with $D.$ Then $(A\otimes_k
B)\otimes_j(C\otimes_k D)$ is adding first and then sorting in the
$i_j$ direction, while $(A\otimes_j C)\otimes_k(B\otimes_j D)$ is
the reverse process. To see that the correct inequality holds we
again compare the results position by position in lexicographic
order. The first differing value is also the first difference
between the two corresponding 2 dimensional sub-arrays which are in
the directions $i_j$ and $i_{n-1}.$ These sub-arrays are the results
of sorting and then pointwise addition and vice versa respectively,
and so by the proof for existence of $\eta^{13}$ in
Theorem~\ref{3dY} the desired result is shown.
\end{proof}

\section{Examples of \texorpdfstring{$n$}{N}-fold operads}

The categories from Section~\ref{section:categoryexamples} give us a
domain in which we can exhibit some concrete examples of operads. To
have an operad with an element ${\cal C}(0)$ we will need to
``compactify'' by adjoining a new initial  object with the desired
properties to the example categories based on totally ordered sets.

\begin{definition} For ${\cal V}$ an $n$-fold monoidal category whose morphisms are the
$\le$ relations of a totally ordered set with least element, we
define its \emph{compactification} $\overline{\cal V}$ by adding a
new
  initial object which we will
denote by $\emptyset.$ The morphisms will still be given by the
original ordering augmented by letting $\emptyset$ be the new least
element. All the original products have their original definition on
objects of ${\cal V}$. However, for $i=1\dots n$ and $A$ an object
of $\overline{\cal V}$ we let $\emptyset\otimes_i A = \emptyset$ and
$A\otimes_i\emptyset = \emptyset.$ Note that even when the original
product in ${\cal V}$ was given by max, the new product in
$\overline{\cal V}$ gives $\emptyset$ when one of the operands is
$\emptyset.$
 \end{definition}

\begin{theorem}
For ${\cal V}$ an $n$-fold monoidal category whose morphisms are the
$\le$ relations of a totally ordered set with least element, the
compactification $\overline{\cal V}$ has the structure of an
$n$-fold monoidal category. Moreover  $(\overline{\cal
V},\coprod,\otimes_1,\dots, \otimes_n)$ is a strict $(n+1)$-fold
monoidal category with distinct units for which forgetting the first
tensor product (given by the coproduct) recovers $\overline{\cal
V}.$
\end{theorem}

\begin{proof}
Our products are all still strictly associative. By definition
$\emptyset\otimes_i I = \emptyset.$ The diagrams will all commute
since the morphisms are the ordering. Therefore we only need to
check that the interchangers still exist when one of the operands is
$\emptyset.$ Indeed in this case the two products in question both
become $\emptyset$ and the interchanger is the identity. Unit
conditions are still obeyed for the  same reason.

 Now
$\overline{\cal V}$ has coproducts given by $\coprod = \max$, where
max is taken with respect to the new total order with $\emptyset$ as
least element. Thus $\emptyset$ is the strict unit for $\coprod.$
Denoting by 0 the common unit of the other products, by definition
we have that $0\coprod 0 = 0$ and $\emptyset \otimes_i \emptyset =
\emptyset.$ Therefore the unit conditions for the interchangers
involving the two units hold as well.
  We have already demonstrated that
$\coprod$ (max) interchanges with any product which preserves
ordering. Our new products of $\overline{\cal V}$ do preserve the
new ordering.
\end{proof}

 In all the following examples the operad
composition is associative since it is based upon ordering, so all
we need check for is the existence of that composition. We will
refer to the example categories developed in the previous section,
but assume that we are dealing with their compactification. Note
that each of the following examples satisfy the hypothesis of
Theorem~\ref{foo} since $\coprod$ (max) distributes over each
$\otimes_i$, since each product preserves the ordering.

\begin{example}
Of course ${\cal C}(j) = \emptyset$ and ${\cal C}(j) = 0$ for all
$j$ are trivially operads, where $0$ is the monoidal unit.  First
we look at the simplest interesting examples: 2-fold operads in an
ordered monoid such as $\nat$, where $\otimes_1$ is max and $\otimes_2$
is $+$. We always set ${\cal C}(0) = \emptyset$ but often only list the later terms.
 A nontrivial
2-fold operad in $\nat$ is a nonzero sequence $\{{\cal C}(j)\}_{j\ge
0}$ of natural numbers which has the property that for any $j_1
\dots j_k$, $\max({\cal C}(k), \sum{{\cal C}(j_i)}) \le {\cal
C}(\sum{j_i})$ and for which ${\cal C}(1) = 0.$  This translates
into saying that for any two whole numbers $x, y$ we have that
${\cal C}(x+y) \ge {\cal C}(x)+{\cal C}(y)$ and that ${\cal C}(1) =
0.$ The latter condition both satisfies the unit axioms and makes it
redundant to also insist that the sequence be monotone increasing.
Perhaps the first example that comes to mind is the Fibonacci
numbers. Minimal examples are formed by choosing a starting term or
terms and then determining each later $n^{th}$ term. These are
minimal in the sense that the principle which determines each of the
later terms in succession is that of choosing the minimal next term
out of all possible such terms. For a starting finite sequence
$0,a_2 , \dots ,a_l$ which obeys the the axioms of a 2-fold operad
so far, the operad ${\cal C}_{0,a_2, \dots ,a_l}$ is found by
defining terms ${\cal C}_{a_1 , \dots ,a_l}(n)$ for $n>l$ to be the
maximum (in general the coproduct!) of all the values of $\max({\cal
C}(k), \sum_{i=1}^k{{\cal C}(j_i)})$ where the sum of the $j_i$ is
$n.$ Some basic examples are the following sequences.

\begin{alignat*}{2}
{\cal C}_{0,1} &= (\emptyset,0,1,1,2,2,3,3,\dots),
    &\qquad {\cal C}_{0,0,1}&=(\emptyset,0,0,1,1,1,2,2,2,3,3,3,\dots) \\
{\cal C}_{0,2} &= (\emptyset,0,2,2,4,4,6,6,\dots),
    &\qquad {\cal C}_{0,0,2} &= (\emptyset,0,0,2,2,2,4,4,4,6,6,6,\dots)
\end{alignat*}
and
\[
{\cal C}_{0,1,2,4,8} =
(\emptyset,0,1,2,4,8,8,9,10,12,16,16,17,18,20,24,\dots).
\]
It is clear that the growth of these sequences oscillates
around linear growth in a predictable way.
\end{example}

\begin{theorem}
If ``arbitrary'' starting terms $0,a_2, \dots ,a_k \in \nat$ are given
(themselves of course obeying the axioms of a 2-fold operad), then
the $n^{th}$ term of the 2-fold operad ${\cal C}_{0,a_2, \dots ,a_k}$ in $\nat$ is given by
\[ a_n  = a_q + pa_k \text{ where } n = pk + q,
    \text{ for } p \in \nat, 0\le q <k.\]
\end{theorem}

\begin{proof}
We need to show that
$$a_n = \max\limits_{j_1+\dots+j_l = n}\{\max(a_l,
\sum_{i=1}^l{a_{j_i}})\} = a_q + pa_k$$ where $n = pk + q$, for $p
\in \nat$, $0\le q <k.$ First we note that $a_q + pa_k$ appears as a
term in the overall max, so that $a_n \ge a_q + pa_k.$

Now we check that $\max(a_l, \sum_{i=1}^l{a_{j_i}})\} $ is always
less than or equal to $ a_q + pa_k$. We need only consider the cases
in which $l < n.$ Since $a_l$ is included  at least once as one of
the $a_{j_i},$ we need to show only that $\sum_{i=1}^l{a_{j_i}}$ is
always less than or equal to $a_q + pa_k$ where the sum of the $j_i$
is $n.$ This follows by strong induction on $n$. The base cases $n=
1\dots k$ hold by definition. We can assume $j_i > 0$ since ${\cal
C}(0) = \emptyset.$ Let $j_i = p_ik+q_i$ for $p_i \in \nat$, and
$0\le q_i <k.$ We may assume without loss of generality that at
least one of the $j_i \ge k,$ since if not then the sum of the
$a_{j_i}$ is less than another sum with $a_{j_1+j_2}$ replacing
$a_{j_1}+ a_{j_2},$ and $k<n.$  Then $\sum{q_i} = n - k\sum{p_i} =
pk+q -k\sum{p_i} < n .$ Thus we have:
\begin{align*}
\sum_{i=1}^l{a_{j_i}} &=\sum{a_{q_i}} + a_k\sum{p_i}    \\
&\le a_{(k(p-\sum{p_i})+q)} + a_k\sum{p_i}      \\
&=a_q + (p-\sum{p_i})a_k + a_k\sum{p_i}         \\
&=a_q+ pa_k.\qedhere
\end{align*}
The first inequality is by the assumption that the terms in the
sequence do form an operad, and the following equality is by our
induction assumption.
\end{proof}

\begin{example}\label{bee}
Consider the 3-fold monoidal category  $\seq(\nat,+)$ of lexicographically
ordered finitely nonzero sequences of the natural numbers (here we
use $\nat$ considered as an example of an ordered monoid), with
products $\otimes_1$ the lexicographic max , $\otimes_2$ the
concatenation and $\otimes_3$ the pointwise addition.  An example
of a 2-fold operad in  $\seq(\nat,+)$ that is not a 3-fold operad
is the following:

Let ${\cal B}(0) = \emptyset$ and
let ${\cal B}(j)_i = 1$ for $i<j~,~0$ otherwise. We can picture
these as follows:
\[
{\cal B}(1) = \xymatrix@W=1.2pc @H=1.2pc @R=1.2pc @C=1.2pc {*=0{}\ar@{-}[r] &*=0{}\ar@{-}[r] &*=0{}}~,~~
{\cal B}(2) = \xymatrix@W=.5pc @H=.5pc @R=0pc @C=0pc @*[F-]{~}~,~~
{\cal B}(3) = \xymatrix@W=.5pc @H=.5pc @R=0pc @C=0pc @*[F-]{~&~}~,~~
{\cal B}(4) = \xymatrix@W=.5pc @H=.5pc @R=0pc @C=0pc @*[F-]{~&~&~}~,~~
{\cal B}(5) = \xymatrix@W=.5pc @H=.5pc @R=0pc @C=0pc @*[F-]{~&~&~&~}~,~~ \dots
\]
This is a 2-fold operad, with respect to the lexicographic max and
concatenation.  For instance the instance of composition
$\gamma^{12}:{\cal B}(3)\otimes_1({\cal B}(2)\otimes_2{\cal B}(1)
    \otimes_2{\cal B}(3)) \to {\cal B}(6)$
appears as the relation:
\[
\xymatrix@W=1.2pc @H=1.2pc @R=1.2pc @C=1.2pc {*=0{}\ar@{-}[r]\ar@{-}[d] &*=0{}\ar@{-}[r]\ar@{-}[d] &*=0{}\ar@{-}[r]\ar@{-}[d] &*=0{}\ar@{-}[d]\\
*=0{}\ar@{-}[r] &*=0{}\ar@{-}[r] &*=0{}\ar@{-}[r] &*=0{}\\
}
\raisebox{-0.5em}{ $<$ } \xymatrix@W=1.2pc @H=1.2pc @R=1.2pc @C=1.2pc {*=0{}\ar@{-}[r]\ar@{-}[d] &*=0{}\ar@{-}[r]\ar@{-}[d] &*=0{}\ar@{-}[r]\ar@{-}[d] &*=0{}\ar@{-}[r]\ar@{-}[d] &*=0{}\ar@{-}[r]\ar@{-}[d] &*=0{}\ar@{-}[d]\\
*=0{}\ar@{-}[r] &*=0{}\ar@{-}[r] &*=0{}\ar@{-}[r] &*=0{}\ar@{-}[r] &*=0{}\ar@{-}[r] &*=0{}}
\]
However, the relation
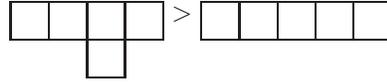
\[
\xymatrix@W=1.2pc @H=1.2pc @R=1.2pc @C=1.2pc {*=0{}\ar@{-}[r]\ar@{-}[d] &*=0{}\ar@{-}[r]\ar@{-}[d] &*=0{}\ar@{-}[r]\ar@{-}[d] &*=0{}\ar@{-}[r]\ar@{-}[d] &*=0{}\ar@{-}[d]\\
*=0{}\ar@{-}[r] &*=0{}\ar@{-}[r] &*=0{}\ar@{-}[r]\ar@{-}[d] &*=0{}\ar@{-}[r]\ar@{-}[d] &*=0{}\\
*=0{} &*=0{} &*=0{}\ar@{-}[r] &*=0{}}
\raisebox{-0.5em}{ $>$ } \xymatrix@W=1.2pc @H=1.2pc @R=1.2pc @C=1.2pc {*=0{}\ar@{-}[r]\ar@{-}[d] &*=0{}\ar@{-}[r]\ar@{-}[d] &*=0{}\ar@{-}[r]\ar@{-}[d] &*=0{}\ar@{-}[r]\ar@{-}[d] &*=0{}\ar@{-}[r]\ar@{-}[d] &*=0{}\ar@{-}[d]\\
*=0{}\ar@{-}[r] &*=0{}\ar@{-}[r] &*=0{}\ar@{-}[r] &*=0{}\ar@{-}[r] &*=0{}\ar@{-}[r] &*=0{}}
\]
shows that
$\gamma^{23}:{\cal B}(3)\otimes_2 ({\cal B}(1)\otimes_3{\cal B}(3)
    \otimes_3{\cal B}(2)) \to {\cal B}(6)$
does not exist, so that ${\cal B}$ is not a 3-fold operad.
\end{example}

\begin{example}
Next we give an example of a 3-fold operad in $\seq(\nat,+)$.
Let ${\cal C}(0) = \emptyset$ and
let ${\cal C}(j) = (j-1,0 \dots).$ We can picture these as follows:
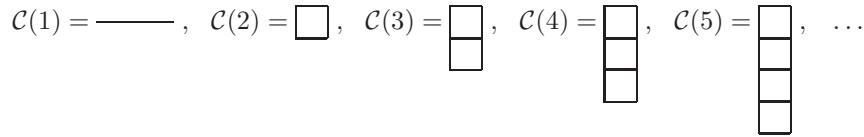
\[
{\cal C}(1) = \xymatrix@W=1.2pc @H=1.2pc @R=1.2pc @C=1.2pc {*=0{}\ar@{-}[r] &*=0{}\ar@{-}[r] &*=0{}}~,~~
{\cal C}(2) = \xymatrix@W=.5pc @H=.5pc @R=0pc @C=0pc @*[F-]{~}~,~~
{\cal C}(3) = \xymatrix@W=.5pc @H=.5pc @R=0pc @C=0pc @*[F-]{~\\~}~,~~
{\cal C}(4) = \xymatrix@W=.5pc @H=.5pc @R=0pc @C=0pc @*[F-]{~\\~\\~}~,~~
{\cal C}(5) = \xymatrix@W=.5pc @H=.5pc @R=0pc @C=0pc @*[F-]{~\\~\\~\\~}~,~~\dots
\]

First we note that the operad ${\cal C}$ just given is a 3-fold
operad since we have that the
$\gamma^{23}:{\cal C}(k)\otimes_2({\cal C}(j_i)\otimes_3\dots
    \otimes_3{\cal C}(j_k)) \to {\cal C}(j)$
exists.  For instance
$\gamma^{23}:{\cal C}(3)\otimes_2 ({\cal C}(1)\otimes_3
    {\cal C}(3)\otimes_3{\cal C}(2)) \to {\cal C}(6)$
appears as the relation
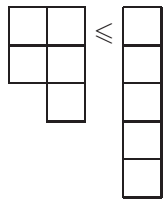
\[
\xymatrix@W=1.2pc @H=1.2pc @R=1.2pc @C=1.2pc {*=0{}\ar@{-}[r]\ar@{-}[d] &*=0{}\ar@{-}[r]\ar@{-}[d] &*=0{}\ar@{-}[d]\\
*=0{}\ar@{-}[r]\ar@{-}[d] &*=0{}\ar@{-}[r]\ar@{-}[d] &*=0{}\ar@{-}[d]\\
*=0{}\ar@{-}[r] &*=0{}\ar@{-}[r]\ar@{-}[d] &*=0{}\ar@{-}[d]\\
*=0{} &*=0{}\ar@{-}[r] &*=0{}}
\raisebox{-1em}{ $\le$ }
\xymatrix@W=1.2pc @H=1.2pc @R=1.2pc @C=1.2pc {*=0{}\ar@{-}[r]\ar@{-}[d] &*=0{}\ar@{-}[d]\\
*=0{}\ar@{-}[r]\ar@{-}[d] &*=0{}\ar@{-}[d]\\
*=0{}\ar@{-}[r]\ar@{-}[d] &*=0{}\ar@{-}[d]\\
*=0{}\ar@{-}[r]\ar@{-}[d] &*=0{}\ar@{-}[d]\\
*=0{}\ar@{-}[r]\ar@{-}[d] &*=0{}\ar@{-}[d]\\
*=0{}\ar@{-}[r] &*=0{}}
\]

Then we remark that as expected the composition
$\gamma^{12}:{\cal C}(k)\otimes_1({\cal C}(j_i)\otimes_2\dots
    \otimes_2{\cal C}(j_k)) \to {\cal C}(j)$
also exists.  For instance
$\gamma^{12}:{\cal C}(3)\otimes_1 ({\cal C}(1)\otimes_2{\cal C}(2)
    \otimes_2{\cal C}(3)) \to {\cal C}(6)$
appears as the relation
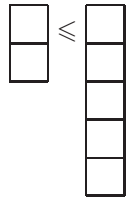
\[
\xymatrix@W=1.2pc @H=1.2pc @R=1.2pc @C=1.2pc {*=0{}\ar@{-}[r]\ar@{-}[d] &*=0{}\ar@{-}[d]\\
*=0{}\ar@{-}[r]\ar@{-}[d] &*=0{}\ar@{-}[d]\\
*=0{}\ar@{-}[r] &*=0{}}
\raisebox{-1em}{ $\le$ }
\xymatrix@W=1.2pc @H=1.2pc @R=1.2pc @C=1.2pc {*=0{}\ar@{-}[r]\ar@{-}[d] &*=0{}\ar@{-}[d]\\
*=0{}\ar@{-}[r]\ar@{-}[d] &*=0{}\ar@{-}[d]\\
*=0{}\ar@{-}[r]\ar@{-}[d] &*=0{}\ar@{-}[d]\\
*=0{}\ar@{-}[r]\ar@{-}[d] &*=0{}\ar@{-}[d]\\
*=0{}\ar@{-}[r]\ar@{-}[d] &*=0{}\ar@{-}[d]\\
*=0{}\ar@{-}[r] &*=0{}}
\]
\end{example}

\begin{example}
Now we consider some products of the previous two described operads
in $\seq(\nat,+)$. We expect ${\cal B} {\otimes'}{\cal C}$ given by
$({\cal B} {\otimes'}{\cal C})(j) = {\cal B}(j) \otimes_3{\cal C}(j)$
to be a 2-fold operad and it is. It appears thus:
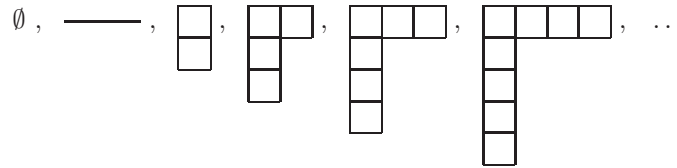
\[
\emptyset~,~~ \xymatrix@W=1.2pc @H=1.2pc @R=1.2pc @C=1.2pc {*=0{}\ar@{-}[r] &*=0{}\ar@{-}[r] &*=0{}}~,~~
\xymatrix@W=.5pc @H=.5pc @R=0pc @C=0pc @*[F-]{~\\~}~,~~
\xymatrix@W=.5pc @H=.5pc @R=0pc @C=0pc @*[F-]{~&~\\~\\~}~,~~
\xymatrix@W=.5pc @H=.5pc @R=0pc @C=0pc @*[F-]{~&~&~\\~\\~\\~}~,~~
\xymatrix@W=.5pc @H=.5pc @R=0pc @C=0pc @*[F-]{~&~&~&~\\~\\~\\~\\~}~,~~
\dots
\]

We demonstrate the tightness of the existence of products of operads
by pointing out that $D(j) = {\cal B}(j) \otimes_2{\cal C}(j)$ does
not form an operad. We leave it to the reader to demonstrate this
fact.
\end{example}

Now we pass to the categories of Young diagrams in which the
interesting products are given by horizontal and vertical stacking.
It is important that we do not restrict the morphisms to those
between diagrams of the same total number of blocks in order that
all  the operad compositions exist.

\begin{theorem}\label{suff}
A sequence of Young diagrams ${\cal C}(n),~ n\in\nat$, in the
category $\modseq(\nat,+)$, is a 2-fold operad if $C(0) = \emptyset$
and for $n \ge 1$, $h({\cal C}(n)) = f(n)$ where $f\colon{\mathbb
Z_+}\to\nat$ is a function such that $f(1) = 0$ and $f(i+j)\ge f(i)
+ f(j)$.
\end{theorem}

\begin{proof}
These conditions are not necessary, but they are sufficient since
the first implies that ${\cal C}(1) =0$ which shows that the unit
conditions are satisfied; and the second implies that the maps
$\gamma$ exist. We see existence of $\gamma^{12}$ since for $j_i > 0$,
$h({\cal C}(k)\otimes_1({\cal C}(j_1)\otimes_2\dots
    \otimes_2{\cal C}(j_k))) =\max(f(k),\max(f(j_i))) \le f(j).$
We have existence of $\gamma^{13}$  and $\gamma^{23}$ since
$\max(f(k), \sum{f(j_i)}) \le f(j)$.
\end{proof}

\begin{example}
Examples of $f$ include  $(x-1)P(x)$ where $P$ is a nonzero
polynomial with coefficients in $\nat$. This is easy to show since
then $P$ will be monotone increasing for $x\ge 1$ and thus
$(i+j-1)P(i+j) = (i-1)P(i+j)+jP(i+j) > (i-1)P(i)+jP(j) - P(j).$ By
this argument we can also use any $f= (x-1)g(x)$ where
$g\colon\nat\to\nat$ is monotone increasing for $x \ge 1$.
\end{example}

For a specific example with a handy picture that also illustrates
again the nontrivial use of the interchange $\eta$ we simply let
$f = x-1.$ Then we have to actually describe the elements of
$\modseq(\nat)$ that make up the operad. One nice choice is the
operad ${\cal C}$ where ${\cal C}(n) = \{n-1,n-1,...,n-1\},$ the
$(n-1)\times (n-1)$ square Young diagram.
\[ {\cal C}(1) = 0,
\ {\cal C}(2) = \xymatrix@W=.6pc @H=.6pc @R=0pc @C=0pc @*[F-]{~},
\ {\cal C}(3) = \xymatrix@W=.6pc @H=.6pc @R=0pc @C=0pc @*[F-]{~&~\\~&~}\ ,
    \ldots
\]
For instance $\gamma^{23}:{\cal C}(3)\otimes_2 ({\cal C}(1)\otimes_3{\cal C}(3)\otimes_3{\cal C}(2)) \to {\cal C}(6)$
appears as the relation
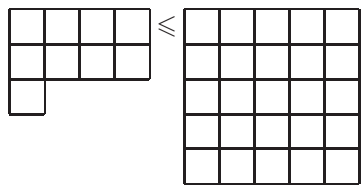
\[
\xymatrix@W=.6pc @H=.6pc @R=0pc @C=0pc @*[F-]{~&~&~&~\\~&~&~&~\\~} \le \xymatrix@W=.6pc @H=.6pc @R=0pc @C=0pc @*[F-]{~&~&~&~&~\\~&~&~&~&~\\~&~&~&~&~\\~&~&~&~&~\\~&~&~&~&~}
\]
An instance of the associativity diagram with upper left position
${\cal C}(2)\otimes_2({\cal C}(3)\otimes_3{\cal C}(2))\otimes_2
({\cal C}(2)\otimes_3{\cal C}(2)\otimes_3{\cal C}(4)\otimes_3{\cal C}(5)\otimes_3{\cal C}(3))$ is as follows:

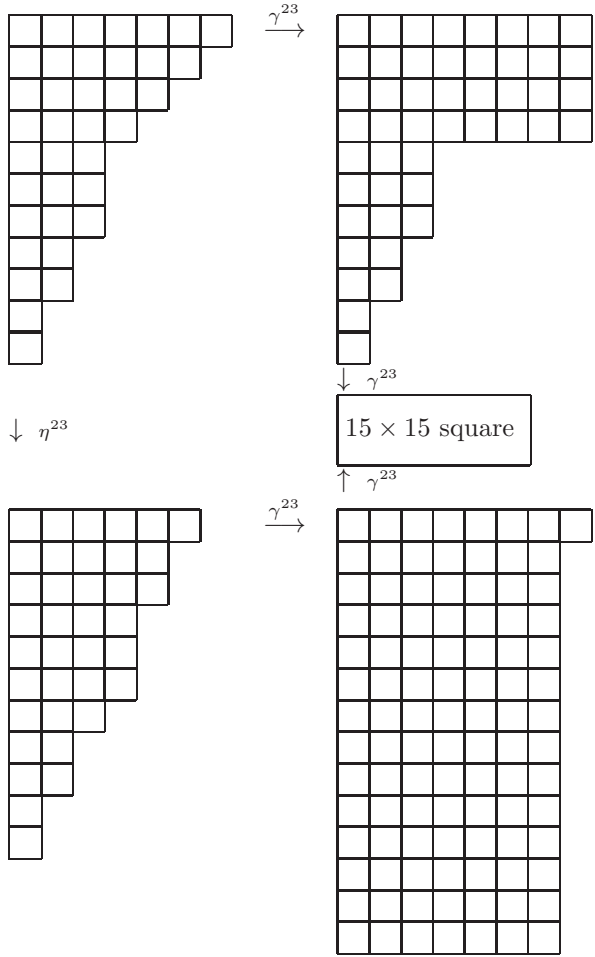
\begin{tabular}{llll}
&$\xymatrix@W=.5pc @H=.5pc @R=0pc @C=0pc @*[F-]{~&~&~&~&~&~&~\\~&~&~&~&~&~\\~&~&~&~&~\\~&~&~&~\\~&~&~\\~&~&~\\~&~&~\\~&~\\~&~\\~\\~}$&$\buildrel\gamma^{23}\over\longrightarrow$&$\xymatrix@W=.5pc @H=.5pc @R=0pc @C=0pc @*[F-]{~&~&~&~&~&~&~&~\\~&~&~&~&~&~&~&~\\~&~&~&~&~&~&~&~\\~&~&~&~&~&~&~&~\\~&~&~\\~&~&~\\~&~&~\\~&~\\~&~\\~\\~}$ \\
&&&$\downarrow~{\scriptstyle\gamma^{23}}$\\
&$\downarrow~{\scriptstyle\eta^{23}}$&&$\xymatrix@W=1.7pc @H=1.7pc @R=0pc @C=0pc @*[F-]{15\times 15 \text{ square }}$\\
&&&$\uparrow~{\scriptstyle\gamma^{23}}$\\
&$\xymatrix@W=.5pc @H=.5pc @R=0pc @C=0pc @*[F-]{~&~&~&~&~&~\\~&~&~&~&~\\~&~&~&~&~\\~&~&~&~\\~&~&~&~\\~&~&~&~\\~&~&~\\~&~\\~&~\\~\\~}$&$\buildrel\gamma^{23}\over\longrightarrow$&$\xymatrix@W=.5pc @H=.5pc @R=0pc @C=0pc @*[F-]{~&~&~&~&~&~&~&~\\~&~&~&~&~&~&~\\~&~&~&~&~&~&~\\~&~&~&~&~&~&~\\~&~&~&~&~&~&~\\~&~&~&~&~&~&~\\~&~&~&~&~&~&~\\~&~&~&~&~&~&~\\~&~&~&~&~&~&~\\~&~&~&~&~&~&~\\~&~&~&~&~&~&~\\~&~&~&~&~&~&~\\~&~&~&~&~&~&~\\~&~&~&~&~&~&~}$\\
&\\
\end{tabular}

\begin{example}
Again we note that the conditions in Theorem~\ref{suff} are not
necessary ones. In fact, given any Young diagram $B$ we can
construct a unique operad that is minimal in each term with respect
to ordering of the diagrams. Again by minimal we mean that the
principle which determines each of the later terms in succession is
that of choosing the minimal next term out of all possible such
terms.
\end{example}
\begin{definition}
The 2-fold operad in the category of Young diagrams generated by a
Young diagram $B$ is  denoted by ${\cal C}_B$ and defined as
follows: ${\cal C}_B(1) =0$ and ${\cal C}_B(2) = B.$ Each successive
term is defined to be the lexicographic maximum of all the products
of prior terms which compose to the term in question; for $n>2$ and
over $\sum{j_i} = n$:
\[ {\cal C}_B(n) = \max\{{\cal C}_B(k)\otimes_2({\cal C}_B(j_1)
    \otimes_3\dots\otimes_3{\cal C}_B(j_k))\}. \]
\end{definition}

\begin{theorem}
If a Young diagram $B$ consists of total number of blocks $q$, then the term ${\cal C}_B(n)$ of the operad generated by $B$ consists of $q(n-1)$ blocks.
\end{theorem}

\begin{proof}
The proof is by strong induction. The number of blocks is given for
${\cal C}_B(1)$ and ${\cal C}_B(2)$. Since the definition is in
terms of a maximum over composable products, if the number of blocks
in each piece of any such a product is assumed by induction to be
respectively $q(k-1)$, and  $q(j_1 -1) \dots q(j_k -1)$, then the
total number of blocks in each product (and thus the maximum) is
$q(n-1)$ since $\sum{j_i} = n$.
\end{proof}

Here are the first few terms of the operad thus generated by
$B = 
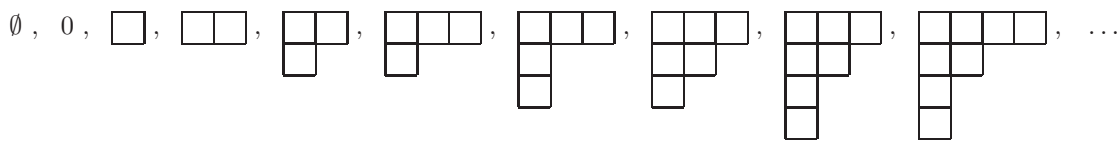
\xymatrix@W=.5pc @H=.5pc @R=0pc @C=0pc @*[F-]{~}~~$.
\[
\emptyset~,~~0~,~~\xymatrix@W=.5pc @H=.5pc @R=0pc @C=0pc @*[F-]{~}~,~~
\xymatrix@W=.5pc @H=.5pc @R=0pc @C=0pc @*[F-]{~&~}~,~~
\xymatrix@W=.5pc @H=.5pc @R=0pc @C=0pc @*[F-]{~&~\\~}~,~~
\xymatrix@W=.5pc @H=.5pc @R=0pc @C=0pc @*[F-]{~&~&~\\~}~,~~
\xymatrix@W=.5pc @H=.5pc @R=0pc @C=0pc @*[F-]{~&~&~\\~\\~}~,~~
\xymatrix@W=.5pc @H=.5pc @R=0pc @C=0pc @*[F-]{~&~&~\\~&~\\~}~,~~
\xymatrix@W=.5pc @H=.5pc @R=0pc @C=0pc @*[F-]{~&~&~\\~&~\\~\\~}~,~~
\xymatrix@W=.5pc @H=.5pc @R=0pc @C=0pc @*[F-]{~&~&~&~\\~&~\\~\\~}~,~~\dots
\]
Note that  height of any given column grows linearly, but that the
length of any row grows logarithmically.
\begin{theorem}
The minimal operad ${\cal C}_{\xymatrix@W=.35pc @H=.35pc @R=.35pc @C=.35pc {*=0{}\ar@{-}[r]\ar@{-}[d] &*=0{}\ar@{-}[d]\\*=0{}\ar@{-}[r]&*=0{}}}$ of Young diagrams which begins with ${\cal C}_{\xymatrix@W=.35pc @H=.35pc @R=.35pc @C=.35pc {*=0{}\ar@{-}[r]\ar@{-}[d] &*=0{}\ar@{-}[d]\\*=0{}\ar@{-}[r]&*=0{}}}(1) =0$
and ${\cal C}_{\xymatrix@W=.35pc @H=.35pc @R=.35pc @C=.35pc {*=0{}\ar@{-}[r]\ar@{-}[d] &*=0{}\ar@{-}[d]\\*=0{}\ar@{-}[r]&*=0{}}}(2) = \xymatrix@W=.5pc @H=.5pc @R=0pc @C=0pc @*[F-]{~}~~$,
has terms ${\cal C}_{\xymatrix@W=.35pc @H=.35pc @R=.35pc @C=.35pc {*=0{}\ar@{-}[r]\ar@{-}[d] &*=0{}\ar@{-}[d]\\*=0{}\ar@{-}[r]&*=0{}}}(n)$ that are built of $n-1$ blocks each, and whose monotone decreasing sequence representation
is given by the formula
\[{\cal C}_{\xymatrix@W=.35pc @H=.35pc @R=.35pc @C=.35pc {*=0{}\ar@{-}[r]\ar@{-}[d] &*=0{}\ar@{-}[d]\\*=0{}\ar@{-}[r]&*=0{}}}(n)_k = \operatorname{Round}\left(n/2^k\right); k=1,2,\dots\]
where rounding is done to the nearest integer and .5 is rounded to zero.
\end{theorem}

\begin{proof}
The proof of the formula for the column heights is by way of first showing that
each term in
${\cal C}_{\xymatrix@W=.35pc @H=.35pc @R=.35pc @C=.35pc {*=0{}\ar@{-}[r]\ar@{-}[d] &*=0{}\ar@{-}[d]\\*=0{}\ar@{-}[r]&*=0{}}}$
can be built canonically as follows:
\[
{\cal C}_{\xymatrix@W=.35pc @H=.35pc @R=.35pc @C=.35pc {*=0{}\ar@{-}[r]\ar@{-}[d] &*=0{}\ar@{-}[d]\\*=0{}\ar@{-}[r]&*=0{}}}(n)
= \begin{array}{cc}  & \lceil\frac{n}{2}\rceil \\
{\cal C}_{\xymatrix@W=.35pc @H=.35pc @R=.35pc @C=.35pc {*=0{}\ar@{-}[r]\ar@{-}[d] &*=0{}\ar@{-}[d]\\*=0{}\ar@{-}[r]&*=0{}}}(\lceil\frac{n}{2}\rceil)\otimes_2(& \overbrace{\underbrace{{\cal C}_{\xymatrix@W=.35pc @H=.35pc @R=.35pc @C=.35pc {*=0{}\ar@{-}[r]\ar@{-}[d] &*=0{}\ar@{-}[d]\\*=0{}\ar@{-}[r]&*=0{}}}(2)\otimes_3 \dots \otimes_3{\cal C}_{\xymatrix@W=.35pc @H=.35pc @R=.35pc @C=.35pc {*=0{}\ar@{-}[r]\ar@{-}[d] &*=0{}\ar@{-}[d]\\*=0{}\ar@{-}[r]&*=0{}}}(2)} \otimes_3{\cal C}_{\xymatrix@W=.35pc @H=.35pc @R=.35pc @C=.35pc {*=0{}\ar@{-}[r]\ar@{-}[d] &*=0{}\ar@{-}[d]\\*=0{}\ar@{-}[r]&*=0{}}}(1)}~~) \\
&\lfloor\frac{n}{2}\rfloor \end{array}
\]
We must demonstrate that the
maximum of all ${\cal C}_{\xymatrix@W=.35pc @H=.35pc @R=.35pc @C=.35pc {*=0{}\ar@{-}[r]\ar@{-}[d] &*=0{}\ar@{-}[d]\\*=0{}\ar@{-}[r]&*=0{}}}(k)\otimes_2({\cal C}_{\xymatrix@W=.35pc @H=.35pc @R=.35pc @C=.35pc {*=0{}\ar@{-}[r]\ar@{-}[d] &*=0{}\ar@{-}[d]\\*=0{}\ar@{-}[r]&*=0{}}}(j_1)\otimes_3\dots\otimes_3{\cal C}_{\xymatrix@W=.35pc @H=.35pc @R=.35pc @C=.35pc {*=0{}\ar@{-}[r]\ar@{-}[d] &*=0{}\ar@{-}[d]\\*=0{}\ar@{-}[r]&*=0{}}}(j_k))$
where $\sum{j_i} = n$
is precisely given by the above canonical construction. We make the assumption (of strong induction)
that this holds for terms before the $n^{th}$ term, and check for the inequality
${\cal C}_{\xymatrix@W=.35pc @H=.35pc @R=.35pc @C=.35pc {*=0{}\ar@{-}[r]\ar@{-}[d] &*=0{}\ar@{-}[d]\\*=0{}\ar@{-}[r]&*=0{}}}(k)\otimes_2({\cal C}_{\xymatrix@W=.35pc @H=.35pc @R=.35pc @C=.35pc {*=0{}\ar@{-}[r]\ar@{-}[d] &*=0{}\ar@{-}[d]\\*=0{}\ar@{-}[r]&*=0{}}}(j_1)\otimes_3\dots\otimes_3{\cal C}_{\xymatrix@W=.35pc @H=.35pc @R=.35pc @C=.35pc {*=0{}\ar@{-}[r]\ar@{-}[d] &*=0{}\ar@{-}[d]\\*=0{}\ar@{-}[r]&*=0{}}}(j_k))$
less than or equal to the canonical construction.
The case in which there are only 0 or 1 odd integers among the $j_k$'s is
directly observed using the strong induction. If there are two or more
odd integers among the $j_k$'s  and the first column of the diagram they help determine
is greater than or equal to the first column of
${\cal C}_{\xymatrix@W=.35pc @H=.35pc @R=.35pc @C=.35pc {*=0{}\ar@{-}[r]\ar@{-}[d] &*=0{}\ar@{-}[d]\\*=0{}\ar@{-}[r]&*=0{}}}(k)$
then the inequality holds by induction on the size of the first column.
If there are two or more
odd integers among the $j_k$'s  and the first column of the diagram they help determine
is less than the first column of
${\cal C}_{\xymatrix@W=.35pc @H=.35pc @R=.35pc @C=.35pc {*=0{}\ar@{-}[r]\ar@{-}[d] &*=0{}\ar@{-}[d]\\*=0{}\ar@{-}[r]&*=0{}}}(k)$
then we check the sub-cases $n$ odd and $n$ even. For $n$ even the result is seen directly, and for
$n$ odd we again rely on induction.
\end{proof}

For comparison to the previous example of the operad with square terms,  the
instance of the associativity diagram with upper left position
${\cal C}_{\xymatrix@W=.35pc @H=.35pc @R=.35pc @C=.35pc {*=0{}\ar@{-}[r]\ar@{-}[d] &*=0{}\ar@{-}[d]\\*=0{}\ar@{-}[r]&*=0{}}}(2)\otimes_2({\cal C}_{\xymatrix@W=.35pc @H=.35pc @R=.35pc @C=.35pc {*=0{}\ar@{-}[r]\ar@{-}[d] &*=0{}\ar@{-}[d]\\*=0{}\ar@{-}[r]&*=0{}}}(3)\otimes_3{\cal C}_{\xymatrix@W=.35pc @H=.35pc @R=.35pc @C=.35pc {*=0{}\ar@{-}[r]\ar@{-}[d] &*=0{}\ar@{-}[d]\\*=0{}\ar@{-}[r]&*=0{}}}(2))\otimes_2
({\cal C}_{\xymatrix@W=.35pc @H=.35pc @R=.35pc @C=.35pc {*=0{}\ar@{-}[r]\ar@{-}[d] &*=0{}\ar@{-}[d]\\*=0{}\ar@{-}[r]&*=0{}}}(2)\otimes_3{\cal C}_{\xymatrix@W=.35pc @H=.35pc @R=.35pc @C=.35pc {*=0{}\ar@{-}[r]\ar@{-}[d] &*=0{}\ar@{-}[d]\\*=0{}\ar@{-}[r]&*=0{}}}(2)\otimes_3{\cal C}_{\xymatrix@W=.35pc @H=.35pc @R=.35pc @C=.35pc {*=0{}\ar@{-}[r]\ar@{-}[d] &*=0{}\ar@{-}[d]\\*=0{}\ar@{-}[r]&*=0{}}}(4)\otimes_3{\cal C}_{\xymatrix@W=.35pc @H=.35pc @R=.35pc @C=.35pc {*=0{}\ar@{-}[r]\ar@{-}[d] &*=0{}\ar@{-}[d]\\*=0{}\ar@{-}[r]&*=0{}}}(5)\otimes_3{\cal C}_{\xymatrix@W=.35pc @H=.35pc @R=.35pc @C=.35pc {*=0{}\ar@{-}[r]\ar@{-}[d] &*=0{}\ar@{-}[d]\\*=0{}\ar@{-}[r]&*=0{}}}(3))$ is as follows:

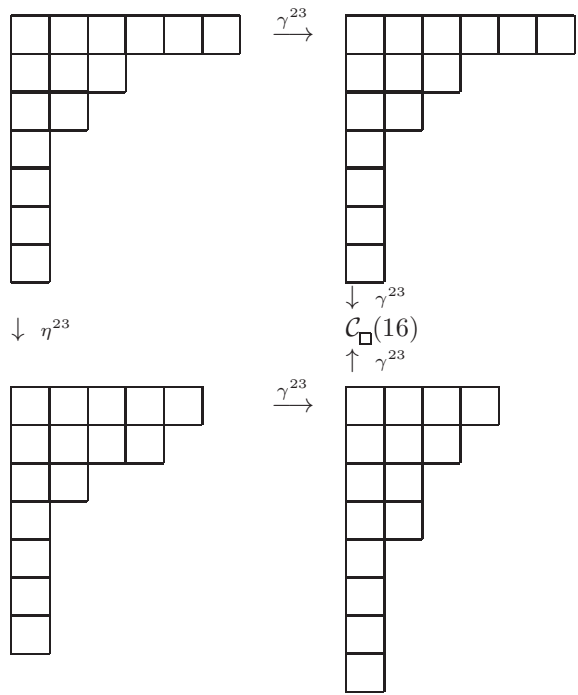
\begin{tabular}{llll}
&$\xymatrix@W=.7pc @H=.7pc @R=0pc @C=0pc @*[F-]{~&~&~&~&~&~\\~&~&~\\~&~\\~\\~\\~\\~}$&$\buildrel\gamma^{23}\over\longrightarrow$&$\xymatrix@W=.7pc @H=.7pc @R=0pc @C=0pc @*[F-]{~&~&~&~&~&~\\~&~&~\\~&~\\~\\~\\~\\~}$ \\
&&&$\downarrow~{\scriptstyle\gamma^{23}}$\\
&$\downarrow~{\scriptstyle\eta^{23}}$&&${\cal C}_{\xymatrix@W=.35pc @H=.35pc @R=.35pc @C=.35pc {*=0{}\ar@{-}[r]\ar@{-}[d] &*=0{}\ar@{-}[d]\\*=0{}\ar@{-}[r]&*=0{}}}(16)$\\
&&&$\uparrow~{\scriptstyle\gamma^{23}}$\\
&$\xymatrix@W=.7pc @H=.7pc @R=0pc @C=0pc @*[F-]{~&~&~&~&~\\~&~&~&~\\~&~\\~\\~\\~\\~}$&$\buildrel\gamma^{23}\over\longrightarrow$&$\xymatrix@W=.7pc @H=.7pc @R=0pc @C=0pc @*[F-]{~&~&~&~\\~&~&~\\~&~\\~&~\\~\\~\\~\\~}$\\
&\\
\end{tabular}

There may be interesting applications of the type of growth modeled
by operads in iterated monoidal categories.  Since the growth is in
multiple dimensions it suggests applications to studies of
allometric measurements. Broadly this refers to any $n$
characteristics of a system which grow in tandem. These measurements
are often used in biological sciences to try to predict values of
one characteristic from others, such as tree height from trunk
diameter or crown diameter, or skeletal mass from total body mass or
dimensions, or even genomic diversity from various geographical
features.  Allometric comparisons are often used in geology and
chemistry, for instance when predicting the growth of speleothems or
crystals. There are also potential applications to networks, where
the growth of diameter or linking distance of a network is related
logarithmically to the growth in number of nodes. In computational
geometry, the number of vertices of the convex hull of $n$ uniformly
scattered points in a polygon grows as the log of $n$.

This sort of minimal growth in the terms of the operad could be
perturbed, for example by replacing the term $\xymatrix@W=.5pc
@H=.5pc @R=0pc @C=0pc @*[F-]{~&~}\ $ in the above with the alternate
term \raisebox{0.5em} {$\xymatrix@W=.5pc @H=.5pc @R=0pc @C=0pc
@*[F-]{~\\~}\ $}, which would affect the later terms in turn. An
interesting avenue for further investigation would be the comparison
of such perturbations to determine the relative effects of a given
perturbation's size and position of occurrence in the sequence. What
we really want is a formula for minimal operads in Young diagrams
analogous to the one found above for operads in $\nat.$


We conclude with a description of the concepts of $n$-fold operad
algebra and of the tensor products of operad algebras.

\begin{definition}\label{opalg}
Let ${\cal C}$ be  an $n$-fold operad in ${\cal V}$. A
${\cal C}$-algebra is an object $A\in{\cal V}$ and maps
\[
\theta^{pq}:{\cal C}(j)\otimes_p(\otimes_q^j A)\to A
\]
for $n\ge q>p \ge 1$, $j\ge 0$.

\begin{enumerate}
\item Associativity: The following diagram is required to commute
for all $n\ge q>p \ge 1$, $k\ge 1$, $j_s\ge 0$ , where
$j = \sum\limits_{s=1}^k j_s$.

\[
\xymatrix{{\cal C}(k) \otimes_p ({\cal C}(j_1) \otimes_q \dots \otimes_q {\cal C}(j_k))\otimes_p(\otimes_q^j A)
\ar[rr]^>>>>>>>>>>>>{\gamma^{pq} \otimes_p \text{id}}
\ar[dd]_{\text{id} \otimes_p \eta^{pq}}
&& {\cal C}(j)\otimes_p(\otimes_q^j A)
\ar[d]_{\theta^{pq}}\\
&&A\\
{\cal C}(k) \otimes_p (({\cal C}(j_1)\otimes_p (\otimes_q^{j_1} A))\otimes_q \dots \otimes_q ({\cal C}(j_k)\otimes_p (\otimes_q^{j_k} A)))
\ar[rr]^>>>>>>>>>>>>{\text{id} \otimes_p (\otimes_q^k \theta^{pq})}
&&{\cal C}(k)\otimes_p(\otimes_q^k A)
\ar[u]^{\theta^{pq}}
}
\]

\item Units: The following diagram is required to commute for all
$n\ge q>p \ge 1$.
\[
\xymatrix{
I\otimes_p A
\ar[d]_{{\cal J}\otimes_p 1}
\ar@{=}[r]^{}
&A\\
{\cal C}(1)\otimes_p A
\ar[ur]^{\theta^{pq}}
}
\]
\end{enumerate}
\end{definition}

\begin{example}
Of course the initial object is always an algebra for every operad,
and every object is an algebra for the initial operad. For a slightly
less trivial example we turn to the height preordered category of
Remark~\ref{hpre}. Define the operad ${\cal B}(j)$ as in
Example~\ref{bee}. Then any nonzero sequence $A$ is an algebra for
this operad.
\end{example}

\begin{remark}
Let ${\cal C}$ and  ${\cal D}$ be $m$-fold operads in an $n$-fold
monoidal category.  If $A$ is an algebra of ${\cal C}$ and $B$ is
an algebra of ${\cal D}$ then $A\otimes_{i+m} B$ is an algebra for
${\cal C}\otimes'_i{\cal D}.$

That the product of $n$-fold operad algebras is itself an $n$-fold
operad algebra is easy to verify once we note that the new $\theta$
is in terms of the two old ones:
\[
\theta^{pq}_{A\otimes_{i+m} B} =
(\theta^{pq}_{A}\otimes_{i+m}\theta^{pq}_{B})\circ \eta^{p(i+m)}
    \circ (1 \otimes_p \eta^{q(i+m)})
\]
Maps of operad algebras are straightforward to describe--they are
required to preserve  structure; that is to commute with  $\theta.$
\end{remark}


\end{document}